\documentclass[12pt,authoryear]{elsarticle}
\journal{}

\makeatletter
\def\ps@pprintTitle{%
 \let\@oddhead\@empty
 \let\@evenhead\@empty
 \def\@oddfoot{\today\hfill\thepage}%
 \let\@evenfoot\@oddfoot}
\makeatother

\usepackage{amsmath,amssymb,amsfonts,amsthm,graphicx,appendix}
\usepackage[bookmarksnumbered,colorlinks=true]{hyperref}
\usepackage[labelfont=bf]{caption}

\providecommand{\doi}[1]{\href{https://doi.org/#1}{DOI:#1}}
\usepackage{xurl} 
\renewcommand{\doi}[1]{%
 \href{https://doi.org/#1}{\nolinkurl{DOI:#1}}%
}

\usepackage[OT1]{fontenc}
\usepackage{natbib,enumitem,mathtools,mathrsfs,dsfont}
\usepackage{tikz,float,orcidlink}
\usetikzlibrary{arrows.meta}
\usepackage{geometry}
\newtheorem{theorem}{Theorem}[section]
\newtheorem{lemma}[theorem]{Lemma}
\newtheorem{corollary}[theorem]{Corollary}
\newtheorem{proposition}[theorem]{Proposition}
\newtheorem{conjecture}[theorem]{Conjecture}

\theoremstyle{definition}

\newtheorem{remark}[theorem]{Remark}

\numberwithin{equation}{section}

\newcommand{\N}{\mathbb{N}}
\newcommand{\Z}{\mathbb{Z}}
\newcommand{\R}{\mathbb{R}}
\newcommand{\C}{\mathbb{C}}
\newcommand{\PP}{\mathsf{P}}
\newcommand{\EE}{\mathsf{E}}
\newcommand{\Var}{\mathsf{Var}}
\newcommand{\Cov}{\mathsf{Cov}}
\newcommand{\bb}[1]{\boldsymbol{#1}}
\newcommand{\ind}{\mathds{1}}
\newcommand{\rd}{\mathrm{d}}
\newcommand{\e}{\varepsilon}
\newcommand{\ii}{\mathrm{i}}
\renewcommand{\Re}{\mathrm{Re}\, }

\newcommand{\OO}{\mathcal{O}}
\newcommand{\oo}{\mathrm{o}}
\newcommand{\leqdef}{\vcentcolon=}

\newcommand{\Leb}{\mathrm{Leb}}
\newcommand{\sig}{\sigma}

\allowdisplaybreaks

\begin{document}

\begin{frontmatter}

\title{Convergence in law of the maximum of the Riemann zeta function on expanding intervals near the critical line}

\author{Fr\'{e}d\'{e}ric Ouimet~\texorpdfstring{\orcidlink{0000-0001-7933-5265}}{}}

\ead{frederic.ouimet2@uqtr.ca}
\address{Universit\'{e} du Qu\'{e}bec \`{a} Trois-Rivi\`{e}res, Canada.\vspace{-5mm}}

\begin{abstract}
Let $T \to \infty$, let $\tau$ be uniformly distributed on $[T,2T]$, and let $\sig = \sig(T)$ be a deterministic real part approaching $1/2$ from the right. We establish the convergence in law of the recentered maximum of $\log|\zeta(\sig + \ii\tau + \ii h)|$ over $I_{\theta} = \smash{[ - (\log T)^{\theta},(\log T)^{\theta}]}$, for every fixed $\theta > 0$. Under explicit conditions on the displacement $\sig - 1/2$, the recentered maximum converges to a Gumbel distribution with scale parameter $1/(2\sqrt{1 + \theta})$. The recentering is an exact upper quantile of the one-point law of the associated Haar field $\smash{F_{h}^{\mathrm{H}}} = \smash{\Re\sum_{p^{j}\leq X}U_{p}^{j}/(jp^{j(\sig + \ii h)})}$, where the $U_{p}$ are independent and uniformly distributed on the complex unit circle. The positive constant in the limiting distribution is characterized by a Gaussian process describing the limiting microscopic increments of the Haar field around a high point. The result is established without the Riemann hypothesis. The case $-1 < \theta \leq 0$ will be treated separately in an upcoming paper.
\end{abstract}


\end{frontmatter}

\setcounter{tocdepth}{2}
\tableofcontents

\section{Introduction}

\subsection{The logarithmic zeta field and its large values}

For a complex argument $s$ with $\Re(s) > 1$, the Riemann zeta function is defined by
\[
\zeta(s) \leqdef \sum_{n=1}^{\infty}n^{-s} = \prod_{p\ \mathrm{prime}}(1 - p^{-s})^{-1}.
\]
It extends meromorphically to $\C$, with a simple pole at $s = 1$. Let $T > 0$ be large. Throughout the paper, limits are taken as $T \to \infty$ unless otherwise specified. Let $\tau$ be uniformly distributed on $[T,2T]$ under the probability law $\PP$, with expectation $\EE$, and write $t \leqdef \log\log T$. For a deterministic real part $\sig = \sig(T)$, the random function $h \longmapsto \log|\zeta(\sig + \ii\tau + \ii h)|$ records the log-modulus of $\zeta$ along a randomly translated vertical segment. Its maximum depends on the frequency of large values at a single position and on the dependence between values at nearby positions. The present paper identifies a limiting distribution for this maximum on
\[
I_{\theta} \leqdef [ - (\log T)^{\theta},(\log T)^{\theta}], \qquad \theta > - 1,
\]
for fixed $\theta > 0$, when $\sig$ approaches $1/2$ from the right at the rates specified in Section~\ref{sec:main.results}.

At a single random height on the critical line, \citet{MR0020594,MR1220477} proved the central limit theorem
\[
\frac{\log|\zeta(1/2 + \ii\tau)|}{\sqrt{t/2}} \ \xrightarrow{\mathrm{d}}\ \mathcal{N}(0,1),
\]
where $\mathcal{N}(0,1)$ denotes the standard normal distribution. A short proof was provided by \citet{MR3832861}. The multidimensional results of \citet{MR2678896} described how the limiting Gaussian variables become correlated when their vertical shifts approach one another. These results identified typical fluctuations of order $\sqrt{t}$. Maxima over intervals of bounded length or of length growing as a power of $\log T$ instead occur at heights of order $t$, so their analysis requires information beyond the central limit scale.

The relevant large-deviation estimates included both the Gaussian exponential factor and the smaller prefactor. For every fixed $0 < \alpha < 2$, \citet{ArguinBailey2023} proved the unconditional upper bound
\[
\PP\big(\log|\zeta(1/2 + \ii\tau)| > \alpha t\big) \ll _{\alpha} t^{-1/2}\exp( - \alpha^{2}t).
\]
For every fixed $\alpha > 0$, \citet{ArguinBailey2025} obtained a lower bound of the same Gaussian order on the line $\sig = 1/2 + \delta/\log T$, with a fixed $\delta$ sufficiently large depending on $\alpha$. More recently, \citet{ArguinCreighton2026} proved Gaussian-order lower bounds for the Lebesgue measure of large values on the critical line itself. Assuming the Riemann hypothesis, \citet{ArguinBaileyRoberts2026} established the corresponding upper bounds for every fixed $\alpha > 0$. These upper and lower bounds did not identify the leading multiplicative factor in the one-point tail probability with a relative error tending to zero. Moreover, even exact one-point tail probabilities would not determine the distribution of the maximum over $I_{\theta}$: they describe how often the logarithmic zeta function exceeds a level at one random ordinate, whereas the maximum concerns whether it exceeds that level anywhere on the interval. Passing from the former probability to the latter requires additional information about the dependence between values at different ordinates.

The bounded terms in the recentering must also be determined to identify the location of the limiting maximum law. The moment conjectures of \citet{KeatingSnaith2000} contained an arithmetic Euler-product factor and a random-matrix factor. The arithmetic factor is $\mathfrak{a}(z)$ defined in \eqref{eq:maximum.arithmetic.factor} of Theorem~\ref{thm:maximum}. Proposition~\ref{prop:maximum.recentering} shows that this factor contributes to the recentering $u_{T}$ of $\max_{h\in I_{\theta}}\log|\zeta(\sig + \ii\tau + \ii h)|$, with $u_{T}$ defined precisely in \eqref{eq:height.quantile}.

\subsection{Maxima on the critical line and the length of the interval}

For comparison with the present result, define the critical-line maximum
\[
\mathcal{M}_{T}(\theta) \leqdef \max_{h\in I_{\theta}} \log|\zeta(1/2 + \ii\tau + \ii h)|.
\]
The predictions of \citet{FyodorovHiaryKeating2012} and \citet{MR3151088} related maxima over intervals of fixed length to extremes of logarithmically correlated fields and characteristic polynomials of random matrices. For $\theta = 0$, their predicted recentering was
\[
t - \frac{3}{4}\log t.
\]
The coefficient $3/4$ in the second term reflects the dependence between nearby high values. Replacing the critical-line logarithmic zeta field by independent Gaussian variables with the same variance and the same number of microscopic sampling positions would instead give the coefficient $1/4$, as explained by \citet{arXiv:1601.00582}.

The leading-order relation $\mathcal{M}_{T}(0)/t \to 1$ in probability was proved by \citet{MR3851835}, with the Riemann hypothesis required for the lower bound, and unconditionally by \citet{ABBRS_2019}. The partition-function estimates of \citet{arXiv:1906.05783} subsequently gave an upper bound with the predicted second-order term and a smaller third-order error. The sharper results of \citet{ABR-2020,ArguinBourgadeRadziwill2023} established
\begin{equation}\label{eq:intro.fixed.tightness}
\mathcal{M}_{T}(0) - t + \frac{3}{4}\log t = \OO_{\PP}(1),
\end{equation}
where $\OO_{\PP}(1)$ denotes a tight family of random variables. In particular, the upper-tail estimate of \citet{ABR-2020} took the form
\[
\PP\left(\mathcal{M}_{T}(0) > t - \frac{3}{4}\log t + y\right) \leq C y\exp( - 2y), \qquad y \geq 1,
\]
and \citet{ArguinBourgadeRadziwill2023} provided lower-tail control and a matching right-tail lower bound with a Gaussian correction in an extended range of $y$; the corresponding upper bound follows from the proof of \citet[Theorem~1]{ABR-2020}. Thus these results identified the correct recentering in \eqref{eq:intro.fixed.tightness} to bounded order. The more precise question of whether $\mathcal{M}_{T}(0) - (t - \frac{3}{4}\log t)$ converges in distribution remains open.

The behavior changes when the observation interval expands. For fixed $\theta > 0$, let
\[
m(\theta) \leqdef \sqrt{1 + \theta}.
\]
\citet{AOR2019} conjectured tightness of $\mathcal{M}_{T}(\theta)$ after subtracting
\begin{equation}\label{eq:intro.expanding.recentering}
m(\theta)t - \frac{1}{4m(\theta)}\log t.
\end{equation}
They proved the corresponding leading-order asymptotic
\[
\mathcal{M}_{T}(\theta) = m(\theta)t + \oo_{\PP}(t).
\]
Their lower bound was unconditional for every $\theta > 0$, while their matching upper bound was unconditional for $\theta \leq 3$ and assumed the Riemann hypothesis for $\theta > 3$. At the finer recentering in \eqref{eq:intro.expanding.recentering}, \citet{ArguinBailey2023} subsequently established upper-tail tightness for $0 < \theta < 3$. \citet{ArguinBaileyRoberts2026} extended this conclusion to every fixed $\theta > 0$ under the Riemann hypothesis.

A more precise conjecture was formulated and investigated numerically by \citet{arXiv:2104.07403} on the critical line ($\sig = 1/2$). For an interval with half-length $\pi(\log T)^{\theta}$, they predicted convergence to a Gumbel distribution after the recentering in \eqref{eq:intro.expanding.recentering}, with scale parameter $1/(2m(\theta))$ and a specified bounded location correction. The interval factor $\pi$ and the one-point arithmetic corrections affect the bounded correction to the recentering. Moving away from the critical line can also change the coefficient of $\log t$, as shown in \eqref{eq:maximum.displacement.recentering} below.

The coefficient $1/(4m(\theta))$ in \eqref{eq:intro.expanding.recentering} can be motivated by a simple heuristic calculation. The interval $I_{\theta}$ contains approximately $e^{(1 + \theta) t}$ positions separated by the critical-line microscopic scale $(\log T)^{-1}$. For a centered Gaussian variable $Z$ of variance $t/2$, the probability of exceeding a threshold $v \asymp t$ satisfies $\PP(Z > v) \asymp t^{-1/2} e^{- v^{2}/t}$. Requiring the expected number of exceedances to be nontrivial gives
\begin{equation}\label{eq:expected.exceedances}
e^{(1 + \theta) t} \times t^{-1/2} e^{- v^{2}/t} \asymp 1
\qquad\Longrightarrow\qquad
v = m(\theta)t - \frac{\log t}{4m(\theta)} + \OO(1),
\end{equation}
which recovers \eqref{eq:intro.expanding.recentering} to bounded order. This calculation explains the proposed normalization, but it does not justify treating the positions as independent or determine a limiting law. Correlations between nearby high values still affect the constant that determines the location of the predicted limiting Gumbel distribution.

Shrinking intervals ($- 1 < \theta < 0$) exhibit a further contribution from primes whose phases change little across the whole interval; see, e.g., Fig.~1 of \citet{AOR2019} and the Gaussian term in Conjecture~\ref{conj:maximum.mesoscopic} below. In this setting, \citet{ArguinHamdan2026} obtained upper bounds to bounded-order precision after removing an explicit small-prime contribution with Gaussian fluctuations, complementing the prior works of \citet{ABR-2020,ArguinBourgadeRadziwill2023} for $\theta = 0$. The present paper treats the case $\theta > 0$.

\subsection{Randomized prime sums and Gaussian approximations}

Randomized Euler products separate the prime contributions while preserving the spatial dependence generated by the factors $p^{-\ii h}$. Let $(U_{p})_{p}$ and $(G_{p})_{p}$ be two families indexed by primes. The variables $U_{p}$ are independent and uniformly distributed on the complex unit circle, and the variables $G_{p}$ are independent centered complex Gaussian variables whose real and imaginary parts are independent with variance $1/2$. Two frequently studied prime-sum models are
\begin{equation}\label{eq:toy.models}
F_{h}^{\mathrm{rand}} \leqdef \Re\sum_{p\leq T}\frac{U_{p}p^{-\ii h}}{\sqrt{p}}, \qquad F_{h}^{\mathrm{G}} \leqdef \Re\sum_{p\leq T}\frac{G_{p}p^{-\ii h}}{\sqrt{p}},
\end{equation}
where both sums run over primes $p \leq T$. Their covariance functions agree:
\[
\Cov(F_{h}^{\mathrm{rand}},F_{h'}^{\mathrm{rand}})
= \Cov(F_{h}^{\mathrm{G}},F_{h'}^{\mathrm{G}})
= \frac{1}{2}\sum_{p\leq T}\frac{\cos((h - h')\log p)}{p}.
\]
For $h$ and $h'$ in a fixed bounded interval, this sum grows logarithmically as $|h - h'|$ decreases, until the variance scale $t/2 + \OO(1)$ is reached. Organizing the primes according to the size of $\log p$ produces an approximate branching structure, as illustrated in Fig.~1 of \cite{AOR2019}: nearby positions share contributions from low frequencies and separate at higher frequencies.

For the model with uniform prime phases on $I_{0} = [-1,1]$, \citet{MR3619786} proved the two-term asymptotic
\[
\max_{h\in I_{0}}F_{h}^{\mathrm{rand}} = t - \frac{3}{4}\log t + \oo_{\PP}(\log t).
\]
For the Gaussian prime model on expanding intervals ($\theta > 0$), \citet{ADH_2021} established
\[
\max_{h\in I_{\theta}}F_{h}^{\mathrm{G}} = m(\theta)t - \frac{1}{4m(\theta)}\log t + \oo_{\PP}(\log t).
\]
They also analyzed the intermediate regime $\theta = t^{-\alpha}$ with $0 < \alpha < 1$, in which the coefficient of the logarithmic correction interpolated between the expanding-interval and fixed-interval values. In the same regime, \citet{Chang2024} subsequently established tightness and refined tail estimates after the recentering $m(\theta)t - \{(1 + 2\alpha)/(4 m(\theta))\}\log t$.

\subsection{Limiting laws for maxima in branching, Gaussian, and random-matrix models}

Several related models have limiting maximum laws that go beyond two-term asymptotics or tightness. Under suitable nonlattice and integrability assumptions, \citet{Aidekon2013} proved convergence in law for the minimum of a branching random walk in its boundary normalization. Conditional on nonextinction, the limiting law was expressed through a positive derivative-martingale limit. For logarithmically correlated Gaussian fields, \citet{BramsonDingZeitouni2016} established convergence of the recentered maximum of the two-dimensional discrete Gaussian free field, while \citet{Madaule2015} treated continuous Gaussian fields on $[0,1]^{d}$ with covariance obtained by integrating dilations of a compactly supported $C^{1}$ kernel. \citet{DingRoyZeitouni2017} proved a general convergence theorem for discrete log-correlated Gaussian fields satisfying precise covariance and convergence assumptions.

For maxima, or negatives of minima, these results yielded limiting distribution functions of the form
\begin{equation}\label{eq:intro.random.shift}
\EE\big[\exp\{ - C Z\exp( - \gamma x)\}\big], \qquad x \in \R, \quad C, \gamma > 0,
\end{equation}
where $Z$ is a positive random variable arising from a derivative martingale or an analogue of such a martingale. The parameter $\gamma$ depends on the reproduction and displacement laws in branching models, and on the dimension and covariance normalization in Gaussian fields.

After deterministic centering and scaling, the distribution functions of maxima of certain correlated Gaussian fields and branching models also converge to the simpler Gumbel form
\begin{equation}\label{eq:intro.gumbel}
\exp\{ - C\exp( - \gamma x)\}, \qquad x \in \R, \quad C, \gamma > 0,
\end{equation}
with no random derivative-martingale factor. \citet{ChiariniCiprianiHazra2016} proved such a limit for the discrete Gaussian free field in dimensions $d \geq 3$. For fixed $0 \leq \alpha < 1$, \citet{KistlerSchmidt2015} obtained the simpler Gumbel limit after recentering for Gaussian fields with $N^{\alpha}$ hierarchical levels, interpolating between the random energy model ($\alpha = 0$) and the branching random walk ($\alpha = 1$). Further examples include massive Gaussian free fields, certain discrete fractional Gaussian fields, and membrane models in dimensions $d \geq 5$, treated by \citet{ChiariniCiprianiHazraInterfaces2016}, as well as Gaussian branching random walks conditioned to have all terminal values nonnegative, treated by \citet{FelsHartungLouidor2024}.

Limiting laws for recentered maxima have also been established for the log-modulus of characteristic polynomials of random matrices. Let $U_{N}$ be an $N \times N$ random matrix from the circular unitary ensemble (CUE), that is, a Haar-distributed unitary matrix, and define
\[
F_{h}^{\mathrm{CUE}} \leqdef \log|\det(\mathrm{Id}_{N} - U_{N}e^{-\ii h})|, \qquad h \in [0,2\pi].
\]
The moment calculations of \citet{KeatingSnaith2000} helped establish the comparison with $\log|\zeta|$, and \citet{FyodorovHiaryKeating2012} and \citet{MR3151088} predicted convergence in law of the maximum after recentering by $\log N - (3/4)\log\log N$. The leading term was proved by \citet{MR3594368}, and the second-order term by \citet{MR3848227}. In the broader C$\beta$E setting (the CUE case is $\beta = 2$), \citet{MR3848391} proved tightness with recentering
\[
\sqrt{2/\beta}(\log N - (3/4)\log\log N),
\]
for every $\beta > 0$. Subsequently, \citet{PaquetteZeitouni2025} proved convergence of the distribution function of the recentered maximum to an expression of the form \eqref{eq:intro.random.shift}, with $\gamma = \sqrt{2\beta}$ and $Z$ the positive total mass of a limiting derivative-martingale measure.

The main result of the present paper, Theorem~\ref{thm:maximum}, establishes a limiting distribution function of the simpler form \eqref{eq:intro.gumbel} for the recentered maxima of $\log|\zeta|$ and of the Euler field $F^{\mathrm{A}}$ in \eqref{eq:Haar.field} on expanding intervals ($\theta > 0$) near the critical line. This is the first convergence-in-law result of this type for the local maxima of $\log|\zeta|$ or of a Dirichlet-polynomial approximation to this field.

\newpage
For bounded or shrinking intervals on the critical line, we formulate the conjecture of \citet[Section~2]{MR3151088} below in terms of joint convergence with the small-prime contribution. This formulation is motivated by \citet[Eq.~(1.15) and the subsequent discussion]{AOR2019} and \citet[Conjecture~1.1 and Theorem~1.2]{ArguinHamdan2026}.

\begin{conjecture}[{\citet{MR3151088}}]\label{conj:maximum.mesoscopic}
Let $ - 1 < \theta \leq 0$ and set $t_{\theta} \leqdef (1 + \theta)t$ and $G_T(\theta) = \ind_{\{\theta<0\}}\Re\sum_{\log p\leq\exp(|\theta|t)} p^{-1/2 - \ii\tau}$. Then, as $T \to \infty$,
\[
\left(\frac{G_T(\theta)}{\sqrt{t}}, \max_{h\in I_{\theta}}\log|\zeta(1/2 + \ii\tau + \ii h)| - G_T(\theta) - t_{\theta} + \frac{3}{4}\log(t_{\theta})\right)
\xrightarrow{\mathrm{d}} \left(\sqrt{\frac{|\theta|}{2}}\,Z_{0}, \mathcal{M}_{\theta}\right),
\]
where $Z_{0} \sim \mathcal{N}(0,1)$ is independent of $\mathcal{M}_{\theta}$, whose distribution function has the mixed Gumbel form \eqref{eq:intro.random.shift}, with $\gamma = 2$, and whose upper tail satisfies
\[
\PP(\mathcal{M}_{\theta} > y) \sim C_{\theta}y\exp( - 2y), \qquad y \to \infty,
\]
for some constant $C_{\theta} > 0$. The Gaussian contribution vanishes when $\theta = 0$.
\end{conjecture}

\subsection{Organization}

Section~\ref{sec:main.results} gives the parameter assumptions, defines the Euler and Haar fields, introduces the quantile recentering, and states the main theorem. Section~\ref{sec:high.points} proves the joint conditional local-profile limit and localizes high excursions under simultaneous conditioning. Section~\ref{sec:void} establishes the limiting probabilities of no exceedance on random grids. Section~\ref{sec:continuum} passes to the continuum and establishes the limiting constant. Section~\ref{sec:recentering} evaluates the recentering explicitly and completes the proof of the main theorem, including the transfer to the zeta function.

\section{Parameters, notation, and the main result}\label{sec:main.results}

\subsection{Parameter assumptions}

Let $t \leqdef \log\log T$, and fix constants $\theta > 0$, $\phi > 7/2$, $\Phi > \phi + 1$, and $C_{\sig} > 0$. The real part $\sig = \sig(T)$ and the prime-power cutoff $X = X(T)$ are deterministic. Write $\eta_{\sig} \leqdef \sig - 1/2$, and assume throughout the paper that
\begin{equation}\label{eq:fixed.parameters}
T^{1/t^{\Phi}} \leq X \leq T^{1/t^{\phi}}, \qquad 0 < \eta_{\sig} \leq \frac{C_{\sig}t^{\Phi}}{\log T}, \qquad \frac{\eta_{\sig}\log X}{t} \to \infty.
\end{equation}
For example, $X = T^{1/t^{4}}$ and $\eta_{\sig} = \delta t^{6}/\log T$, with a fixed $\delta > 0$, satisfy \eqref{eq:fixed.parameters} with $\phi = 4$, $\Phi = 6$, and $C_{\sig} \geq \delta$.

Define the microscopic scale, effective logarithmic depth, and effective number of microscopic intervals by
\begin{equation}\label{eq:scales}
\ell_{T} \leqdef 2\eta_{\sig}, \qquad s_{T} \leqdef \log(1/\ell_{T}), \qquad N_{T} \leqdef \frac{|I_{\theta}|}{\ell_{T}}.
\end{equation}
The assumptions in \eqref{eq:fixed.parameters} imply
\begin{equation}\label{eq:scale.asymptotics}
s_{T} = t + \OO(\log t), \qquad \log N_{T} = (1 + \theta)t + \OO(\log t), \qquad \ell_{T}\log X \to \infty.
\end{equation}
Indeed, $\eta_{\sig}\log X \geq t$ eventually, and the upper bound on $X$ then gives $\eta_{\sig} \geq t^{\phi + 1}/\log T$. Combining this with the upper bound on $\eta_{\sig}$ proves the first assertion of \eqref{eq:scale.asymptotics}. The other assertions follow from \eqref{eq:scales} and \eqref{eq:fixed.parameters}.

\subsection{Fields, maxima, and recentering}

Probabilities and expectations with respect to $\tau$ are denoted by $\PP$ and $\EE$, respectively. To distinguish this arithmetic randomness from other sources, such as the prime phases $(U_{p})_{p}$ defining the Haar field in \eqref{eq:Haar.field}, the notation $\PP_{\mathrm{A}} \leqdef \PP$ and $\EE_{\mathrm{A}} \leqdef \EE$ is also used, where $\mathrm{A}$ stands for arithmetic. Define the Euler polynomial by
\[
\mathcal{Q}(s) \leqdef \sum_{p^{j}\leq X}\frac{p^{-js}}{j}, \qquad s \in \C,
\]
where the sum runs over primes $p$ and positive integers $j$ satisfying $p^{j} \leq X$. With the same prime-indexed sequence of uniforms $(U_{p})_{p}$ as in \eqref{eq:toy.models}, define the Euler field (also called the arithmetic field) and the Haar field by
\begin{equation}\label{eq:Haar.field}
F_{h}^{\mathrm{A}} \leqdef \Re\mathcal{Q}(\sig + \ii\tau + \ii h), \qquad
F_{h}^{\mathrm{H}} \leqdef \Re\sum_{p^{j}\leq X}\frac{U_{p}^{j} p^{-\ii jh}}{j p^{j\sig}}, \qquad h \in I_{\theta}.
\end{equation}
Probabilities and expectations with respect to the prime phases $(U_{p})_{p}$ are denoted by $\PP_{\mathrm{H}}$ and $\EE_{\mathrm{H}}$, respectively. Independence and rotational invariance of these phases imply that the distribution of $F_{h}^{\mathrm{H}}$ does not depend on $h$. In particular, the pointwise variance of the Haar field is
\begin{equation}\label{eq:field.variance}
V_{T} \leqdef \EE_{\mathrm{H}}[(F_{0}^{\mathrm{H}})^{2}] = \frac{1}{2}\sum_{p^{j}\leq X}\frac{1}{j^{2}p^{2j\sig}} = \frac{s_{T}}{2} + \OO(1),
\end{equation}
which follows from the prime number theorem; see, e.g., \citet[Eqs.~(4.12) and~(4.14)]{OuimetHaar2026} or \citet[proof of Lemma~4.3, Eq.~(4.17)]{OuimetExpanding2026}.

For $\star \in \{\mathrm{A}, \mathrm{H}\}$, write
\[
M_{T}^{\star} \leqdef \max_{h\in I_{\theta}}F_{h}^{\star}, \qquad M_{T}^{\zeta} \leqdef \max_{h\in I_{\theta}}\log|\zeta(\sig + \ii\tau + \ii h)|.
\]
At a zero of $\zeta$ on the segment $\{\sig + \ii\tau + \ii h : h \in I_{\theta}\}$, let $\log|\zeta(\sig + \ii\tau + \ii h)| = - \infty$. This convention has no effect on $\smash{M_{T}^{\zeta}}$, which can equivalently be defined as the logarithm of the maximum of the continuous function $h \mapsto |\zeta(\sig + \ii\tau + \ii h)|$.

Let $u_{T}$ be the unique upper quantile determined, for all sufficiently large $T$, by
\begin{equation}\label{eq:height.quantile}
N_{T}\PP_{\mathrm{H}}(F_{0}^{\mathrm{H}} > u_{T}) = 1.
\end{equation}
This normalization is the analogue of the heuristic balance in \eqref{eq:expected.exceedances}, with $N_{T}$ microscopic intervals and the Haar tail probability in place of its Gaussian approximation. Theorem~\ref{thm:maximum} shows that $u_{T}$ provides the recentering for the Haar maximum $M_{T}^{\mathrm{H}}$, the Euler maximum $M_{T}^{\mathrm{A}}$, and the logarithmic zeta maximum $\smash{M_{T}^{\zeta}}$. Applying the one-point estimates proved in \citet[Lemma~4.7, Eq.~(4.79)]{OuimetExpanding2026} gives, as $T \to \infty$,
\begin{equation}\label{eq:height.asymptotics}
\frac{u_{T}}{s_{T}} \to \sqrt{1 + \theta}, \qquad \lambda_{T} \leqdef \frac{u_{T}}{V_{T}} \to \beta_{c}(\theta) \leqdef 2m(\theta) = 2\sqrt{1 + \theta}.
\end{equation}

\subsection{The main result}

Let $B_{1}$ and $B_{2}$ be independent standard Brownian motions, and define
\[
G(v) \leqdef \frac{1}{\sqrt{2}}\int_{0}^{\infty} \frac{e^{-y/2}}{\sqrt{y}} \left\{(\cos(vy) - 1)\rd{B}_{1}(y) + \sin(vy)\rd{B}_{2}(y)\right\}, \qquad v \in \R.
\]
The integral exists in $L^{2}$ because its integrands are square integrable. It defines a centered Gaussian process anchored by $G(0) = 0$, with covariance
\begin{equation}\label{eq:main.profile.covariance}
\mathfrak{C}(v,w) \leqdef \EE[G(v)G(w)] = \frac{1}{4}\log\left( \frac{(1 + v^{2})(1 + w^{2})}{1 + (v - w)^{2}} \right), \qquad v,w \in \R.
\end{equation}
By \eqref{eq:main.profile.covariance} and the Gaussianity of $G$, $\EE[|G(v) - G(w)|^{4}] = 3 \, \EE[|G(v) - G(w)|^{2}]^{2} \leq C|v - w|^{4}$, so Kolmogorov's continuity theorem \citep[Theorem~4.23, p.~95]{Kallenberg2021} gives a continuous version, which is used throughout.

Define
\begin{equation}\label{eq:main.profile}
D_{\theta}(v) \leqdef G(v) - \frac{\beta_{c}(\theta)}{4}\log(1 + v^{2}), \qquad v \in \R.
\end{equation}
As proved in Lemma~\ref{lem:local.profile}, after conditioning $F_{h}^{\star} - u_{T}$ to lie in a fixed bounded interval of positive length, the process $v \mapsto \smash{F_{h + \ell_{T}v}^{\star} - F_{h}^{\star}}$ converges in distribution to $D_{\theta}$ in $C([ - R,R])$ for every fixed $R > 0$. Thus, $D_{\theta}$ describes the limiting local fluctuations of either the Euler field $F^{\mathrm{A}}$ or the Haar field $F^{\mathrm{H}}$ relative to its height at the conditioning point.

For $q > 0$, define
\begin{equation}\label{eq:lattice.constant}
p(q) \leqdef \PP\big(D_{\theta}(qj) \leq 0\ \text{for every }j \in \Z\setminus\{0\}\big), \qquad
\kappa_{\theta}(q) \leqdef \frac{p(q)}{q}, \qquad
\kappa_{\theta} \leqdef \lim_{q\downarrow0}\kappa_{\theta}(q).
\end{equation}
Since $D_{\theta}(0) = 0$, the probability $p(q)$ measures how often the origin is the highest point of the limiting Gaussian profile $D_{\theta}$ sampled on the lattice $q\Z$. The normalization by $q$ in $\kappa_{\theta}(q)$ accounts for the lattice spacing.

\begin{theorem}[Limiting distribution of the maximum]\label{thm:maximum}
Assume \eqref{eq:fixed.parameters} and fix $\theta > 0$. The limit as $q \downarrow 0$ in \eqref{eq:lattice.constant} exists and satisfies $0 < \kappa_{\theta} < \infty$. For every $x \in \R$, as $T \to \infty$,
\begin{equation}\label{eq:maximum.limit}
\PP(M_{T}^{\zeta} - u_{T} \leq x) \to \exp\big\{ - \kappa_{\theta}\exp( - \beta_{c}(\theta)x)\big\}.
\end{equation}
Thus the limiting law is Gumbel with location parameter $\log\kappa_{\theta}/\beta_{c}(\theta)$ and scale parameter $1/\beta_{c}(\theta) = 1/(2\sqrt{1 + \theta})$. The same limit holds as $T \to \infty$ for $\PP_{\star}(M_{T}^{\star} - u_{T} \leq x)$, for $\star \in \{\mathrm{A}, \mathrm{H}\}$. For a displacement $\eta_{\sig} = c\,t^{A}/\log T$ satisfying \eqref{eq:fixed.parameters}, with fixed $c > 0$ and $\phi + 1 < A \leq \Phi$ (and $c \leq C_{\sig}$ if $A = \Phi$), the quantile satisfies, as $T \to \infty$,
\begin{equation}\label{eq:maximum.displacement.recentering}
u_{T} = m(\theta) t - \frac{2A(2 + \theta) + 1}{4m(\theta)}\log t - \frac{2 + \theta}{2m(\theta)}\log(2c) + \frac{1}{2m(\theta)}\log\left(\frac{\mathfrak{a}(m(\theta))}{m(\theta)\sqrt{\pi}}\right) + \oo(1),
\end{equation}
where $\mathfrak{a}$ is the arithmetic factor appearing in the moment conjecture of \citet[Eq.~(5)]{KeatingSnaith2000}, defined for $z > 0$ by
\begin{equation}\label{eq:maximum.arithmetic.factor}
\mathfrak{a}(z) \leqdef \prod_{p} \left(1 - \frac{1}{p}\right)^{z^{2}} \sum_{j=0}^{\infty}\frac{(z^{(j)})^{2}}{(j!)^{2}p^{j}}, \qquad z^{(j)} \leqdef \frac{\Gamma(z + j)}{\Gamma(z)} = \prod_{r=0}^{j - 1}(z + r), \quad j \in \N,
\end{equation}
and $z^{(0)} \leqdef 1$ by convention.
\end{theorem}

\begin{remark}
The coefficient of $\log t$ in \eqref{eq:maximum.displacement.recentering} differs from that in the critical-line prediction \eqref{eq:intro.expanding.recentering}: the displacement contributes the additional term $ - A(2 + \theta)\log t/(2m(\theta))$.
\end{remark}

\subsection{Outline and ideas}

The main task is to turn the estimates of \citet{OuimetExpanding2026} for exceedance probabilities near $u_{T}$ at one or several positions into a limiting law for $M_{T}^{\zeta}$. This requires counting each cluster of nearby exceedances once and controlling the dependence between different clusters.

{\bf Step 1. Locate the relevant heights.} One can picture $I_{\theta}$ as a forest, as illustrated in Fig.~1 of \citet{AOR2019}: each spatial interval of length one, called a unit block, is a tree. Its leaves are microscopic positions separated by distances of order $\ell_{T} = 2\eta_{\sig}$, while shared branches represent prime contributions that vary little between nearby positions. There are therefore about $N_{T} = |I_{\theta}|/\ell_{T}$ such positions across the forest.

The quantile equation $N_{T}\PP_{\mathrm{H}}(F_{0}^{\mathrm{H}} > u_{T}) = 1$ balances this number of positions with the probability of an exceedance at one position. The high-value tail has the Gaussian leading term $\log\PP_{\mathrm{H}}(F_{0}^{\mathrm{H}} > v) = - v^{2}/s_{T} + \oo(s_{T})$ for $v \asymp s_{T}$. The defining quantile equation therefore gives $u_{T}^{2}/s_{T} \sim \log N_{T}$, hence $u_{T} \sim \sqrt{s_{T}\log N_{T}}$. Moreover, the one-point estimate of \citet[Lemma~4.7]{OuimetExpanding2026} gives, for every fixed $x \in \R$, uniformly in $h \in I_{\theta}$, as $T \to \infty$,
\[
N_{T}\PP_{\star}(F_{h}^{\star} > u_{T} + x) \to e^{-\beta_{c}(\theta)x}, \qquad \beta_{c}(\theta) = 2\sqrt{1 + \theta}.
\]
This identifies the exponential dependence on $x$, but counting every sampled position where $F_{h}^{\star} > u_{T} + x$ would count the same cluster several times.

{\bf Step 2. Describe and localize a high peak.} To understand this clustering, Section~\ref{sec:high.points} conditions on $F_{h}^{\star} \in [u_{T},u_{T} + 1]$ at a fixed position $h$ and identifies the resulting local profile. The variable $v$ describes the displacement from $h$ in units of $\ell_{T}$:
\[
F_{h + \ell_{T}v}^{\star} - F_{h}^{\star} \mid \{F_{h}^{\star} - u_{T} \in [0,1]\} \ \xrightarrow{\mathrm{d}}\ D_{\theta}(v) = G(v) - \frac{\beta_{c}(\theta)}{4}\log(1 + v^{2}).
\]
To use this local description for counting peaks, we must also show that comparable heights are unlikely farther away within the same unit block, even under simultaneous conditioning at several separated points. For $r > 0$ with $d = \ell_{T}e^{r} \leq 1$, consider positions $h'$ in the same block satisfying $|h' - h| \asymp d$. Let $\smash{F_{d}^{\star,\mathrm{lo}}}(h)$ denote the contribution at $h$ of the prime-power terms with $\log(p^{j}) \leq 1/d$, and subtract this same contribution from both field values:
\[
X^{\star} = F_{h}^{\star} - F_{d}^{\star,\mathrm{lo}}(h), \qquad Y^{\star}(h') = F_{h'}^{\star} - F_{d}^{\star,\mathrm{lo}}(h).
\]
Heuristically, under the high-point conditioning, $X^{\star}$ has mean $m(\theta)r$ and variance at most about $r/2$, while the values $Y^{\star}(h')$ at the roughly $e^{r}$ competing microscopic positions are approximately centered with variance about $r/2$. Fix $C > 0$. Since $\theta > 0$ implies $m(\theta) > 1$, choose $1 < a < m(\theta)$. Since $F_{h'}^{\star} - F_{h}^{\star} = Y^{\star}(h') - X^{\star}$, a competitor can reach $F_{h}^{\star} - C$ only if $X^{\star} \leq ar + C$ or $Y^{\star}(h') \geq ar$. Gaussian tail estimates and a union bound suggest the total conditional probability bound
\[
e^{-(m(\theta) - a)^{2}r} + e^{r}e^{-a^{2}r} = e^{-(m(\theta) - a)^{2}r} + e^{-(a^{2} - 1)r} \to 0.
\]
Both events therefore become unlikely as the relative distance $d/\ell_{T} = e^{r}$ increases.

{\bf Step 3. Count each peak once.} The local profile and localization estimate now allow us to count high blocks through their grid maximizers. For this purpose, Section~\ref{sec:void} uses randomly shifted grids of spacing $q\ell_{T}$ in trimmed unit blocks separated by small gaps. These gaps keep points in distinct blocks a fixed positive distance apart, allowing the joint estimates for separated high points to apply. A high grid point represents its block when no other grid value there is larger. Away from block endpoints, high-point conditioning gives the limiting probability $p(q) = \PP(D_{\theta}(qj) \leq 0\text{ for all }j \in \Z\setminus\{0\})$. As the gaps vanish, the expected number of high blocks is approximately
\[
\frac{N_{T}}{q}\times\frac{e^{-\beta_{c}(\theta)x}}{N_{T}}\times p(q) = \kappa_{\theta}(q)e^{-\beta_{c}(\theta)x}.
\]
To obtain a Poisson law from this expected count, we must also control joint occurrences of high blocks. The random shifts turn sums over possible maximizers into spatial integrals. Across the expanding forest, averaging over separated positions makes joint high-value probabilities behave like products of one-point probabilities. Counting distinct high blocks in bounded height bands then gives, as $T \to \infty$ and then the gaps vanish, the corresponding Poisson factorial moments. Removing the upper height cutoff yields the zero-count probability $\exp\{ - \kappa_{\theta}(q)e^{-\beta_{c}(\theta)x}\}$, which is the limiting Gumbel distribution function of the recentered grid maximum: no high block means that every grid value is at most $u_{T} + x$.

{\bf Step 4. Recover the continuum and zeta maxima.} To pass from this grid law to the continuum, Section~\ref{sec:continuum} shows that grids with spacing $q\ell_{T}$ capture the relevant peaks as $q \downarrow 0$, and that the omitted gaps have negligible effect. Together with the tightness of the recentered Euler and Haar maxima established in Lemma~\ref{max:quantile.tightness}, this proves $\kappa_{\theta}(q) \to \kappa_{\theta} \in (0,\infty)$. The approximation $|\smash{M_{T}^{\zeta}} - M_{T}^{\mathrm{A}}| = \oo_{\PP}(1)$ from \citet{OuimetExpanding2026} transfers the Gumbel law to zeta.

With the limiting law established, Section~\ref{sec:recentering} makes the quantile recentering explicit by evaluating the Haar cumulant generating function: for fixed $\lambda > 0$, $\log\EE_{\mathrm{H}}[e^{\lambda F_{0}^{\mathrm{H}}}] = \lambda^{2}s_{T}/4 + \log\mathfrak{a}(\lambda/2) + \oo(1)$. Combining this expansion with the one-point tail estimate via an exponential tilting argument, we are able to deduce
\[
1 = N_{T}\PP_{\mathrm{H}}(F_{0}^{\mathrm{H}} > u_{T}) \sim \frac{N_{T}\sqrt{s_{T}}}{2u_{T}\sqrt{\pi}}\mathfrak{a}(u_{T}/s_{T})e^{-u_{T}^{2}/s_{T}}.
\]
Solving this quantile equation identifies the arithmetic correction to $u_{T}$, while $\kappa_{\theta}$ accounts for the clustering of high values. Together, these corrections determine the location of the limiting law beyond tightness.

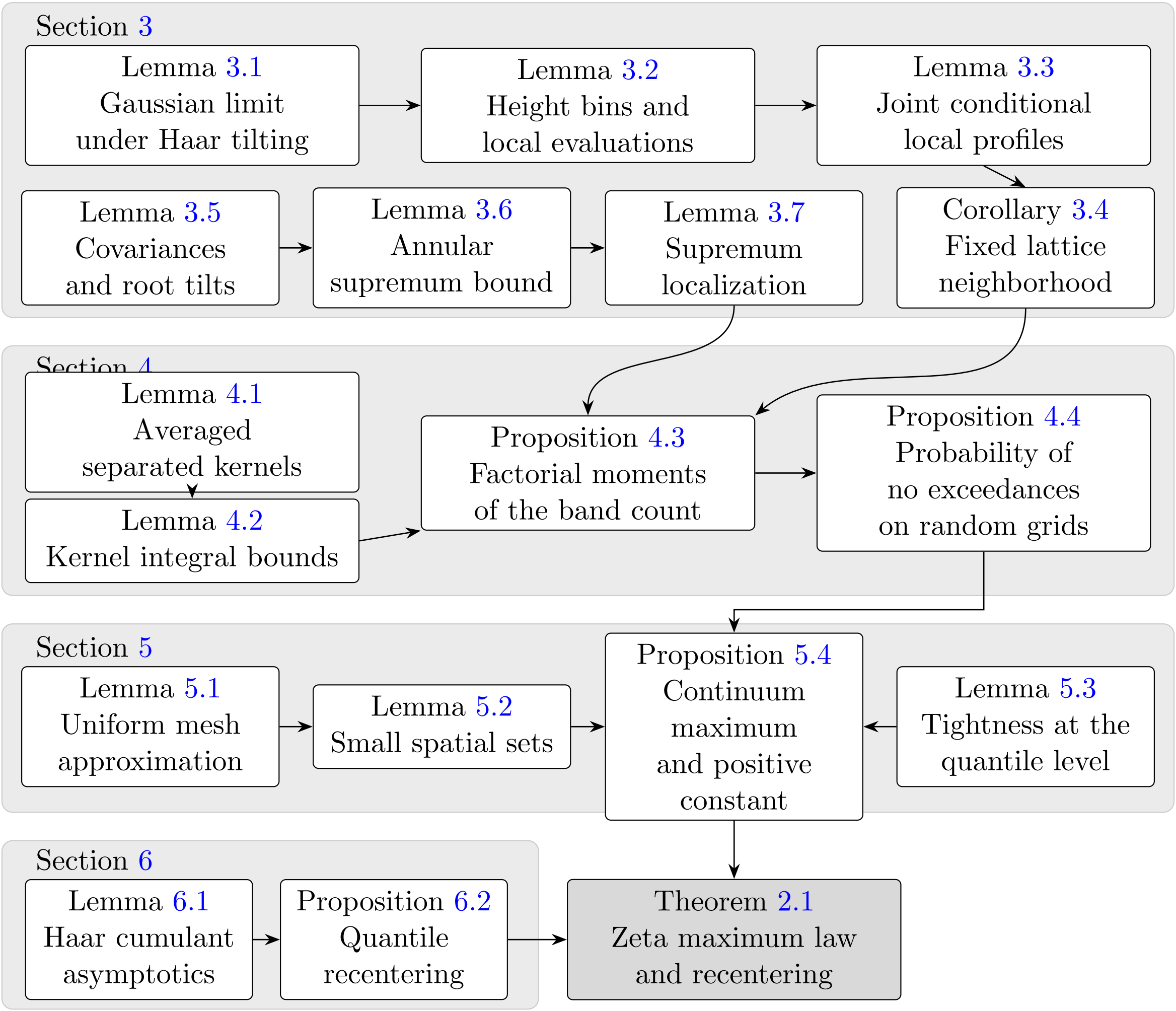
\begin{figure}[!ht]
\centering
\begin{tikzpicture}[x = 1cm, y = 1cm, every node/.style = {font = \small, align = center, execute at begin node = {\hyphenpenalty = 10000\exhyphenpenalty = 10000\relax}}, result/.style = {draw, rounded corners = 2pt, fill = white, text width = 4.1cm, minimum height = 1.05cm, inner sep = 4pt}, main/.style = {result, fill = black!15}, dependency/.style = {-{Stealth[length = 2mm]}, line width = 0.5pt}]
\path[fill = black!8, draw = black!20, rounded corners = 4pt] ( - 2.5,3.15) rectangle (12.9, - 0.9854);
\path[fill = black!8, draw = black!20, rounded corners = 4pt] ( - 2.5, - 1.3582) rectangle (12.9, - 4.6388);
\path[fill = black!8, draw = black!20, rounded corners = 4pt] ( - 2.5, - 5.0116) rectangle (12.9, - 7.4861);
\path[fill = black!8, draw = black!20, rounded corners = 4pt] ( - 2.5, - 7.8589) rectangle (4.55, - 10.0725);
\node[font = \small, anchor = west] at ( - 2.2,2.85) {Section~\ref{sec:high.points}};
\node[font = \small, anchor = west] at ( - 2.2, - 1.6382) {Section~\ref{sec:void}};
\node[font = \small, anchor = west] at ( - 2.2, - 5.3116) {Section~\ref{sec:continuum}};
\node[font = \small, anchor = west] at ( - 2.2, - 8.1089) {Section~\ref{sec:recentering}};

\node[result] (tilted) at (0,1.8) {Lemma~\ref{lem:local.profile.tilted} \\Gaussian limit \\under Haar tilting};
\node[result] (bins) at (5.2,1.8) {Lemma~\ref{lem:local.profile.bins} \\Height bins and \\local evaluations};
\node[result] (profile) at (10.4,1.8) {Lemma~\ref{lem:local.profile} \\Joint conditional \\local profiles};
\node[result, text width = 3.1cm] (covariance) at ( - 0.55, - 0.0718) {Lemma~\ref{lem:annular.covariance} \\Covariances \\and root tilts};
\node[result, text width = 3.1cm] (annular) at (3.28, - 0.0718) {Lemma~\ref{lem:annular.supremum} \\Annular \\supremum bound};
\node[result, text width = 3.1cm] (localization) at (7.12, - 0.0718) {Lemma~\ref{lem:conditional.supremum.localization} \\Supremum \\localization};
\node[result, text width = 3.1cm] (lattice) at (10.95, - 0.0718) {Corollary~\ref{cor:local.profile.lattice} \\Fixed lattice \\neighborhood};
\node[result] (average) at (0, - 2.4930) {Lemma~\ref{max:kernel.average} \\Averaged separated kernels};
\node[result] (retained) at (0, - 3.9204) {Lemma~\ref{max:retained.kernel} \\Kernel integral bounds};
\node[result] (factorial) at (5.2, - 3.0296) {Proposition~\ref{max:factorial.moments} \\Factorial moments \\of the band count};
\node[result] (void) at (10.4, - 3.0296) {Proposition~\ref{max:grid.void} \\Probability of no exceedances \\on random grids};
\node[result, text width = 3.1cm] (mesh) at ( - 0.55, - 6.3616) {Lemma~\ref{max:mesh.error} \\Uniform mesh \\approximation};
\node[result, text width = 3.1cm] (smallsets) at (3.28, - 6.3616) {Lemma~\ref{max:small.spatial.sets} \\Small spatial sets};
\node[result, text width = 3.1cm] (tightness) at (10.95, - 6.3616) {Lemma~\ref{max:quantile.tightness} \\Tightness at the \\quantile level};
\node[result, text width = 3.1cm] (continuum) at (7.12, - 6.3616) {Proposition~\ref{max:mesh.to.continuum} \\Continuum maximum \\and positive constant};
\node[result, text width = 2.7cm] (cumulants) at ( - 0.7, - 9.1589) {Lemma~\ref{lem:haar.cumulants} \\Haar cumulant \\asymptotics};
\node[result, text width = 2.7cm] (recentering) at (2.65, - 9.1589) {Proposition~\ref{prop:maximum.recentering} \\Quantile \\recentering};
\node[main] (maximum) at (7.12, - 9.1589) {Theorem~\ref{thm:maximum} \\Zeta maximum law \\and recentering};

\draw[dependency] (tilted.east) -- (bins.west);
\draw[dependency] (bins.east) -- (profile.west);
\draw[dependency] (profile.south) -- (lattice.north);
\draw[dependency] (covariance.east) -- (annular.west);
\draw[dependency] (annular.east) -- (localization.west);
\draw[dependency] (localization.south) to[out = -90, in = 90] (factorial.north);
\draw[dependency] (lattice.south) to[out = -90, in = 40] (factorial.north east);
\draw[dependency] (average.south) -- (retained.north);
\draw[dependency] (retained.east) -- (factorial.south west);
\draw[dependency] (factorial.east) -- (void.west);
\draw[dependency] (void.south) -- (10.4, - 4.8252) -- (7.12, - 4.8252) -- (continuum.north);
\draw[dependency] (mesh.east) -- (smallsets.west);
\draw[dependency] (smallsets.east) -- (continuum.west);
\draw[dependency] (tightness.west) -- (continuum.east);
\draw[dependency] (continuum.south) -- (7.12, - 7.6725) -- (maximum.north);
\draw[dependency] (cumulants.east) -- (recentering.west);
\draw[dependency] (recentering.east) -- (maximum.west);
\end{tikzpicture}
\caption{Selected principal dependencies among the proved results of the present paper. An arrow indicates that the target proof uses the source result or an estimate established in its proof. Auxiliary dependencies and arrows already represented by directed paths are omitted. White boxes are intermediate results, grouped by section in shaded regions, and the gray box is the main theorem.}
\label{fig:result.dependencies}
\end{figure}

{\bf Dependencies.} Figure~\ref{fig:result.dependencies} summarizes selected principal dependencies among the proved results in the present paper.

\subsection{Additional notation}

Write $\N \leqdef \{1,2,\ldots\}$ and $\N_{0} \leqdef \{0,1,2,\ldots\}$. For a collection of sets $\mathcal{C}$, write $\bigcup\mathcal{C} \leqdef \bigcup_{A\in\mathcal{C}}A$. For auxiliary random variables, $\PP$ and $\EE$ also denote probability and expectation when the law is clear from the context. The symbol $\star$ always belongs to $\{\mathrm{A}, \mathrm{H}\}$: $\star = \mathrm{A}$ specifies the Euler field $F^{\mathrm{A}}$ under $\PP_{\mathrm{A}}$, and $\star = \mathrm{H}$ specifies the Haar field $F^{\mathrm{H}}$ under $\PP_{\mathrm{H}}$. For $x \in \R$, write $x_{+} \leqdef \max\{x,0\}$. For $d \in \N$, $\bb{0}_{d}$ and $\bb{1}_{d}$ denote the $d$-dimensional vectors of zeros and ones, respectively. For a vector $\bb{x}$, $\|\bb{x}\|_{2}$ denotes its Euclidean norm. For a matrix $A$, $\|A\|_{\mathrm{op}}$ denotes the operator norm induced by the Euclidean norm, and $\mathrm{Id}_{k}$ denotes the $k \times k$ identity matrix. In the proofs, limits are taken as $T \to \infty$ along deterministic functions $(\sig,X)$ satisfying \eqref{eq:fixed.parameters}, unless another limit is specified. The notation $f = \OO(g)$ means that $|f|/g$ is bounded, and $f = \oo(g)$ means that $f/g \to 0$, for positive $g$. The relation $f \asymp g$ means that $cg \leq f \leq Cg$ for some constants $0 < c \leq C < \infty$ and all sufficiently large $T$. Positive finite constants, usually denoted by $c$, $C$, $C'$, etc., may change from one occurrence to the next. They may depend on the fixed parameters $\theta$, $\phi$, $\Phi$, and $C_{\sig}$, and on explicitly fixed spatial, height, and lattice cutoffs, but not on $T$ or on the particular values of $\sig$ and $X$ in the stated range.

\section{The Euler and Haar fields near high points}\label{sec:high.points}

\subsection{Joint conditional local profiles}\label{sec:local.profile}

Throughout this subsection, fix $k \in \N$, $d_{0},R,A > 0$, and $b_{-} > 0$. Let $\bb{h} = (h_{1},\ldots,h_{k})$ be a tuple of root locations satisfying $|h_{i} - h_{j}| \geq d_{0}$ for $i \neq j$ and $[h_{i} - R\ell_{T},h_{i} + R\ell_{T}] \subseteq I_{\theta}$ for every $i \in \{1,\ldots,k\}$. Let $\bb{J} = (J_{1},\ldots,J_{k})$ be the corresponding height bins, where $J_{i} = [a_{i},b_{i}] \subseteq [ - A,A]$ and $b_{i} - a_{i} \geq b_{-}$. For $\star \in \{\mathrm{A}, \mathrm{H}\}$ and $i \in \{1,\ldots,k\}$, define the joint height-bin event and the local increments around each root by
\begin{equation}\label{eq:local.profile.event.increments}
E_{T}(\bb{h}, \bb{J}) \leqdef \bigcap_{i=1}^{k}\{F_{h_{i}}^{\star} - u_{T} \in J_{i}\}, \qquad
Z_{T,h_{i}}^{\star}(v) \leqdef F_{h_{i} + \ell_{T}v}^{\star} - F_{h_{i}}^{\star}, \qquad v \in [ - R,R].
\end{equation}

For each $i \in \{1,\ldots,k\}$, fix $m_{i} \in \N$ and offsets $v_{i,a} \in [ - R,R]$, for $a \in \{1,\ldots,m_{i}\}$, and let $m \leqdef \smash{\sum_{i=1}^{k}m_{i}}$. Define the root and increment vectors by
\begin{equation}\label{eq:local.profile.root.increment.vectors}
\bb{Y}^{\star} \leqdef (Y_{i}^{\star})_{i=1}^{k} = (F_{h_{i}}^{\star})_{i=1}^{k}, \qquad
\bb{Z}^{\star} \leqdef (Z_{i,a}^{\star})_{i,a} = (Z_{T,h_{i}}^{\star}(v_{i,a}))_{i,a},
\end{equation}
where the coordinates of $\bb{Z}^{\star}$ are ordered first by $i$ and then by $a$. For Haar quantities, omit the superscript and write $\bb{Y} \leqdef \bb{Y}^{\mathrm{H}}$ and $\bb{Z} \leqdef \bb{Z}^{\mathrm{H}}$. For every bounded complex-valued random variable $Q$ measurable with respect to the Haar phases $(U_{p})_{p}$, define the positively tilted Haar expectation by
\begin{equation}\label{eq:local.profile.tilt.definition}
\EE_{\mathrm{H}}^{\lambda_{T}\bb{1}_{k},\bb{0}_{m}}[Q]
\leqdef \frac{\EE_{\mathrm{H}}[Qe^{\lambda_{T}\sum_{i=1}^{k}Y_{i}}]}{\EE_{\mathrm{H}}[e^{\lambda_{T}\sum_{i=1}^{k}Y_{i}}]}.
\end{equation}
The next two lemmas identify the local profiles under this tilt and then under conditioning on the height bins. Lemma~\ref{lem:local.profile} subsequently establishes convergence of the entire local paths.

\begin{lemma}[Joint Gaussian limit under Haar tilting]\label{lem:local.profile.tilted}
For the roots $h_{1},\ldots,h_{k}$ and offsets $(v_{i,a})_{1\leq i\leq k,\,1\leq a\leq m_{i}}$ specified above, as $T \to \infty$,
\begin{equation}\label{eq:local.profile.covariance.blocks}
\Cov_{\mathrm{H}}(Y_{i},Z_{j,a}) \to \ind_{\{i=j\}}\mathfrak{c}(v_{j,a}), \qquad \Cov_{\mathrm{H}}(Z_{i,a},Z_{j,b}) \to \ind_{\{i=j\}}\mathfrak{C}(v_{i,a},v_{j,b}),
\end{equation}
where $\mathfrak{c}(v) \leqdef - \frac{1}{4}\log(1 + v^{2})$ and $\mathfrak{C}(v,w) = \mathfrak{c}(v - w) - \mathfrak{c}(v) - \mathfrak{c}(w)$ is the covariance function of $G$ from \eqref{eq:main.profile.covariance}. The root covariance matrix $\Sigma_{T} \leqdef \smash{(\Cov_{\mathrm{H}}(Y_{i},Y_{j}))_{i,j=1}^{k}}$ satisfies
\begin{equation}\label{eq:local.profile.root.matrix}
\Sigma_{T} = V_{T}\mathrm{Id}_{k} + \OO_{k,d_{0}}(\log t), \qquad V_{T} = \frac{s_{T}}{2} + \OO(1), \qquad s_{T} \sim t,
\end{equation}
where the matrix error is in operator norm. In particular, for some constants $0 < c \leq C < \infty$ and all sufficiently large $T$,
\begin{equation}\label{eq:local.profile.root.matrix.bounds}
ct\mathrm{Id}_{k} \preceq \Sigma_{T} \preceq Ct\mathrm{Id}_{k},
\end{equation}
where $\preceq$ denotes the Loewner order on real symmetric matrices. For $\bb{y} \in \R^{k}$ and $\bb{\xi} \in \R^{m}$, let $\bb{\omega} \leqdef \bb{y}/\sqrt{V_{T}}$. Then, as $T \to \infty$,
\begin{equation}\label{eq:local.profile.central.transform}
\EE_{\mathrm{H}}^{\lambda_{T}\bb{1}_{k},\bb{0}_{m}}\left[\exp\{\ii\bb{\omega}\cdot(\bb{Y} - u_{T}\bb{1}_{k}) + \ii\bb{\xi}\cdot\bb{Z}\}\right]
\to \exp( - \|\bb{y}\|_{2}^{2}/2)\varphi_{D_{\theta}}(\bb{\xi}),
\end{equation}
where
\[
\varphi_{D_{\theta}}(\bb{\xi}) \leqdef \exp\left\{\ii\beta_{c}(\theta)\sum_{i,a}\xi_{i,a}\mathfrak{c}(v_{i,a}) - \frac{1}{2}\sum_{i}\sum_{a,b}\xi_{i,a}\xi_{i,b}\mathfrak{C}(v_{i,a},v_{i,b})\right\}.
\]
Let $D_{\theta,1}, \ldots,D_{\theta,k}$ be independent copies of the Gaussian profile $D_{\theta}$ defined in \eqref{eq:main.profile}. Then $\varphi_{D_{\theta}}$ is the characteristic function of $(D_{\theta,i}(v_{i,a}))_{i,a}$. Thus, under the tilted Haar law, the normalized root heights $(\bb{Y} - u_{T}\bb{1}_{k})/\sqrt{V_{T}}$ and the local increments $\bb{Z}$ converge jointly to a standard Gaussian vector and the indicated evaluations of the independent profiles $D_{\theta,1},\ldots,D_{\theta,k}$, with the Gaussian vector independent of all the profiles. The covariance limits are uniform over the permitted roots, and the characteristic-function convergence is also uniform for $\bb{y}$ and $\bb{\xi}$ in fixed compact sets.
\end{lemma}

\begin{proof}[Proof of Lemma~\ref{lem:local.profile.tilted}]
The proof is divided into two steps. Step~1 establishes the covariance limits and root covariance bounds for the Haar field. Step~2 uses a primewise cumulant expansion to derive the joint Gaussian limit under the root-height tilt and identify the independent local profiles $D_{\theta,1},\ldots,D_{\theta,k}$.

\medskip
{\bf Step 1. Microscopic covariances.} Extend the Haar field $F^{\mathrm{H}}$ to all real spatial arguments by the same finite sum in \eqref{eq:Haar.field}. For $d \in \R$, let $\mathfrak{C}_{T}(d) \leqdef \EE_{\mathrm{H}}[F_{0}^{\mathrm{H}}F_{d}^{\mathrm{H}}]$. Orthogonality of the prime powers gives
\begin{equation}\label{eq:local.profile.covariance.sum}
\mathfrak{C}_{T}(d) = \frac{1}{2}\sum_{p^{j}\leq X}\frac{\cos(jd\log p)}{j^{2}p^{2j\sig}}, \qquad \mathfrak{C}_{T}(0) = V_{T}.
\end{equation}

The first required limit, obtained from \eqref{eq:local.profile.covariance.sum}, is
\begin{equation}\label{eq:local.profile.covariance.drop}
\mathfrak{C}_{T}(\ell_{T}v) - V_{T} = \frac{1}{2}\sum_{j=1}^{\infty}\sum_{p\leq X^{1/j}}\frac{\cos(j\ell_{T}v\log p) - 1}{j^{2}p^{2j\sig}} \to \mathfrak{c}(v) = - \frac{1}{4}\log(1 + v^{2}),
\end{equation}
uniformly for $|v| \leq R$. The contribution of the powers $j \geq 2$ tends uniformly to zero, since
\[
\left|\frac{\cos(j\ell_{T}v\log p) - 1}{j^{2}p^{2j\sig}}\right|
\leq \frac{2}{j^{2}p^{j}}, \qquad \sum_{p}\sum_{j\geq2}\frac{1}{j^{2}p^{j}} < \infty,
\]
and the supremum over $|v| \leq R$ of every fixed summand tends to zero. The prime number theorem with its classical error term, as stated in \citet[Theorem~6.9, Eq.~(6.14), p.~179]{MR2378655}, gives
\begin{equation}\label{eq:local.profile.prime.number.theorem}
\pi(y) = \operatorname{Li}(y) + \OO\big(y\exp( - c\sqrt{\log y})\big), \qquad y \to \infty,
\end{equation}
where $\pi(y)$ counts primes not exceeding $y$, $\operatorname{Li}(y) = \int_{2}^{y}\rd{u}/\log u$, and $c > 0$ is an absolute constant. For the terms with $j = 1$ in \eqref{eq:local.profile.covariance.drop}, partial summation using \eqref{eq:local.profile.prime.number.theorem} gives, uniformly for $|v| \leq R$,
\[
\begin{aligned}
&\frac{1}{2}\sum_{p\leq X}\frac{e^{-\ell_{T}\log p}}{p}\{\cos(\ell_{T}v\log p) - 1\}
= \frac{1}{2}\int_{2}^{X}\frac{e^{-\ell_{T}\log u}}{u\log u}\{\cos(\ell_{T}v\log u) - 1\}\rd{u} + \oo(1) \\
&\quad= \frac{1}{2}\int_{\ell_{T}\log2}^{\ell_{T}\log X}e^{-y}\frac{\cos(vy) - 1}{y}\rd{y} + \oo(1) = \frac{1}{2}\int_{0}^{\infty}e^{-y}\frac{\cos(vy) - 1}{y}\rd{y} + \oo(1).
\end{aligned}
\]
Its value follows by differentiating with respect to $v$:
\[
\frac{\rd}{\rd{v}}\int_{0}^{\infty}e^{-y}\frac{\cos(vy) - 1}{y}\rd{y} = - \int_{0}^{\infty}e^{-y}\sin(vy)\rd{y} = - \frac{v}{1 + v^{2}},
\]
so that
\[
\int_{0}^{\infty}e^{-y}\frac{\cos(vy) - 1}{y}\rd{y} = - \int_{0}^{v}\frac{z}{1 + z^{2}}\rd{z} = - \frac{1}{2}\log(1 + v^{2}).
\]
This proves \eqref{eq:local.profile.covariance.drop}.

The second required limit is
\begin{equation}\label{eq:local.profile.remote.drop}
\mathfrak{C}_{T}(d + \ell_{T}v) - \mathfrak{C}_{T}(d) = \oo(1),
\end{equation}
uniformly for $d_{0} \leq |d| \leq 2|I_{\theta}|$ and $|v| \leq R$. For $j = 1$ and $\log p \leq t^{3}$, use $|\cos a - \cos b| \leq |a - b|$ and $p^{-2\sig} \leq p^{-1}$ to obtain
\[
\left|\frac{1}{2}\sum_{\substack{p\leq X \\ \log p\leq t^{3}}}\frac{\cos((d + \ell_{T}v)\log p) - \cos(d\log p)}{p^{2\sig}}\right|
\leq \frac{\ell_{T}|v|}{2}\sum_{\substack{p\leq X \\ \log p\leq t^{3}}}\frac{\log p}{p^{2\sig}}
\leq C_{R}\ell_{T}t^{3}\log t = \oo(1),
\]
where the last inequality follows from Mertens' estimate $\smash{\sum_{p\leq\exp(t^{3})}p^{-1} = 3\log t + \OO(1)}$. For $\log p > t^{3}$, the remote covariance estimate in \citet[proof of Lemma~4.8, Eq.~(4.116)]{OuimetExpanding2026} bounds each covariance by $\oo(1)$, uniformly at separations at least $d_{0}/2$. This applies to $d$ and $d + \ell_{T}v$ for all sufficiently large $T$. The remaining powers $j \geq 2$ are handled by dominated convergence, now using the bound
\[
|\cos(j(d + \ell_{T}v)\log p) - \cos(jd\log p)| \leq \min\{2,jR\ell_{T}\log p\}.
\]
This proves \eqref{eq:local.profile.remote.drop}, including its uniformity.

The arguments proving \eqref{eq:local.profile.covariance.drop} and \eqref{eq:local.profile.remote.drop} apply with any fixed bound on the offsets and any fixed positive lower bound on the separation. In particular, both estimates also hold uniformly for $|v| \leq 2R$, with $d_{0}$ replaced by $d_{0}/2$ in the latter estimate.

For the roots and increments in \eqref{eq:local.profile.root.increment.vectors}, let $d = h_{j} - h_{i}$. Their covariance blocks are exactly, for $i,j \in \{1,\ldots,k\}$,
\begin{equation}\label{eq:local.profile.covariance.identities}
\begin{aligned}
\Cov_{\mathrm{H}}(Y_{i},Z_{j,a})
&= \mathfrak{C}_{T}(d + \ell_{T}v_{j,a}) - \mathfrak{C}_{T}(d), \\
\Cov_{\mathrm{H}}(Z_{i,a},Z_{j,b})
&= \mathfrak{C}_{T}(d + \ell_{T}(v_{j,b} - v_{i,a})) - \mathfrak{C}_{T}(d - \ell_{T}v_{i,a}) - \mathfrak{C}_{T}(d + \ell_{T}v_{j,b}) + \mathfrak{C}_{T}(d).
\end{aligned}
\end{equation}
For $i = j$, applying \eqref{eq:local.profile.covariance.drop} to the second line of \eqref{eq:local.profile.covariance.identities} gives the limit $\mathfrak{c}(v_{j,b} - v_{i,a}) - \mathfrak{c}(v_{i,a}) - \mathfrak{c}(v_{j,b}) = \mathfrak{C}(v_{i,a},v_{j,b})$. For $i \neq j$, grouping the second line of \eqref{eq:local.profile.covariance.identities} into two covariance differences and applying \eqref{eq:local.profile.remote.drop} shows that it tends to zero. This proves \eqref{eq:local.profile.covariance.blocks}. To prove \eqref{eq:local.profile.root.matrix}, split the off-diagonal covariance sum at $\log p = t^{3}$. Its lower-prime part is $\OO(\log t)$, its upper-prime part is $\oo(1)$, and its higher-power part is $\OO(1)$. Since the dimension $k$ is fixed and every diagonal entry is $V_{T}$, the matrix error is $\OO(\log t)$ in operator norm. The remaining asymptotics in \eqref{eq:local.profile.root.matrix} follow from \eqref{eq:field.variance} and \eqref{eq:scale.asymptotics}, and imply \eqref{eq:local.profile.root.matrix.bounds}.

\medskip
{\bf Step 2. The tilted local characteristic function.} For complex $\bb{z}$ close to $\lambda_{T}\bb{1}_{k}$ and bounded complex $\bb{w}$, independence over primes gives
\[
\EE_{\mathrm{H}}[e^{\bb{z}\cdot\bb{Y} + \bb{w}\cdot\bb{Z}}]
= \prod_{p}\EE_{\mathrm{H}}[e^{\bb{z}\cdot\bb{Y}_{p} + \bb{w}\cdot\bb{Z}_{p}}],
\]
where $\bb{Y}_{p}$ and $\bb{Z}_{p}$ denote the contributions of the prime $p$.

Let $A_{p} \leqdef \bb{z}\cdot\bb{Y}_{p} + \bb{w}\cdot\bb{Z}_{p}$. Each coordinate of $\bb{Y}_{p}$ is centered because $\EE_{\mathrm{H}}[U_{p}^{j}] = 0$ for every $j \in \N$, and the same holds for $\bb{Z}_{p}$, whose coordinates are differences of prime contributions. Since the vectors have fixed dimensions and $\bb{z},\bb{w}$ range over bounded sets,
\[
\EE_{\mathrm{H}}[A_{p}] = 0, \qquad
|A_{p}| \leq C\sum_{j=1}^{\infty}\frac{p^{-j\sig}}{j} \leq C\frac{p^{-\sig}}{1 - p^{-\sig}} \leq Cp^{-1/2}.
\]
Taylor expansion of the exponential gives $|e^{A_{p}} - 1 - A_{p} - A_{p}^{2}/2| \leq |A_{p}|^{3}e^{|A_{p}|}/6$. Taking expectations and using $\EE_{\mathrm{H}}[A_{p}] = 0$ yields
\[
\EE_{\mathrm{H}}[e^{A_{p}}]
= 1 + \frac{1}{2}\EE_{\mathrm{H}}[A_{p}^{2}] + \OO(p^{-3/2}).
\]
Applying the Taylor expansion $\log(1 + x) = x + \OO(|x|^{2})$ with $x = \EE_{\mathrm{H}}[e^{A_{p}}] - 1 = \OO(p^{-1})$ then gives, for all sufficiently large primes $p$,
\begin{equation}\label{eq:local.profile.prime.expansion}
\log\EE_{\mathrm{H}}[e^{\bb{z}\cdot\bb{Y}_{p} + \bb{w}\cdot\bb{Z}_{p}}]
= \frac{1}{2}\EE_{\mathrm{H}}[(\bb{z}\cdot\bb{Y}_{p} + \bb{w}\cdot\bb{Z}_{p})^{2}] + r_{T,p}(\bb{z},\bb{w}), \qquad
|r_{T,p}(\bb{z},\bb{w})| \leq Cp^{-3/2}.
\end{equation}
For every fixed prime $p$, using $|v_{i,a}| \leq R$ and $|e^{-\ii x} - 1| \leq |x|$ gives
\[
|Z_{i,a,p}| \leq \ell_{T}R\log p\sum_{j=1}^{\infty}p^{-j\sig} \leq \frac{\ell_{T}R\log p}{\sqrt{p} - 1} \to 0.
\]
Thus, for each fixed prime $p$, $\bb{Z}_{p} \to \bb{0}_{m}$ uniformly in the root locations $\bb{h}$ and in the value of $U_{p}$ on the unit circle. Consequently, $r_{T,p}(\bb{z},\bb{w}) - r_{T,p}(\bb{z},\bb{0}_{m}) \to 0$, while $\smash{\sum_{p}p^{-3/2}} < \infty$ supplies a summable majorant.

Choose a fixed threshold $P$, independent of $T$, such that the expansion and remainder bound in \eqref{eq:local.profile.prime.expansion} hold for every prime $p > P$. For each prime $p \leq P$, Jensen's inequality gives $\EE_{\mathrm{H}}[e^{\lambda_{T}\bb{1}_{k}\cdot\bb{Y}_{p}}] \geq 1$. Since $\bb{Y}_{p}$ is uniformly bounded and $\bb{Z}_{p} \to \bb{0}_{m}$ uniformly, there exists a fixed $\delta > 0$, independent of $T$, such that the transforms remain bounded away from zero for $\max_{i\in\{1,\ldots,k\}}|z_{i} - \lambda_{T}| < 2\delta$ and bounded $\bb{w}$, for all sufficiently large $T$. Their logarithmic differences and the quadratic terms involving $\bb{Z}_{p}$ are therefore $\oo(1)$.

Combining these finitely many contributions with dominated convergence for the remainder differences over $p > P$ gives a total remainder of $\oo(1)$. Independence over primes allows the logarithmic transforms to be summed, while the quadratic terms sum to the covariances of $\bb{Y}$ and $\bb{Z}$, giving
\begin{equation}\label{eq:local.profile.cumulant.difference}
\begin{aligned}
&\log\EE_{\mathrm{H}}[e^{\bb{z}\cdot\bb{Y} + \bb{w}\cdot\bb{Z}}] - \log\EE_{\mathrm{H}}[e^{\bb{z}\cdot\bb{Y}}] \\
&\qquad= \sum_{i,j,a}z_{i}w_{j,a}\Cov_{\mathrm{H}}(Y_{i},Z_{j,a}) + \frac{1}{2}\sum_{i,a,j,b}w_{i,a}w_{j,b}\Cov_{\mathrm{H}}(Z_{i,a},Z_{j,b}) + \oo(1).
\end{aligned}
\end{equation}
Only purely imaginary $\bb{w}$ in fixed compact sets are needed below. The identity can equivalently be exponentiated, so no global choice of logarithm is required.

Let $H_{T}(\bb{z}) \leqdef \log\EE_{\mathrm{H}}[e^{\bb{z}\cdot\bb{Y}}]$ denote the root logarithmic transform. Choose the analytic branch of the logarithm defining $H_{T}$ on $\{\bb{z} \in \C^{k}: \max_{i\in\{1,\ldots,k\}}|z_{i} - \lambda_{T}| < 2\delta\}$ so that $H_{T}(\lambda_{T}\bb{1}_{k}) = \log\EE_{\mathrm{H}}[\smash{e^{\lambda_{T}\sum_{i}Y_{i}}}]$, where the logarithm on the right is the ordinary real logarithm of the positive expectation. Summing \eqref{eq:local.profile.prime.expansion} with $\bb{w} = \bb{0}_{m}$ and including the terms with $p \leq P$ gives
\begin{equation}\label{eq:local.profile.root.cumulant}
H_{T}(\bb{z}) = \frac{1}{2}\bb{z}^{\mathsf{T}}\Sigma_{T}\bb{z} + R_{T}(\bb{z}), \qquad
\sup_{\substack{\bb{z}\in\C^{k} \\ \max_{i\in\{1,\ldots,k\}}|z_{i} - \lambda_{T}|<2\delta}}|R_{T}(\bb{z})| \leq C.
\end{equation}
The remainder $R_{T}$ is analytic there, and its uniform bound follows from the bounded contributions of the finitely many primes $p \leq P$ and the summable bound $Cp^{-3/2}$ for $p > P$. For $\max_{i\in\{1,\ldots,k\}}|z_{i} - \lambda_{T}| < \delta$, coordinate circles of radius $\delta$ remain inside the larger neighborhood where $|R_{T}| \leq C$. Cauchy's estimates therefore give $|\partial_{i}R_{T}(\bb{z})| \leq C/\delta$ and $|\partial_{i}\partial_{j}R_{T}(\bb{z})| \leq 2C/\delta^{2}$ for $i,j \in \{1,\ldots,k\}$. Consequently,
\begin{equation}\label{eq:local.profile.remainder.derivatives}
\begin{aligned}
\|\nabla R_{T}(\bb{z})\|_{2}
&\leq \sqrt{k}\max_{i}|\partial_{i}R_{T}(\bb{z})| \leq \frac{C\sqrt{k}}{\delta} = \OO(1), \\
\|\nabla^{2}R_{T}(\bb{z})\|_{\mathrm{op}}
&\leq k\max_{i,j}|\partial_{i}\partial_{j}R_{T}(\bb{z})| \leq \frac{2Ck}{\delta^{2}} = \OO(1),
\end{aligned}
\end{equation}
since $k$, $\delta$, and $C$ are independent of $T$.

Since $\lambda_{T}V_{T} = u_{T}$ and $\lambda_{T} = \OO(1)$ by \eqref{eq:height.asymptotics}, differentiating $H_{T}$ in \eqref{eq:local.profile.root.cumulant} and using \eqref{eq:local.profile.root.matrix} and \eqref{eq:local.profile.remainder.derivatives} yields
\begin{equation}\label{eq:local.profile.root.gradient}
\begin{aligned}
\|\nabla H_{T}(\lambda_{T}\bb{1}_{k}) - u_{T}\bb{1}_{k}\|_{2}
&= \|\lambda_{T}(\Sigma_{T} - V_{T}\mathrm{Id}_{k})\bb{1}_{k} + \nabla R_{T}(\lambda_{T}\bb{1}_{k})\|_{2} \\
&\leq |\lambda_{T}|\sqrt{k}\|\Sigma_{T} - V_{T}\mathrm{Id}_{k}\|_{\mathrm{op}} + \|\nabla R_{T}(\lambda_{T}\bb{1}_{k})\|_{2}
= \OO(\log t).
\end{aligned}
\end{equation}
For bounded $\bb{y}$, let $\bb{\omega} = \bb{y}/\sqrt{V_{T}}$, and recall that $\lambda_{T} = u_{T}/V_{T}$ by \eqref{eq:height.asymptotics}. Substituting the quadratic decomposition \eqref{eq:local.profile.root.cumulant} of $H_{T}$ and expanding $R_{T}$ to first order, with remainder controlled by its Hessian bound in \eqref{eq:local.profile.remainder.derivatives}, gives
\begin{equation}\label{eq:local.profile.root.central.expansion}
\begin{aligned}
&H_{T}(\lambda_{T}\bb{1}_{k} + \ii\bb{\omega}) - H_{T}(\lambda_{T}\bb{1}_{k}) - \ii u_{T}\bb{\omega}\cdot\bb{1}_{k} \\
&\qquad= \frac{1}{2}(\lambda_{T}\bb{1}_{k} + \ii\bb{\omega})^{\mathsf{T}}\Sigma_{T}(\lambda_{T}\bb{1}_{k} + \ii\bb{\omega}) - \frac{\lambda_{T}^{2}}{2}\bb{1}_{k}^{\mathsf{T}}\Sigma_{T}\bb{1}_{k} \\
&\qquad\qquad + R_{T}(\lambda_{T}\bb{1}_{k} + \ii\bb{\omega}) - R_{T}(\lambda_{T}\bb{1}_{k}) - \ii u_{T}\bb{\omega}\cdot\bb{1}_{k} \\
&\qquad= \ii\bb{\omega}\cdot\{\lambda_{T}\Sigma_{T}\bb{1}_{k} + \nabla R_{T}(\lambda_{T}\bb{1}_{k}) - u_{T}\bb{1}_{k}\} - \frac{1}{2}\bb{\omega}^{\mathsf{T}}\Sigma_{T}\bb{\omega} + \OO(\|\bb{\omega}\|_{2}^{2}) \\
&\qquad= \ii\bb{\omega}\cdot\{\nabla H_{T}(\lambda_{T}\bb{1}_{k}) - u_{T}\bb{1}_{k}\} - \frac{1}{2}\bb{\omega}^{\mathsf{T}}\Sigma_{T}\bb{\omega} + \OO(\|\bb{\omega}\|_{2}^{2}) \\
&\qquad= - \frac{1}{2}\|\bb{y}\|_{2}^{2} + \OO\left(\frac{\log t}{\sqrt{t}}\right) + \OO\left(\frac{\log t}{t}\right),
\end{aligned}
\end{equation}
where the last equality uses the gradient estimate \eqref{eq:local.profile.root.gradient} and the covariance-matrix asymptotics in \eqref{eq:local.profile.root.matrix}.

Applying \eqref{eq:local.profile.cumulant.difference} with $\bb{z} = \lambda_{T}\bb{1}_{k} + \ii\bb{\omega}$ and $\bb{w} = \ii\bb{\xi}$ gives, as $T \to \infty$,
\begin{equation}\label{eq:local.profile.increment.exponent}
\begin{aligned}
&\log\EE_{\mathrm{H}}\left[e^{(\lambda_{T}\bb{1}_{k} + \ii\bb{\omega})\cdot\bb{Y} + \ii\bb{\xi}\cdot\bb{Z}}\right] - \log\EE_{\mathrm{H}}\left[e^{(\lambda_{T}\bb{1}_{k} + \ii\bb{\omega})\cdot\bb{Y}}\right] \\
&\quad\to \ii\beta_{c}(\theta)\sum_{i,a}\xi_{i,a}\mathfrak{c}(v_{i,a}) - \frac{1}{2}\sum_{i}\sum_{a,b}\xi_{i,a}\xi_{i,b}\mathfrak{C}(v_{i,a},v_{i,b}),
\end{aligned}
\end{equation}
where the limit follows from \eqref{eq:local.profile.covariance.blocks}, $\lambda_{T} \to \beta_{c}(\theta)$ from \eqref{eq:height.asymptotics}, and $\bb{\omega} = \bb{y}/\sqrt{V_{T}} \to \bb{0}_{k}$. Combining \eqref{eq:local.profile.root.central.expansion} and \eqref{eq:local.profile.increment.exponent} through the tilt definition \eqref{eq:local.profile.tilt.definition} proves \eqref{eq:local.profile.central.transform}. The bounds above are uniform over the permitted roots and on fixed compact sets of $\bb{y}$ and $\bb{\xi}$, proving the stated uniformity.
\end{proof}

For $k \in \N$ and $\bb{h} = (h_{1}, \ldots,h_{k})$, consider the kernel
\begin{equation}\label{eq:local.profile.kernel}
\mathcal{K}_{T}(\bb{h}) \leqdef \frac{\EE_{\mathrm{H}}\left[\exp\{\lambda_{T}\sum_{i=1}^{k} F_{h_{i}}^{\mathrm{H}}\}\right]}{\{\EE_{\mathrm{H}}[\exp(\lambda_{T} F_{0}^{\mathrm{H}})]\}^{k}}.
\end{equation}
Define the limiting height measure $\Theta(\rd{x}) \leqdef \beta_{c}(\theta) e^{-\beta_{c}(\theta) x}\rd{x}$.

\begin{lemma}[Joint height bins and finite-dimensional local profiles]\label{lem:local.profile.bins}
Consider the roots $h_{1},\ldots,h_{k}$, bins $J_{1},\ldots,J_{k}$, and offsets $(v_{i,a})_{1\leq i\leq k,\,1\leq a\leq m_{i}}$ specified above. Let $f:\R^{m} \to \R$ be a fixed bounded continuous function. For $\star \in \{\mathrm{A}, \mathrm{H}\}$, as $T \to \infty$,
\begin{equation}\label{eq:local.profile.finite.dimensional.bin}
\EE_{\star}\left[\ind_{E_{T}(\bb{h},\bb{J})}f(\bb{Z}^{\star})\right]
= N_{T}^{-k}\mathcal{K}_{T}(\bb{h})\left\{\left(\prod_{i=1}^{k}\Theta(J_{i})\right)\EE\left[f\big((D_{\theta,i}(v_{i,a}))_{i,a}\big)\right] + \oo(1)\right\}.
\end{equation}
In particular,
\begin{equation}\label{eq:local.profile.unmarked.bin}
\PP_{\star}(E_{T}(\bb{h},\bb{J}))
= N_{T}^{-k}\mathcal{K}_{T}(\bb{h})\left\{\prod_{i}\Theta(J_{i}) + \oo(1)\right\},
\end{equation}
and, conditionally on $E_{T}(\bb{h},\bb{J})$, the vector $\bb{Z}^{\star}$ converges in distribution to $(D_{\theta,i}(v_{i,a}))_{i,a}$. The errors are uniform over the permitted roots, bins, and $\star \in \{\mathrm{A}, \mathrm{H}\}$.
\end{lemma}

\begin{proof}[Proof of Lemma~\ref{lem:local.profile.bins}]
The proof is divided into three steps. Step~1 establishes a Fourier bound away from the central frequency range. Step~2 combines this bound with Fourier inversion and the one-point tail normalization to obtain the joint height-bin asymptotic with local characteristic-function marks under the Haar law. Step~3 transfers this asymptotic to the Euler law and deduces the result for bounded continuous functions of the local increments.

The covariance estimates and tilted characteristic-function limit are given by Lemma~\ref{lem:local.profile.tilted}. The transform, Fourier approximation, and one-point tail estimates from \citet[Lemmas~4.4, 4.6, and 4.7]{OuimetExpanding2026} will be used to recover the height bins and transfer the asymptotic to the Euler field.

\medskip
{\bf Step 1. A Fourier bound away from zero.} Fix $0 < \e_{0} < 1/4$ and let $\Omega_{T} \leqdef t^{\e_{0}/2}/8$. To control the contribution of frequencies beyond the central scale $t^{-1/2}$ in the asymptotic evaluation of the Fourier inversion integral, the following bound is needed:
\begin{equation}\label{eq:local.profile.fourier.envelope}
\left|\EE_{\mathrm{H}}^{\lambda_{T}\bb{1}_{k},\bb{0}_{m}}[e^{\ii\bb{\omega}\cdot\bb{Y} + \ii\bb{\xi}\cdot\bb{Z}}]\right|
\leq C_{\bb{\xi}}\exp( - ct\|\bb{\omega}\|_{2}^{2}), \qquad \|\bb{\omega}\|_{\infty} \leq 4\Omega_{T}.
\end{equation}
Choose $P_{T} = t^{B}$, with fixed $B > 2 + \e_{0}$. Under the positive tilt, the prime phases remain independent. Let $\bb{Y}_{p}$ and $\bb{Z}_{p}$ denote the contributions of the prime $p$ to the vectors $\bb{Y}$ and $\bb{Z}$. For fixed $\bb{\xi} \in \R^{m}$ and $p > P_{T}$, the prime contribution $Q_{p} \leqdef \bb{\omega}\cdot\bb{Y}_{p} + \bb{\xi}\cdot\bb{Z}_{p}$ satisfies
\[
|Q_{p}| \leq C(\|\bb{\omega}\|_{2} + \|\bb{\xi}\|_{2})p^{-1/2}
\leq C_{\bb{\xi}}t^{\e_{0}/2}t^{-B/2}
= C_{\bb{\xi}}t^{-(B - \e_{0})/2}
= \oo(1),
\]
uniformly for $\|\bb{\omega}\|_{\infty} \leq 4\Omega_{T}$, since $\|\bb{\omega}\|_{2} \leq 4\sqrt{k}\Omega_{T} = \sqrt{k}t^{\e_{0}/2}/2$. If $\smash{Q_{p}'}$ is an independent copy under the same prime tilt, then $\cos z \leq 1 - z^{2}/4$ for sufficiently small $|z|$ gives
\begin{equation}\label{eq:local.profile.prime.characteristic.bound}
\begin{aligned}
\left|\EE_{\mathrm{H}}^{\lambda_{T}\bb{1}_{k},\bb{0}_{m}}[e^{\ii Q_{p}}]\right|^{2}
&= \EE_{\mathrm{H}}^{\lambda_{T}\bb{1}_{k},\bb{0}_{m}}[\cos(Q_{p} - Q_{p}')] \\
&\leq 1 - \frac{1}{4}\EE_{\mathrm{H}}^{\lambda_{T}\bb{1}_{k},\bb{0}_{m}}[(Q_{p} - Q_{p}')^{2}]
= 1 - \frac{1}{2}\Var_{\mathrm{H}}^{\lambda_{T}\bb{1}_{k},\bb{0}_{m}}(Q_{p}).
\end{aligned}
\end{equation}
Since $|Y_{i,p}| \leq \sum_{j\geq1}p^{-j/2}/j \leq Cp^{-1/2}$ and $\lambda_{T} = \OO(1)$ by \eqref{eq:height.asymptotics}, there exists $C_{0} > 0$ such that the tilt factor and its expectation satisfy
\[
\begin{aligned}
&e^{-C_{0}p^{-1/2}} \leq \exp\left\{\lambda_{T}\sum_{i=1}^{k}Y_{i,p}\right\} \leq e^{C_{0}p^{-1/2}}, \\
&e^{-C_{0}p^{-1/2}} \leq \EE_{\mathrm{H}}\left[\exp\left\{\lambda_{T}\sum_{i=1}^{k}Y_{i,p}\right\}\right] \leq e^{C_{0}p^{-1/2}},
\end{aligned}
\]
so that the prime-tilt density satisfies
\[
e^{-2C_{0}p^{-1/2}} \leq \frac{\exp\left\{\lambda_{T}\sum_{i=1}^{k}Y_{i,p}\right\}}{\EE_{\mathrm{H}}\left[\exp\left\{\lambda_{T}\sum_{i=1}^{k}Y_{i,p}\right\}\right]}
\leq e^{2C_{0}p^{-1/2}}.
\]
Comparing the two-copy formula for the tilted variance with its Haar counterpart gives
\begin{equation}\label{eq:local.profile.prime.variance.comparison}
\Var_{\mathrm{H}}^{\lambda_{T}\bb{1}_{k},\bb{0}_{m}}(Q_{p})
= \frac{1}{2}\EE_{\mathrm{H}}^{\lambda_{T}\bb{1}_{k},\bb{0}_{m}}[(Q_{p} - Q_{p}')^{2}]
\geq e^{-Cp^{-1/2}}\Var_{\mathrm{H}}(Q_{p}).
\end{equation}
Consequently, combining \eqref{eq:local.profile.prime.characteristic.bound} with \eqref{eq:local.profile.prime.variance.comparison} and using $1 - z \leq e^{-z}$ gives
\begin{equation}\label{eq:local.profile.fourier.product.bound}
\prod_{p>P_{T}}\left|\EE_{\mathrm{H}}^{\lambda_{T}\bb{1}_{k},\bb{0}_{m}}[e^{\ii Q_{p}}]\right|
\leq \exp\left\{ - c\sum_{p>P_{T}}\Var_{\mathrm{H}}(Q_{p})\right\}.
\end{equation}
Write $\bb{Y}^{>P_{T}}$ and $\bb{Z}^{>P_{T}}$ for the corresponding sums over primes $p > P_{T}$. The covariance matrix $\Cov_{\mathrm{H}}(\smash{\sum_{p\leq P_{T}}\bb{Y}_{p}})$ of the root contributions from primes $p \leq P_{T}$ has operator norm at most $C\log\log P_{T} = \OO(\log\log t)$, so subtracting it from $\Sigma_{T}$ and applying \eqref{eq:local.profile.root.matrix.bounds} of Lemma~\ref{lem:local.profile.tilted} gives
\[
\Var_{\mathrm{H}}(\bb{\omega}\cdot\bb{Y}^{>P_{T}})
\geq \{ct - C\log\log t\}\|\bb{\omega}\|_{2}^{2}
\geq c't\|\bb{\omega}\|_{2}^{2}.
\]
Deleting independent prime contributions cannot increase the variance, so \eqref{eq:local.profile.covariance.blocks} of Lemma~\ref{lem:local.profile.tilted} also gives $\Var_{\mathrm{H}}(\bb{\xi}\cdot\bb{Z}^{>P_{T}}) \leq C\|\bb{\xi}\|_{2}^{2}$. The inequality $(a + b)^{2} \geq a^{2}/2 - b^{2}$ yields
\begin{equation}\label{eq:local.profile.large.prime.variance}
\begin{aligned}
\sum_{p>P_{T}}\Var_{\mathrm{H}}(Q_{p})
&= \Var_{\mathrm{H}}(\bb{\omega}\cdot\bb{Y}^{>P_{T}} + \bb{\xi}\cdot\bb{Z}^{>P_{T}}) \\[-2mm]
&\geq \frac{1}{2}\Var_{\mathrm{H}}(\bb{\omega}\cdot\bb{Y}^{>P_{T}}) - \Var_{\mathrm{H}}(\bb{\xi}\cdot\bb{Z}^{>P_{T}}) \geq ct\|\bb{\omega}\|_{2}^{2} - C\|\bb{\xi}\|_{2}^{2}.
\end{aligned}
\end{equation}
Substituting \eqref{eq:local.profile.large.prime.variance} into \eqref{eq:local.profile.fourier.product.bound} and using that all remaining prime characteristic factors have modulus at most one proves \eqref{eq:local.profile.fourier.envelope}.

\medskip
{\bf Step 2. Haar root bins with bounded local marks.} Define the weighted bin function and its integral by
\begin{equation}\label{eq:local.profile.weighted.bin.function}
f_{T}(\bb{x}) \leqdef e^{-\lambda_{T}\sum_{i}x_{i}}\prod_{i}\ind_{J_{i}}(x_{i}), \qquad
I_{T}(\bb{J}) \leqdef \int_{\R^{k}}f_{T}(\bb{x})\rd{\bb{x}} = \prod_{i}\int_{J_{i}}e^{-\lambda_{T}x}\rd{x}.
\end{equation}
Use the majorants and minorants in \citet[proof of Lemma~4.6, Eqs.~(4.61)--(4.63)]{OuimetExpanding2026}. For these weighted rectangular bins, they provide integrable continuous functions $f_{T}^{-},f_{T}^{+}$ satisfying
\begin{equation}\label{eq:local.profile.fourier.approximants}
\begin{aligned}
&f_{T}^{-} \leq f_{T} \leq f_{T}^{+}, \qquad \operatorname{supp}(\widehat{f}_{T}^{\pm}) \subseteq [ - 4\Omega_{T},4\Omega_{T}]^{k}, \qquad \|f_{T}^{+} - f_{T}^{-}\|_{1} \leq C\Omega_{T}^{-1}, \\
&\|f_{T}^{\pm}\|_{1} \leq C, \qquad \widehat{f}_{T}^{\pm}(\bb{0}_{k}) = I_{T}(\bb{J}) + \OO(\Omega_{T}^{-1}), \qquad |\widehat{f}_{T}^{\pm}(\bb{\omega}) - \widehat{f}_{T}^{\pm}(\bb{0}_{k})| \leq C(\Omega_{T}^{-1} + \|\bb{\omega}\|_{2}).
\end{aligned}
\end{equation}
The constants are uniform in the permitted bin endpoints and widths. The Fourier convention is $\smash{\widehat{f}}(\bb{\omega}) = \int e^{-\ii\bb{\omega}\cdot\bb{x}}f(\bb{x})\rd{\bb{x}}$.

Let $\chi_{T}(\bb{\omega},\bb{\xi})$ denote the tilted characteristic function on the left-hand side of \eqref{eq:local.profile.central.transform} of Lemma~\ref{lem:local.profile.tilted}. Fourier inversion and Fubini's theorem give
\begin{equation}\label{eq:local.profile.fourier.inversion}
\begin{aligned}
&\EE_{\mathrm{H}}^{\lambda_{T}\bb{1}_{k},\bb{0}_{m}}\left[f_{T}^{+}(\bb{Y} - u_{T}\bb{1}_{k})e^{\ii\bb{\xi}\cdot\bb{Z}}\right]
= \frac{1}{(2\pi)^{k}}\int_{[ - 4\Omega_{T},4\Omega_{T}]^{k}}\widehat{f}_{T}^{+}(\bb{\omega})\chi_{T}(\bb{\omega},\bb{\xi})\rd{\bb{\omega}} \\
&\quad= \frac{V_{T}^{-k/2}}{(2\pi)^{k}}\int_{[ - 4\Omega_{T}\sqrt{V_{T}},4\Omega_{T}\sqrt{V_{T}}]^{k}} \widehat{f}_{T}^{+}\left(\frac{\bb{y}}{\sqrt{V_{T}}}\right)\chi_{T}\left(\frac{\bb{y}}{\sqrt{V_{T}}},\bb{\xi}\right)\rd{\bb{y}}.
\end{aligned}
\end{equation}
By \eqref{eq:local.profile.central.transform} of Lemma~\ref{lem:local.profile.tilted} and \eqref{eq:local.profile.fourier.approximants}, the integrand equals $I_{T}(\bb{J})e^{-\|\bb{y}\|_{2}^{2}/2}\varphi_{D_{\theta}}(\bb{\xi})$, up to a uniform $\oo(1)$ error on bounded $\bb{y}$-sets. On the full integration region, \eqref{eq:local.profile.fourier.envelope} gives
\[
\left|\widehat{f}_{T}^{+}\left(\frac{\bb{y}}{\sqrt{V_{T}}}\right)\chi_{T}\left(\frac{\bb{y}}{\sqrt{V_{T}}},\bb{\xi}\right)\right|
\leq C_{\bb{\xi}}e^{-c\|\bb{y}\|_{2}^{2}}, \qquad
\int_{\R^{k}}e^{-\|\bb{y}\|_{2}^{2}/2}\rd{\bb{y}} = (2\pi)^{k/2}.
\]
Dominated convergence in \eqref{eq:local.profile.fourier.inversion}, with the same majorant for all permitted roots and bins, therefore proves
\[
\EE_{\mathrm{H}}^{\lambda_{T}\bb{1}_{k},\bb{0}_{m}}\left[f_{T}^{+}(\bb{Y} - u_{T}\bb{1}_{k})e^{\ii\bb{\xi}\cdot\bb{Z}}\right]
= (2\pi V_{T})^{-k/2}\{I_{T}(\bb{J})\varphi_{D_{\theta}}(\bb{\xi}) + \oo(1)\}.
\]
To replace $f_{T}^{+}$ by $f_{T}$, let $g_{T} \leqdef f_{T}^{+} - f_{T}^{-} \geq 0$. Since $|e^{\ii\bb{\xi}\cdot\bb{Z}}| = 1$, the replacement error is bounded, using \eqref{eq:local.profile.fourier.envelope} and \eqref{eq:local.profile.fourier.approximants}, by
\begin{equation}\label{eq:local.profile.fourier.replacement}
\EE_{\mathrm{H}}^{\lambda_{T}\bb{1}_{k},\bb{0}_{m}}[g_{T}(\bb{Y} - u_{T}\bb{1}_{k})] \leq \frac{\|g_{T}\|_{1}}{(2\pi)^{k}}\int_{\R^{k}}Ce^{-ct\|\bb{\omega}\|_{2}^{2}}\rd{\bb{\omega}} \leq C\Omega_{T}^{-1}V_{T}^{-k/2} = \oo(V_{T}^{-k/2}).
\end{equation}
This yields the weighted-bin identity
\begin{equation}\label{eq:local.profile.haar.weighted.bin}
\EE_{\mathrm{H}}^{\lambda_{T}\bb{1}_{k},\bb{0}_{m}}\left[f_{T}(\bb{Y} - u_{T}\bb{1}_{k})e^{\ii\bb{\xi}\cdot\bb{Z}}\right]
= (2\pi V_{T})^{-k/2}\{I_{T}(\bb{J})\varphi_{D_{\theta}}(\bb{\xi}) + \oo(1)\}.
\end{equation}
By \eqref{eq:local.profile.weighted.bin.function}, $f_{T}(\bb{Y} - u_{T}\bb{1}_{k}) = e^{-\lambda_{T}\sum_{i}(Y_{i} - u_{T})}\ind_{E_{T}(\bb{h},\bb{J})}$. Substituting this identity and then applying the tilt definition \eqref{eq:local.profile.tilt.definition} gives
\begin{equation}\label{eq:local.profile.undo.tilt}
\begin{aligned}
\EE_{\mathrm{H}}[\ind_{E_{T}(\bb{h},\bb{J})}e^{\ii\bb{\xi}\cdot\bb{Z}}]
&= \EE_{\mathrm{H}}\left[e^{\lambda_{T}\sum_{i}(Y_{i} - u_{T})}f_{T}(\bb{Y} - u_{T}\bb{1}_{k})e^{\ii\bb{\xi}\cdot\bb{Z}}\right] \\
&= e^{-k\lambda_{T}u_{T}}\EE_{\mathrm{H}}\left[e^{\lambda_{T}\sum_{i}Y_{i}}f_{T}(\bb{Y} - u_{T}\bb{1}_{k})e^{\ii\bb{\xi}\cdot\bb{Z}}\right] \\
&= e^{-k\lambda_{T}u_{T}}\EE_{\mathrm{H}}[e^{\lambda_{T}\sum_{i}Y_{i}}]\EE_{\mathrm{H}}^{\lambda_{T}\bb{1}_{k},\bb{0}_{m}}\left[f_{T}(\bb{Y} - u_{T}\bb{1}_{k})e^{\ii\bb{\xi}\cdot\bb{Z}}\right].
\end{aligned}
\end{equation}
The one-point tail estimate and exponential-moment asymptotic in \citet[proof of Lemma~4.7, Eqs.~(4.85) and~(4.88)]{OuimetExpanding2026}, together with the quantile normalization \eqref{eq:height.quantile} and $\lambda_{T} = u_{T}/V_{T}$ from \eqref{eq:height.asymptotics}, give
\[
\frac{1}{N_{T}}
= \PP_{\mathrm{H}}(F_{0}^{\mathrm{H}} > u_{T})
= \frac{e^{-u_{T}^{2}/(2V_{T})}}{\lambda_{T}\sqrt{2\pi V_{T}}} \frac{\EE_{\mathrm{H}}[e^{\lambda_{T}F_{0}^{\mathrm{H}}}]}{e^{\lambda_{T}^{2}V_{T}/2}} \{1 + \oo(1)\}
= \frac{e^{-\lambda_{T}u_{T}}\EE_{\mathrm{H}}[e^{\lambda_{T}F_{0}^{\mathrm{H}}}]}{\lambda_{T}\sqrt{2\pi V_{T}}}\{1 + \oo(1)\}.
\]
Equivalently,
\begin{equation}\label{eq:local.profile.scalar.prefactor}
N_{T}e^{-\lambda_{T}u_{T}}\EE_{\mathrm{H}}[e^{\lambda_{T}F_{0}^{\mathrm{H}}}](2\pi V_{T})^{-1/2}
= \lambda_{T}\{1 + \oo(1)\}.
\end{equation}
The definition \eqref{eq:local.profile.kernel}, together with \eqref{eq:local.profile.scalar.prefactor}, now isolates the dependence between roots:
\begin{equation}\label{eq:local.profile.joint.prefactor}
\begin{aligned}
e^{-k\lambda_{T}u_{T}}\EE_{\mathrm{H}}[e^{\lambda_{T}\sum_{i}Y_{i}}](2\pi V_{T})^{-k/2}
&= \mathcal{K}_{T}(\bb{h})\left\{e^{-\lambda_{T}u_{T}}\EE_{\mathrm{H}}[e^{\lambda_{T}F_{0}^{\mathrm{H}}}](2\pi V_{T})^{-1/2}\right\}^{k} \\
&= N_{T}^{-k}\mathcal{K}_{T}(\bb{h})\lambda_{T}^{k}\{1 + \oo(1)\}.
\end{aligned}
\end{equation}
Finally, \eqref{eq:height.asymptotics} and the definition of $I_{T}(\bb{J})$ in \eqref{eq:local.profile.weighted.bin.function} imply
\begin{equation}\label{eq:local.profile.height.bin.limit}
\lambda_{T}^{k}I_{T}(\bb{J})
= \prod_{i}\int_{J_{i}}\lambda_{T}e^{-\lambda_{T}x}\rd{x}
= \prod_{i}\Theta(J_{i}) + \oo(1),
\end{equation}
uniformly in the permitted bins. Substituting \eqref{eq:local.profile.haar.weighted.bin}, \eqref{eq:local.profile.joint.prefactor}, and \eqref{eq:local.profile.height.bin.limit} into \eqref{eq:local.profile.undo.tilt} proves
\begin{equation}\label{eq:local.profile.marked.characteristic}
\EE_{\mathrm{H}}[\ind_{E_{T}(\bb{h},\bb{J})}e^{\ii\bb{\xi}\cdot\bb{Z}}]
= N_{T}^{-k}\mathcal{K}_{T}(\bb{h})\left\{\left(\prod_{i}\Theta(J_{i})\right)\varphi_{D_{\theta}}(\bb{\xi}) + \oo(1)\right\}.
\end{equation}

\medskip
{\bf Step 3. Arithmetic transfer at the rare-event scale.} Use the root and increment vectors $\bb{Y}^{\star}$ and $\bb{Z}^{\star}$ from \eqref{eq:local.profile.root.increment.vectors} under either law. For each prime power $p^{j} \leq X$, with $j \in \N$, the corresponding coefficient has modulus $1/(jp^{j\sig}) \leq p^{-j/2}$ for a root and at most $2/(jp^{j\sig}) \leq 2p^{-j/2}$ for an increment, since an increment is the difference of two field evaluations.

The transform and exceptional-probability estimates from \citet[Lemma~4.4, Eqs.~(4.24)--(4.25)]{OuimetExpanding2026} therefore apply at root parameters $\lambda_{T} + \ii\omega_{i}$ and increment parameters $\ii\xi_{i,a}$. Their real parts are bounded, and their imaginary parts are bounded in absolute value by $t^{\e_{0}/2}$ for all sufficiently large $T$.

Let $\mathcal{A}_{T}$ be the arithmetic truncation event
\[
\mathcal{A}_{T} \leqdef \left\{\max_{1\leq i\leq k}|Y_{i}^{\mathrm{A}}| \leq t^{1 + \e_{0}}, \quad \max_{\substack{1\leq i\leq k \\ 1\leq a\leq m_{i}}}|Z_{i,a}^{\mathrm{A}}| \leq t^{1 + \e_{0}}\right\}.
\]
Uniformly over the Fourier cube, the cited transform estimate gives
\begin{equation}\label{eq:local.profile.arithmetic.transform}
\left|\EE_{\mathrm{A}}\left[\ind_{\mathcal{A}_{T}}e^{(\lambda_{T}\bb{1}_{k} + \ii\bb{\omega})\cdot\bb{Y}^{\mathrm{A}} + \ii\bb{\xi}\cdot\bb{Z}^{\mathrm{A}}}\right]
- \EE_{\mathrm{H}}[e^{(\lambda_{T}\bb{1}_{k} + \ii\bb{\omega})\cdot\bb{Y}^{\mathrm{H}} + \ii\bb{\xi}\cdot\bb{Z}^{\mathrm{H}}}]\right| \leq \exp\{ - ct^{1 + 2\e_{0}} + Ct\},
\end{equation}
and
\begin{equation}\label{eq:local.profile.arithmetic.truncation}
\PP_{\mathrm{A}}(\mathcal{A}_{T}^{c}) \leq \exp\{ - ct^{1 + 2\e_{0}} + Ct\}.
\end{equation}
Insert \eqref{eq:local.profile.arithmetic.transform} into the Fourier integral \eqref{eq:local.profile.fourier.inversion} before normalizing the exponential tilt. Since $\|\smash{\widehat{f}_{T}^{+}}\|_{\infty} \leq C$ by \eqref{eq:local.profile.fourier.approximants} and the Fourier cube has volume $\OO(\Omega_{T}^{k})$, the absolute error is at most
\[
C\Omega_{T}^{k}e^{-k\lambda_{T}u_{T}}\exp\{ - ct^{1 + 2\e_{0}} + Ct\}
\leq \exp\{ - c't^{1 + 2\e_{0}} + C't\}.
\]
The same argument with $g_{T} = f_{T}^{+} - f_{T}^{-}$ and $\bb{\xi} = \bb{0}_{m}$ controls the majorant replacement. By \eqref{eq:local.profile.fourier.replacement}, its Haar contribution is $\OO(\smash{\Omega_{T}^{-1}})$ times the prefactor
\[
\Xi_{T}(\bb{h}) \leqdef e^{-k\lambda_{T}u_{T}}\EE_{\mathrm{H}}[e^{\lambda_{T}\sum_{i}Y_{i}}](2\pi V_{T})^{-k/2}.
\]
Jensen's inequality gives $\EE_{\mathrm{H}}[e^{\lambda_{T}\sum_{i}Y_{i}}] \geq 1$. Since $u_{T} = \OO(t)$, $\lambda_{T} = \OO(1)$, and $V_{T} \asymp t$ by \eqref{eq:height.asymptotics} and \eqref{eq:local.profile.root.matrix} of Lemma~\ref{lem:local.profile.tilted},
\begin{equation}\label{eq:local.profile.prefactor.lower.bound}
\Xi_{T}(\bb{h}) \geq e^{-Ct}, \qquad
\frac{\exp\{ - c't^{1 + 2\e_{0}} + C't\}}{\Xi_{T}(\bb{h})}
\leq \exp\{ - c't^{1 + 2\e_{0}} + C''t\} \to 0.
\end{equation}
Removing the factor $\ind_{\mathcal{A}_{T}}$ from $\EE_{\mathrm{A}}[\ind_{\mathcal{A}_{T}}\ind_{E_{T}(\bb{h},\bb{J})}e^{\ii\bb{\xi}\cdot\bb{Z}^{\mathrm{A}}}]$ changes the expectation by at most $\PP_{\mathrm{A}}(\mathcal{A}_{T}^{c})$. By \eqref{eq:local.profile.arithmetic.truncation} and \eqref{eq:local.profile.prefactor.lower.bound}, this error is $\oo(\Xi_{T}(\bb{h}))$. Thus \eqref{eq:local.profile.marked.characteristic} holds with $\EE_{\star}$ and $\bb{Z}^{\star}$ in place of $\EE_{\mathrm{H}}$ and $\bb{Z}$, respectively, for each $\star \in \{\mathrm{A}, \mathrm{H}\}$.

Taking $\bb{\xi} = \bb{0}_{m}$ in \eqref{eq:local.profile.marked.characteristic} proves \eqref{eq:local.profile.unmarked.bin}. The bin assumptions imply $0 < c \leq \prod_{i}\Theta(J_{i}) \leq C$. Dividing \eqref{eq:local.profile.marked.characteristic} by \eqref{eq:local.profile.unmarked.bin} therefore gives
\begin{equation}\label{eq:local.profile.conditional.characteristic}
\EE_{\star}[e^{\ii\bb{\xi}\cdot\bb{Z}^{\star}}\mid E_{T}(\bb{h},\bb{J})]
= \frac{\{\prod_{i}\Theta(J_{i})\}\varphi_{D_{\theta}}(\bb{\xi}) + \oo(1)}{\prod_{i}\Theta(J_{i}) + \oo(1)}
\to \varphi_{D_{\theta}}(\bb{\xi}).
\end{equation}
By L\'evy's continuity theorem for characteristic functions, \eqref{eq:local.profile.conditional.characteristic} implies that, conditionally on $E_{T}(\bb{h},\bb{J})$, the vector $\bb{Z}^{\star}$ converges in distribution to $(D_{\theta,i}(v_{i,a}))_{i,a}$. The convergence is uniform: any sequence of roots and bins satisfies the same characteristic-function limit with the same continuous limiting function $\varphi_{D_{\theta}}$.

Consequently, for every fixed bounded continuous $f:\R^{m} \to \R$,
\begin{equation}\label{eq:local.profile.conditional.test.function}
\EE_{\star}[f(\bb{Z}^{\star})\mid E_{T}(\bb{h},\bb{J})]
= \EE\left[f\big((D_{\theta,i}(v_{i,a}))_{i,a}\big)\right] + \oo(1).
\end{equation}
Multiplying \eqref{eq:local.profile.conditional.test.function} by \eqref{eq:local.profile.unmarked.bin} gives \eqref{eq:local.profile.finite.dimensional.bin}.
\end{proof}

\begin{lemma}[Joint conditional local profiles under high-point conditioning]\label{lem:local.profile}
Consider the roots $h_{1},\ldots,h_{k}$ and bins $J_{1},\ldots,J_{k}$ specified above. Let $\Psi:C([ - R,R])^{k} \to \R$ be a fixed bounded continuous function. Then, as $T \to \infty$,
\begin{equation}\label{eq:local.profile.marked.bin}
\begin{aligned}
&\EE_{\star}\left[ \ind_{E_{T}(\bb{h}, \bb{J})} \Psi(Z_{T,h_{1}}^{\star}, \ldots,Z_{T,h_{k}}^{\star})\right] \\
&\qquad= N_{T}^{-k}\mathcal{K}_{T}(\bb{h}) \left\{ \left(\prod_{i=1}^{k}\Theta(J_{i})\right) \EE[\Psi(D_{\theta,1}, \ldots,D_{\theta,k})] + \oo(1) \right\}.
\end{aligned}
\end{equation}
The error is uniform over the roots, bins, and $\star \in \{\mathrm{A}, \mathrm{H}\}$ satisfying the stated conditions.
\end{lemma}

\begin{proof}[Proof of Lemma~\ref{lem:local.profile}]
Lemma~\ref{lem:local.profile.bins} gives the conditional finite-dimensional limits and the bin-probability asymptotic \eqref{eq:local.profile.unmarked.bin}. It remains to prove tightness of the conditional paths in $C([ - R,R])^{k}$. Fix $0 < \e_{0} < 1/4$, as in the proof of that lemma.

Let $\smash{L_{T,h}^{\star}(v) \leqdef \ell_{T}(F^{\star})'(h + \ell_{T}v)}$. For each prime power $p^{j} \leq X$, its coefficient is $a_{p,j}(v) \leqdef - \ii\ell_{T}(\log p)p^{-j\sig}e^{-\ii j(h + \ell_{T}v)\log p}$. Writing $x = j\eta_{\sig}\log p$, this coefficient and its derivative with respect to $v$ satisfy
\[
|a_{p,j}(v)| = \ell_{T}(\log p)p^{-j\sig} = p^{-j/2}\frac{2x}{j}e^{-x}, \quad
|a_{p,j}'(v)| = j\ell_{T}^{2}(\log p)^{2}p^{-j\sig} = p^{-j/2}\frac{4x^{2}}{j}e^{-x}.
\]
Both are bounded by $Cp^{-j/2}$. Differentiating the finite covariance sum \eqref{eq:local.profile.covariance.sum} and applying partial summation also gives
\begin{equation}\label{eq:local.profile.derivative.bounds}
\Var_{\mathrm{H}}(L_{T,h}^{\mathrm{H}}(v)) \leq C, \qquad
|\Cov_{\mathrm{H}}(F_{h_{i}}^{\mathrm{H}},L_{T,h}^{\mathrm{H}}(v))| \leq C.
\end{equation}
Indeed, their contributions with $j = 1$ are bounded respectively by
\[
C\ell_{T}^{2}\int_{\log2}^{\infty}we^{-\ell_{T}w}\rd{w}
\leq C\int_{0}^{\infty}ye^{-y}\rd{y}, \qquad
C\ell_{T}\int_{\log2}^{\infty}e^{-\ell_{T}w}\rd{w}
\leq C\int_{0}^{\infty}e^{-y}\rd{y}.
\]
The higher-prime-power contributions are bounded by the corresponding convergent sums with factors $(\log p)^{2}p^{-j}$ and $(\log p)p^{-j}$, for $j \geq 2$.

Apply \citet[Lemma~4.6, Eq.~(4.51)]{OuimetExpanding2026} with conditioning vector $\smash{(F_{h_{1}}^{\star},\ldots,F_{h_{k}}^{\star})}$ constrained to lie in $\smash{\prod_{j=1}^{k}(u_{T} + J_{j})}$, which is precisely the event $E_{T}(\bb{h},\bb{J})$ defined in \eqref{eq:local.profile.event.increments}. Take $L_{T,h_{i}}^{\star}(v)$ as the additional random variable whose conditional exponential moments are to be bounded. If $\smash{\bb{c}_{T} = (\Cov_{\mathrm{H}}(F_{h_{j}}^{\mathrm{H}},L_{T,h_{i}}^{\mathrm{H}}(v)))_{j=1}^{k}}$ is its Haar covariance vector with the roots and $\bb{a}_{T} = u_{T}\bb{1}_{k} + (a_{1},\ldots,a_{k})$, then the root covariance bound \eqref{eq:local.profile.root.matrix.bounds} from Lemma~\ref{lem:local.profile.tilted}, the derivative bounds \eqref{eq:local.profile.derivative.bounds}, and $u_{T} = \OO(t)$ give
\[
\|\Sigma_{T}^{-1}\|_{\mathrm{op}} = \OO(t^{-1}), \qquad
\|\bb{c}_{T}\|_{2} = \OO(1), \qquad \|\bb{a}_{T}\|_{2} = \OO(t).
\]
Thus the conditional mean and variance terms are uniformly bounded:
\[
|\bb{c}_{T}^{\mathsf{T}}\Sigma_{T}^{-1}\bb{a}_{T}| \leq C, \qquad
0 \leq \Var_{\mathrm{H}}(L_{T,h_{i}}^{\mathrm{H}}(v)) - \bb{c}_{T}^{\mathsf{T}}\Sigma_{T}^{-1}\bb{c}_{T} \leq C.
\]
Applying the estimate at the two extra tilts $1$ and $ - 1$ gives
\begin{equation}\label{eq:local.profile.conditional.derivative.laplace}
\EE_{\star}\left[(\exp\{L_{T,h_{i}}^{\star}(v)\} + \exp\{ - L_{T,h_{i}}^{\star}(v)\})\ind_{\mathcal{A}_{T,\mathrm{der}}^{\star}}\,\middle|\,E_{T}(\bb{h},\bb{J})\right] \leq C.
\end{equation}
For the Euler field, the common spatial truncation event from \citet[Lemma~4.5, Eqs.~(4.36)--(4.37)]{OuimetExpanding2026}, applied jointly to $F^{\mathrm{A}}$ and its scaled derivative, is
\[
\mathcal{A}_{T,\mathrm{der}}^{\mathrm{A}} \leqdef \left\{\sup_{z\in I_{\theta}}\max\{|F_{z}^{\mathrm{A}}|,|\ell_{T}(F^{\mathrm{A}})'(z)|\} \leq t^{1 + \e_{0}}\right\}.
\]
For the Haar field, let $\mathcal{A}_{T,\mathrm{der}}^{\mathrm{H}}$ be the whole probability space, so that $\ind_{\mathcal{A}_{T,\mathrm{der}}^{\mathrm{H}}} \equiv 1$.

By \eqref{eq:height.asymptotics}, $\lambda_{T} \asymp 1$, and the bin assumptions give $\prod_{i}\Theta(J_{i}) \geq c > 0$ uniformly over the permitted bins. Thus \eqref{eq:local.profile.unmarked.bin} of Lemma~\ref{lem:local.profile.bins}, \eqref{eq:local.profile.joint.prefactor}, and \eqref{eq:local.profile.prefactor.lower.bound} give, successively,
\[
\begin{aligned}
\PP_{\mathrm{A}}(E_{T}(\bb{h},\bb{J}))
&= N_{T}^{-k}\mathcal{K}_{T}(\bb{h})\left\{\prod_{i}\Theta(J_{i}) + \oo(1)\right\} \\
&= \frac{\Xi_{T}(\bb{h})}{\lambda_{T}^{k}\{1 + \oo(1)\}}\left\{\prod_{i}\Theta(J_{i}) + \oo(1)\right\}
\geq c_{0}\Xi_{T}(\bb{h}) \geq c_{0}e^{-Ct} \geq e^{-C_{1}t}
\end{aligned}
\]
for all sufficiently large $T$. This denominator bound and the exceptional-probability estimate in \citet[Lemma~4.5, Eq.~(4.36)]{OuimetExpanding2026} give
\begin{equation}\label{eq:local.profile.conditional.truncation}
\begin{aligned}
\PP_{\mathrm{A}}((\mathcal{A}_{T,\mathrm{der}}^{\mathrm{A}})^{c}\mid E_{T}(\bb{h},\bb{J}))
&\leq \frac{\PP_{\mathrm{A}}((\mathcal{A}_{T,\mathrm{der}}^{\mathrm{A}})^{c})}{\PP_{\mathrm{A}}(E_{T}(\bb{h},\bb{J}))}
\leq \exp\{ - ct^{1 + 2\e_{0}} + (C + C_{1})t\} = \oo(1).
\end{aligned}
\end{equation}
Define the modified local paths $\widetilde{Z}_{T,h}^{\star}(v) \leqdef Z_{T,h}^{\star}(v)\ind_{\mathcal{A}_{T,\mathrm{der}}^{\star}}$. Since $|z|^{4} \leq C(e^{z} + e^{-z})$, \eqref{eq:local.profile.conditional.derivative.laplace} gives
\begin{equation}\label{eq:local.profile.conditional.derivative.fourth}
\EE_{\star}\left[|L_{T,h_{i}}^{\star}(v)|^{4}\ind_{\mathcal{A}_{T,\mathrm{der}}^{\star}}\,\middle|\,E_{T}(\bb{h},\bb{J})\right] \leq C.
\end{equation}
For $v \geq w$, the fundamental theorem of calculus and H\"older's inequality yield
\[
|\widetilde{Z}_{T,h_{i}}^{\star}(v) - \widetilde{Z}_{T,h_{i}}^{\star}(w)|^{4} = \left|\int_{w}^{v}L_{T,h_{i}}^{\star}(z)\rd{z}\right|^{4}\ind_{\mathcal{A}_{T,\mathrm{der}}^{\star}} \leq |v - w|^{3}\int_{w}^{v}|L_{T,h_{i}}^{\star}(z)|^{4}\ind_{\mathcal{A}_{T,\mathrm{der}}^{\star}}\rd{z}.
\]
Taking conditional expectations and using \eqref{eq:local.profile.conditional.derivative.fourth} proves
\begin{equation}\label{eq:local.profile.increment.tightness}
\EE_{\star}\left[|\widetilde{Z}_{T,h_{i}}^{\star}(v) - \widetilde{Z}_{T,h_{i}}^{\star}(w)|^{4}\,\middle|\,E_{T}(\bb{h},\bb{J})\right]
\leq C|v - w|^{4}.
\end{equation}
Since $\smash{\widetilde{Z}_{T,h_{i}}^{\star}}(0) = 0$, \eqref{eq:local.profile.increment.tightness} and Kolmogorov's tightness criterion in \citet[Theorem~23.7, p.~511]{Kallenberg2021} give tightness of the conditional laws of the modified paths $\smash{(\widetilde{Z}_{T,h_i}^{\star})_{i=1}^{k}}$ in $C([ - R,R])^{k}$, uniformly over the permitted roots $\bb{h}$, bins $\bb{J}$, and $\star \in \{\mathrm{A}, \mathrm{H}\}$; \eqref{eq:local.profile.conditional.truncation} transfers this to the original paths $(Z_{T,h_i}^{\star})_{i=1}^{k}$. Together with the conditional finite-dimensional convergence of Lemma~\ref{lem:local.profile.bins}, this gives uniform conditional weak convergence of $\smash{(Z_{T,h_i}^{\star})_{i=1}^{k}}$ to $(D_{\theta,i})_{i=1}^{k}$ in $C([ - R,R])^{k}$. Multiplying the resulting convergence of conditional expectations by the bin-probability asymptotic \eqref{eq:local.profile.unmarked.bin} of Lemma~\ref{lem:local.profile.bins} proves \eqref{eq:local.profile.marked.bin}.
\end{proof}

\begin{corollary}[A fixed lattice neighborhood]\label{cor:local.profile.lattice}
Fix $q > 0$ and $r \in \N$, and let $p(q,qr)$ denote the probability that the origin is a maximizer of the limiting local profile $D_{\theta}$ on the finite lattice $\{qj:j \in \{ - r,\ldots,r\}\}$:
\begin{equation}\label{eq:local.profile.finite.lattice.constant}
p(q,qr) \leqdef \PP(D_{\theta}(qj) \leq 0 \text{ for }j \in \mathcal{J}_{r}), \qquad \mathcal{J}_{r} \leqdef \{ - r,\ldots, - 1,1,\ldots,r\}.
\end{equation}
In the setting of Lemma~\ref{lem:local.profile}, with $R \geq qr$, as $T \to \infty$,
\begin{equation}\label{eq:local.profile.lattice.bin}
\begin{aligned}
&\PP_{\star}\left( E_{T}(\bb{h}, \bb{J}) \cap \left\{ F_{h_{i} + \ell_{T}qj}^{\star} \leq F_{h_{i}}^{\star} \text{ for all }i \in \{1,\ldots,k\}, \ j \in \mathcal{J}_{r} \right\} \right) \\
&\qquad= N_{T}^{-k}\mathcal{K}_{T}(\bb{h}) \left\{p(q,qr)^{k}\prod_{i}\Theta(J_{i}) + \oo(1)\right\},
\end{aligned}
\end{equation}
uniformly over the permitted roots and bins.
\end{corollary}

\begin{proof}[Proof of Corollary~\ref{cor:local.profile.lattice}]
By Lemma~\ref{lem:local.profile} and continuity of the finite evaluation maps, conditionally on $E_{T}(\bb{h},\bb{J})$,
\[
(Z_{T,h_{i}}^{\star}(qj))_{1\leq i\leq k,\,j\in\mathcal{J}_{r}} \xrightarrow{\mathrm{d}} (D_{\theta,i}(qj))_{1\leq i\leq k,\,j\in\mathcal{J}_{r}},
\]
uniformly over the permitted roots $\bb{h}$ and bins $\bb{J}$. Each limiting coordinate is Gaussian with variance $\frac{1}{2}\log(1 + q^{2}j^{2}) > 0$ by \eqref{eq:main.profile.covariance}, so it has no atom at zero. Hence the boundary of the nonpositive orthant has zero probability under the limiting law.

The portmanteau theorem and independence of the limiting profiles therefore give
\[
\begin{aligned}
&\PP_{\star}\left(Z_{T,h_{i}}^{\star}(qj) \leq 0 \text{ for all }i \in \{1,\ldots,k\},\ j \in \mathcal{J}_{r}\,\middle|\,E_{T}(\bb{h},\bb{J})\right) \\
&\qquad= \prod_{i=1}^{k}\PP\left(D_{\theta,i}(qj) \leq 0 \text{ for all }j \in \mathcal{J}_{r}\right) + \oo(1) = p(q,qr)^{k} + \oo(1).
\end{aligned}
\]
The error is uniform because every sequence of permitted roots $\bb{h}$ and bins $\bb{J}$ has the same conditional weak limit and the same continuity set.

Since $Z_{T,h_{i}}^{\star}(qj) = F_{h_{i} + \ell_{T}qj}^{\star} - F_{h_{i}}^{\star}$, multiplying this conditional probability by the bin-probability asymptotic \eqref{eq:local.profile.unmarked.bin} of Lemma~\ref{lem:local.profile.bins} proves \eqref{eq:local.profile.lattice.bin}.
\end{proof}

\subsection{Conditional localization of high excursions}\label{sec:localization}

The next three lemmas control high excursions of the Euler field $F^{\mathrm{A}}$ and the Haar field $F^{\mathrm{H}}$ away from a prescribed high point, including under simultaneous conditioning at finitely many separated high points. Here, an annulus around a conditioned high point $h_{i}$ means the positions $h$ satisfying $d \leq |h - h_{i}| \leq 2d$, restricted to the block containing $h_{i}$. In one dimension, this set consists of at most two intervals. The first lemma establishes the covariance and conditional mean-shift estimates needed on these sets. The second bounds the exponential moment of the supremum of $F_{h}^{\star} - F_{h_{i}}^{\star}$ over one such set. The third combines these bounds over successive distance scales to show that further high excursions within the same block become unlikely far from $h_{i}$ relative to the microscopic scale $\ell_{T}$.

Throughout this subsection, fix $\theta > 0$, $0 < \e < 1/4$, a compact interval $J \subseteq \R$, $b > 0$, and $k \in \N$. For $s \in [0,1]$ and $j \in \Z$, let $\smash{B_{j}^{\e,s}} = [s + j + \e,s + j + 1 - \e] \cap I_{\theta}$, omitting empty intervals. Choose $h_{1},\ldots,h_{k}$ in distinct intervals $\smash{B_{j}^{\e,s}}$, write $h_{i} \in \smash{B_{j_{i}}^{\e,s}}$, and choose $x_{1},\ldots,x_{k}$ such that $[x_{i},x_{i} + b] \subseteq J$. For $\star \in \{\mathrm{A}, \mathrm{H}\}$, define
\[
E_{\bb{h}}^{\star} \leqdef \bigcap_{i=1}^{k}\{F_{h_{i}}^{\star} - u_{T} \in [x_{i},x_{i} + b]\}.
\]
Choose $m_{-} \in (1,\sqrt{1 + \theta})$ depending only on $\theta$. Write
\[
\bb{F}^{\mathrm{H}} \leqdef (F_{h_{1}}^{\mathrm{H}},\ldots,F_{h_{k}}^{\mathrm{H}})^{\mathsf{T}}, \qquad \bb{v} \leqdef (u_{T} + x_{1},\ldots,u_{T} + x_{k})^{\mathsf{T}}, \qquad \Sigma_{T} \leqdef \Cov_{\mathrm{H}}(\bb{F}^{\mathrm{H}}).
\]
Define the conditioning tilt by
\[
\bb{\lambda}_{T} \leqdef \Sigma_{T}^{-1}\bb{v},
\]
for all sufficiently large $T$, when $\Sigma_{T}$ is invertible by Lemma~\ref{lem:annular.covariance} below.

For $d > 0$, let $\smash{F_{d}^{\star,\mathrm{lo}}}(h)$ be the sum of the terms in $F_{h}^{\star}$ with $\log(p^{j}) \leq 1/d$. Define the dyadic radii by
\begin{equation}\label{eq:sup.localization.radii}
d_{n} \leqdef \ell_{T}2^{n}, \qquad n \in \{0,\ldots,\lfloor s_{T}/\log 2\rfloor\}.
\end{equation}
For $i \in \{1,\ldots,k\}$ and $d = d_{n} \leq 1$, write
\[
r \leqdef \log(d/\ell_{T}) = n\log 2, \qquad A_{i,d} \leqdef \{h \in B_{j_{i}}^{\e,s}:d \leq |h - h_{i}| \leq 2d\}.
\]
Since $\ell_{T} \leq d \leq 1$ and $s_{T} = \log(1/\ell_{T})$, the scale parameter satisfies $0 \leq r \leq s_{T}$. With $i$ and $d$ fixed, separate the high-frequency root value from the common low-frequency value by setting
\[
X^{\star} \leqdef F_{h_{i}}^{\star} - F_{d}^{\star,\mathrm{lo}}(h_{i}), \qquad Y^{\star}(h) \leqdef F_{h}^{\star} - F_{d}^{\star,\mathrm{lo}}(h_{i}), \qquad D^{\star}(z) \leqdef \ell_{T}(F^{\star})'(z),
\]
for $h \in A_{i,d}$ and $z \in I_{\theta}$. The exact increment identity is
\begin{equation}\label{eq:sup.localization.increment}
F_{h}^{\star} - F_{h_{i}}^{\star} = Y^{\star}(h) - X^{\star}.
\end{equation}
Define the Gaussian regression mean terms associated with the conditioned root heights $\bb{v}$ by
\begin{equation}\label{eq:sup.localization.mean.terms}
\mu_{X} \leqdef \Cov_{\mathrm{H}}(\bb{F}^{\mathrm{H}},X^{\mathrm{H}})^{\mathsf{T}}\bb{\lambda}_{T}, \qquad \mu_{Y}(h) \leqdef \Cov_{\mathrm{H}}(\bb{F}^{\mathrm{H}},Y^{\mathrm{H}}(h))^{\mathsf{T}}\bb{\lambda}_{T}.
\end{equation}
Since $\bb{\lambda}_{T} = \Sigma_{T}^{-1}\bb{v}$, these expressions would be the conditional means of $X^{\mathrm{H}}$ and $Y^{\mathrm{H}}(h)$ given $\smash{\bb{F}^{\mathrm{H}}} = \bb{v}$ if the variables were jointly centered Gaussian. For the Euler and Haar fields considered here, they enter the conditional exponential-moment bound \eqref{eq:sup.localization.conditional.laplace}.

\begin{lemma}[Annular covariance and conditioning tilt estimates]\label{lem:annular.covariance}
Under the preceding assumptions, the eigenvalues of $\Sigma_{T}$ are comparable to $t$. For all sufficiently large $T$, the conditioning tilt satisfies
\begin{equation}\label{eq:sup.localization.tilt.bounds}
\min_{j\in\{1,\ldots,k\}}\lambda_{T,j} \geq 2m_{-}, \qquad \max_{j\in\{1,\ldots,k\}}|\lambda_{T,j}| \leq C_{k,\e,J}.
\end{equation}
For $h \in A_{i,d}$, the annular covariances satisfy
\begin{equation}\label{eq:sup.localization.local.covariances}
\begin{aligned}
&\Var_{\mathrm{H}}(X^{\mathrm{H}}) = \frac{r}{2} + \OO(1), \qquad \Var_{\mathrm{H}}(Y^{\mathrm{H}}(h)) = \frac{r}{2} + \OO(1), \\
&\Cov_{\mathrm{H}}(F_{h_{i}}^{\mathrm{H}},X^{\mathrm{H}}) = \frac{r}{2} + \OO(1), \qquad |\Cov_{\mathrm{H}}(X^{\mathrm{H}},Y^{\mathrm{H}}(h))| + |\Cov_{\mathrm{H}}(F_{h_{i}}^{\mathrm{H}},Y^{\mathrm{H}}(h))| \leq C.
\end{aligned}
\end{equation}
The scaled derivative satisfies, for $z \in I_{\theta}$,
\begin{equation}\label{eq:sup.localization.derivative.covariances}
\begin{aligned}
\Var_{\mathrm{H}}(D^{\mathrm{H}}(z)) &\leq C, \\
|\Cov_{\mathrm{H}}(D^{\mathrm{H}}(z),X^{\mathrm{H}})| + |\Cov_{\mathrm{H}}(D^{\mathrm{H}}(z),Y^{\mathrm{H}}(h))| + \sum_{j=1}^{k}|\Cov_{\mathrm{H}}(D^{\mathrm{H}}(z),F_{h_{j}}^{\mathrm{H}})| &\leq C_{k}.
\end{aligned}
\end{equation}
If $d \leq t^{-3}$, then
\begin{equation}\label{eq:sup.localization.small.radius.means}
\mu_{X} = \frac{\lambda_{T,i}r}{2} + \OO_{k,\e,J}(1) \geq m_{-}r - C_{k,\e,J}, \qquad |\mu_{Y}(h)| \leq C_{k,\e,J}.
\end{equation}
If $d > t^{-3}$, then
\begin{equation}\label{eq:sup.localization.large.radius.means}
\mu_{X} \geq m_{-}r - C_{k,\e,J}\log t, \qquad |\mu_{Y}(h)| \leq C_{k,\e,J}\log t.
\end{equation}
The derivative conditioning tilt term satisfies
\begin{equation}\label{eq:sup.localization.derivative.mean}
\left|\Cov_{\mathrm{H}}(\bb{F}^{\mathrm{H}},D^{\mathrm{H}}(z))^{\mathsf{T}}\bb{\lambda}_{T}\right| \leq C_{k,\e,J}.
\end{equation}
All estimates are uniform in $s \in [0,1]$, the representative positions $h_{1},\ldots,h_{k}$ belonging to distinct blocks $B_{j}^{\e,s}$, the allowed bin endpoints $x_{1},\ldots,x_{k}$, $i \in \{1,\ldots,k\}$, the dyadic radii $d \leq 1$, and the indicated positions $h$ and $z$.
\end{lemma}

\begin{proof}[Proof of Lemma~\ref{lem:annular.covariance}]
The proof is divided into three steps. Step~1 estimates the root covariance matrix and the conditioning tilt. Step~2 computes the covariances on one annulus, including those of the scaled derivative. Step~3 controls the contributions of the other conditioned roots to the conditioning tilt terms.

\medskip
{\bf Step 1. The root covariance matrix and conditioning tilt.} Each block $\smash{B_{j}^{\e,s}} = [s + j + \e,s + j + 1 - \e] \cap I_{\theta}$ is obtained by trimming a segment of length $\e$ from each end of the unit interval $[s + j,s + j + 1]$ and then intersecting with $I_{\theta}$. Consequently, distinct nonempty blocks are separated by at least $2\e$. The covariance estimates in \citet[Lemma~4.3, Eqs.~(4.13)--(4.14)]{OuimetExpanding2026}, together with $V_{T} \asymp t$ from \eqref{eq:field.variance} and \eqref{eq:scale.asymptotics}, therefore give
\begin{equation}\label{eq:sup.localization.root.matrix}
\Sigma_{T} = V_{T}\mathrm{Id}_{k} + \OO_{k,\e}(\log t), \qquad \Sigma_{T}^{-1} = V_{T}^{-1}\mathrm{Id}_{k} + \OO_{k,\e}\left(\frac{\log t}{t^{2}}\right).
\end{equation}
Here and below, the matrix errors are in operator norm. Since \eqref{eq:height.asymptotics} and \eqref{eq:scale.asymptotics} give $u_{T} = \OO(t)$ and the bin endpoints $x_{1},\ldots,x_{k}$ are uniformly bounded because they belong to the fixed compact interval $J$, \eqref{eq:sup.localization.root.matrix} gives the conditioning tilt
\begin{equation}\label{eq:sup.localization.conditioning.tilt}
\bb{\lambda}_{T} = \frac{u_{T}}{V_{T}}\bb{1}_{k} + \OO_{k,\e,J}\left(\frac{\log t}{t}\right).
\end{equation}
The first estimate in \eqref{eq:sup.localization.root.matrix} and $V_{T} \asymp t$ show that the root covariance eigenvalues are comparable to $t$. Equations \eqref{eq:height.asymptotics} and \eqref{eq:sup.localization.conditioning.tilt} then prove \eqref{eq:sup.localization.tilt.bounds} for all sufficiently large $T$.

\medskip
{\bf Step 2. Covariances on one annulus.} The frequency decomposition gives the estimates in \eqref{eq:sup.localization.local.covariances} as follows.

First, orthogonality of distinct prime powers under the Haar law and partial summation give
\begin{equation}\label{eq:sup.localization.root.variance}
\begin{aligned}
\Var_{\mathrm{H}}(X^{\mathrm{H}})
&= \frac{1}{2}\sum_{\substack{p^{j}\leq X \\ j\log p>1/d}}\frac{1}{j^{2}p^{2j\sig}}
= \frac{1}{2}\int_{1/d}^{\log X}e^{-\ell_{T}v}\frac{\rd{v}}{v} + \OO(1) \\
&= \frac{1}{2}\int_{\ell_{T}/d}^{1}\frac{\rd{y}}{y} + \frac{1}{2}\int_{\ell_{T}/d}^{1}\frac{e^{-y} - 1}{y}\rd{y} + \frac{1}{2}\int_{1}^{\ell_{T}\log X}\frac{e^{-y}}{y}\rd{y} + \OO(1)
= \frac{r}{2} + \OO(1).
\end{aligned}
\end{equation}
The second line uses the substitution $y = \ell_{T}v$. For the last equality, the first integral equals $\log(d/\ell_{T}) = r$, the second has absolute value at most $1$ because $|e^{-y} - 1| \leq y$ for $y \geq 0$, and the third is bounded by $\int_{1}^{\infty}e^{-y}\rd{y} = e^{-1}$. The higher-prime-power variance is bounded by $\smash{\sum_{p}\sum_{j\geq2}(j^{2}p^{j})^{-1}} < \infty$. This proves the first variance estimate in \eqref{eq:sup.localization.local.covariances}.

To obtain the other covariances, define the high-frequency Haar field $H_{d}^{\mathrm{H}}(h) \leqdef F_{h}^{\mathrm{H}} - \smash{F_{d}^{\mathrm{H},\mathrm{lo}}}(h)$ and the low-frequency Haar increment $\Delta_{d}^{\mathrm{H}}(h) \leqdef \smash{F_{d}^{\mathrm{H},\mathrm{lo}}}(h) - \smash{F_{d}^{\mathrm{H},\mathrm{lo}}}(h_{i})$. Then
\[
X^{\mathrm{H}} = H_{d}^{\mathrm{H}}(h_{i}), \qquad Y^{\mathrm{H}}(h) = H_{d}^{\mathrm{H}}(h) + \Delta_{d}^{\mathrm{H}}(h), \qquad \Cov_{\mathrm{H}}(H_{d}^{\mathrm{H}}(z),F_{d}^{\mathrm{H},\mathrm{lo}}(w)) = 0.
\]
For $\xi \leqdef h - h_{i}$ with $d \leq |\xi| \leq 2d$, integration by parts gives
\[
\left|\int_{1/d}^{\log X}e^{-\ell_{T}v}\cos(\xi v)\frac{\rd{v}}{v}\right| \leq \frac{2d}{|\xi|} \leq 2.
\]
The corresponding prime-sum estimate, including the bounded higher-power contribution, is therefore
\begin{equation}\label{eq:sup.localization.high.frequency.covariance}
|\Cov_{\mathrm{H}}(H_{d}^{\mathrm{H}}(h_{i}),H_{d}^{\mathrm{H}}(h))| \leq C.
\end{equation}

For the low-frequency increment $\Delta_{d}^{\mathrm{H}}(h)$, use $1 - \cos y \leq y^{2}/2$ and $|\xi| \leq 2d$ to obtain
\begin{equation}\label{eq:sup.localization.low.frequency.increment}
\Var_{\mathrm{H}}(\Delta_{d}^{\mathrm{H}}(h))
\leq Cd^{2}\sum_{\substack{p\leq X \\ \log p\leq1/d}}\frac{(\log p)^{2}}{p} + C
\leq Cd^{2}\int_{0}^{1/d}v\rd{v} + C \leq C.
\end{equation}
Stationarity of the low- and high-frequency Haar fields and orthogonality between these two parts, together with \eqref{eq:sup.localization.root.variance} and \eqref{eq:sup.localization.low.frequency.increment}, now give the individual identities
\[
\begin{aligned}
&\Var_{\mathrm{H}}(Y^{\mathrm{H}}(h)) = \Var_{\mathrm{H}}(H_{d}^{\mathrm{H}}(h)) + \Var_{\mathrm{H}}(\Delta_{d}^{\mathrm{H}}(h)) = \frac{r}{2} + \OO(1), \\
&\Cov_{\mathrm{H}}(F_{h_{i}}^{\mathrm{H}},X^{\mathrm{H}}) = \Var_{\mathrm{H}}(X^{\mathrm{H}}) = \frac{r}{2} + \OO(1), \\
&\Cov_{\mathrm{H}}(X^{\mathrm{H}},Y^{\mathrm{H}}(h)) = \Cov_{\mathrm{H}}(H_{d}^{\mathrm{H}}(h_{i}),H_{d}^{\mathrm{H}}(h)), \\
&\Cov_{\mathrm{H}}(F_{h_{i}}^{\mathrm{H}},Y^{\mathrm{H}}(h)) = \Cov_{\mathrm{H}}(H_{d}^{\mathrm{H}}(h_{i}),H_{d}^{\mathrm{H}}(h)) - \frac{1}{2}\Var_{\mathrm{H}}(\Delta_{d}^{\mathrm{H}}(h)).
\end{aligned}
\]
These identities, \eqref{eq:sup.localization.high.frequency.covariance}, and \eqref{eq:sup.localization.low.frequency.increment} prove \eqref{eq:sup.localization.local.covariances} uniformly for $h \in A_{i,d}$ and $0 \leq r \leq s_{T}$.

The scaled derivative $D^{\mathrm{H}}(z)$ has bounded variance and bounded covariances with the coordinates just introduced. Indeed, for $Q^{\mathrm{H}}$ equal to $X^{\mathrm{H}}$, $Y^{\mathrm{H}}(h)$, or any root height $F_{h_{j}}^{\mathrm{H}}$, orthogonality of the prime powers and absolute-value bounds give
\[
\begin{aligned}
\Var_{\mathrm{H}}(D^{\mathrm{H}}(z))
&= \frac{\ell_{T}^{2}}{2}\sum_{p^{j}\leq X}\frac{(\log p)^{2}}{p^{2j\sig}}
\leq C\ell_{T}^{2}\int_{0}^{\log X}ve^{-\ell_{T}v}\rd{v} + C, \\
|\Cov_{\mathrm{H}}(D^{\mathrm{H}}(z),Q^{\mathrm{H}})|
&\leq C\ell_{T}\sum_{p^{j}\leq X}\frac{\log p}{jp^{2j\sig}}
\leq C\ell_{T}\int_{0}^{\log X}e^{-\ell_{T}v}\rd{v} + C.
\end{aligned}
\]
The integral bounds follow by partial summation for the terms with $j = 1$, while the terms with $j \geq 2$ contribute bounded amounts. The substitution $y = \ell_{T}v$ gives
\[
\ell_{T}^{2}\int_{0}^{\log X}ve^{-\ell_{T}v}\rd{v}
\leq \int_{0}^{\infty}ye^{-y}\rd{y} = 1, \qquad
\ell_{T}\int_{0}^{\log X}e^{-\ell_{T}v}\rd{v}
\leq \int_{0}^{\infty}e^{-y}\rd{y} = 1.
\]
These estimates prove \eqref{eq:sup.localization.derivative.covariances}, uniformly for $h \in A_{i,d}$ and $z \in I_{\theta}$.

\medskip
{\bf Step 3. The other conditioned roots.} If $d \leq t^{-3}$ and $j \neq i$, let $z \in \{h_{i},h\}$ and write $\xi = z - h_{j}$. Since $z$ and $h_{j}$ lie in distinct blocks, $|\xi| \geq 2\e$. Partial summation gives
\[
\Cov_{\mathrm{H}}(H_{d}^{\mathrm{H}}(z),F_{h_{j}}^{\mathrm{H}})
= \frac{1}{2}\int_{1/d}^{\log X}e^{-\ell_{T}v}\cos(\xi v)\frac{\rd{v}}{v} + \OO(1).
\]
Integration by parts, using that $v \mapsto e^{-\ell_{T}v}/v$ is decreasing, yields
\[
\begin{aligned}
&
\left|\int_{1/d}^{\log X}e^{-\ell_{T}v}\cos(\xi v)\frac{\rd{v}}{v}\right| \\
&= \left|\left[\frac{e^{-\ell_{T}v}\sin(\xi v)}{\xi v}\right]_{1/d}^{\log X}
+ \frac{1}{\xi}\int_{1/d}^{\log X}e^{-\ell_{T}v}\left(\frac{\ell_{T}}{v} + \frac{1}{v^{2}}\right)\sin(\xi v)\rd{v}\right| \\
&\leq \frac{1}{|\xi|}\left\{d e^{-\ell_{T}/d} + \frac{e^{-\ell_{T}\log X}}{\log X} + \int_{1/d}^{\log X}e^{-\ell_{T}v}\left(\frac{\ell_{T}}{v} + \frac{1}{v^{2}}\right)\rd{v}\right\}
= \frac{2d e^{-\ell_{T}/d}}{|\xi|}
\leq \frac{d}{\e}.
\end{aligned}
\]
The prime number theorem error is $\OO((1 + |I_{\theta}|)t^{C}e^{-c_{2}t^{3/2}}) = \oo(1)$, since $1/d \geq t^{3}$. The higher-power contribution is bounded by $\smash{\frac{1}{2}\sum_{p}\sum_{a\geq2}(a^{2}p^{a})^{-1}} < \infty$. Thus both $\Cov_{\mathrm{H}}(X^{\mathrm{H}},\allowbreak F_{h_{j}}^{\mathrm{H}})$ and $\Cov_{\mathrm{H}}(H_{d}^{\mathrm{H}}(h),F_{h_{j}}^{\mathrm{H}})$ are uniformly bounded.

The low-frequency increment has bounded covariance with every root height $\smash{F_{h_{j}}^{\mathrm{H}}}$, $j \in \{1,\ldots,k\}$, as well. Indeed, $|e^{-\ii a(h - h_{i})\log p} - 1| \leq 2da\log p$ implies
\[
|\Cov_{\mathrm{H}}(F_{h_{j}}^{\mathrm{H}},\Delta_{d}^{\mathrm{H}}(h))|
\leq Cd\sum_{\log p\leq1/d}\frac{\log p}{p} + C
\leq Cd\int_{0}^{1/d}\rd{v} + C
\leq C.
\]
For $d \leq t^{-3}$, substituting $\Cov_{\mathrm{H}}(F_{h_{i}}^{\mathrm{H}},X^{\mathrm{H}}) = r/2 + \OO(1)$ and $|\Cov_{\mathrm{H}}(F_{h_{i}}^{\mathrm{H}},Y^{\mathrm{H}}(h))| \leq C$ from \eqref{eq:sup.localization.local.covariances} into the definitions \eqref{eq:sup.localization.mean.terms}, and using the tilt bounds \eqref{eq:sup.localization.tilt.bounds} and the bounded covariances with the other root heights $F_{h_{j}}^{\mathrm{H}}$, $j \neq i$, gives
\[
\mu_{X} = \frac{\lambda_{T,i}r}{2} + \OO_{k,\e,J}(1), \qquad |\mu_{Y}(h)| \leq C_{k,\e,J}.
\]
The lower bound $\lambda_{T,i} \geq 2m_{-}$ from \eqref{eq:sup.localization.tilt.bounds} proves \eqref{eq:sup.localization.small.radius.means}.

For $d > t^{-3}$, taking absolute values in the prime-power covariance sum gives
\[
|\Cov_{\mathrm{H}}(F_{h_{j}}^{\mathrm{H}},F_{d}^{\mathrm{H},\mathrm{lo}}(z))|
\leq \frac{1}{2}\sum_{\log p\leq1/d}\frac{1}{p} + C
\leq C\int_{\log 2}^{1/d}\frac{\rd{v}}{v} + C
\leq C\log(2 + 1/d) \leq C\log t.
\]
Both $h$ and $h_{i}$ are at distance at least $2\e$ from the other roots. For $j \neq i$, the definitions of $X^{\mathrm{H}}$ and $Y^{\mathrm{H}}(h)$ give
\[
\begin{aligned}
\Cov_{\mathrm{H}}(F_{h_{j}}^{\mathrm{H}},X^{\mathrm{H}})
&= \Cov_{\mathrm{H}}(F_{h_{j}}^{\mathrm{H}},F_{h_{i}}^{\mathrm{H}}) - \Cov_{\mathrm{H}}(F_{h_{j}}^{\mathrm{H}},F_{d}^{\mathrm{H},\mathrm{lo}}(h_{i}))
= \OO_{\e}(\log t), \\
\Cov_{\mathrm{H}}(F_{h_{j}}^{\mathrm{H}},Y^{\mathrm{H}}(h))
&= \Cov_{\mathrm{H}}(F_{h_{j}}^{\mathrm{H}},F_{h}^{\mathrm{H}}) - \Cov_{\mathrm{H}}(F_{h_{j}}^{\mathrm{H}},F_{d}^{\mathrm{H},\mathrm{lo}}(h_{i}))
= \OO_{\e}(\log t).
\end{aligned}
\]
In each line, the full-field covariance is $\OO_{\e}(\log t)$ by spatial separation, and the low-frequency covariance is $\OO(\log t)$ by the preceding bound. Substituting into \eqref{eq:sup.localization.mean.terms}, and using \eqref{eq:sup.localization.local.covariances} for the contribution of the root $h_{i}$ and \eqref{eq:sup.localization.tilt.bounds} for the tilt coefficients, gives
\[
\mu_{X} = \frac{\lambda_{T,i}r}{2} + \OO_{k,\e,J}(\log t) \geq m_{-}r - C_{k,\e,J}\log t, \qquad
|\mu_{Y}(h)| \leq C_{k,\e,J}\log t.
\]
This proves \eqref{eq:sup.localization.large.radius.means}.

For the derivative conditioning tilt term, \eqref{eq:sup.localization.derivative.covariances} and \eqref{eq:sup.localization.tilt.bounds} give
\[
\left|\Cov_{\mathrm{H}}(\bb{F}^{\mathrm{H}},D^{\mathrm{H}}(z))^{\mathsf{T}}\bb{\lambda}_{T}\right| \leq \max_{j\in\{1,\ldots,k\}}|\lambda_{T,j}| \sum_{j=1}^{k}|\Cov_{\mathrm{H}}(F_{h_{j}}^{\mathrm{H}},D^{\mathrm{H}}(z))| \leq C_{k,\e,J},
\]
which proves \eqref{eq:sup.localization.derivative.mean}.
\end{proof}

With the fixed $m_{-} \in (1,\sqrt{1 + \theta})$ chosen above, choose $a$ depending only on $\theta$ so that
\begin{equation}\label{eq:sup.localization.exponents}
0 < a < \min\{2,2(m_{-} - 1)\}, \qquad c_{0} \leqdef a(m_{-} - 1) - \frac{a^{2}}{4} > 0.
\end{equation}
The term $am_{-}$ will come from the conditioned root height, the term $a$ from the combined contribution of the spatial cell count and the exponential-moment bound for each cell, and the term $a^{2}/4$ from the remaining fluctuation of that height.

Fix $0 < \e_{0} < 1/4$. Define a common event for all annuli by
\begin{equation}\label{eq:sup.localization.common.event}
\mathcal{A}_{T,\mathrm{sp}}^{\mathrm{A}} \leqdef \left\{\sup_{h\in I_{\theta}}\max\{|F_{h}^{\mathrm{A}}|,|\ell_{T}(F^{\mathrm{A}})'(h)|\} \leq t^{1 + \e_{0}}\right\} \cap \bigcap_{n=0}^{\lfloor s_{T}/\log 2\rfloor} \left\{\sup_{h\in I_{\theta}}|F_{d_{n}}^{\mathrm{A},\mathrm{lo}}(h)| \leq t^{1 + \e_{0}}\right\}.
\end{equation}
For the Haar field, let $\mathcal{A}_{T,\mathrm{sp}}^{\mathrm{H}}$ be the whole probability space, so that $\ind_{\mathcal{A}_{T,\mathrm{sp}}^{\mathrm{H}}} \equiv 1$.

\begin{lemma}[Conditional exponential bound for an annular supremum]\label{lem:annular.supremum}
Under the preceding assumptions, let $c \leqdef c_{0}/2 > 0$. For all sufficiently large $T$ and $\star \in \{\mathrm{A}, \mathrm{H}\}$,
\begin{equation}\label{eq:sup.localization.uniform.annulus}
\EE_{\star}\left[\exp\left\{a\sup_{h\in A_{i,d}}(F_{h}^{\star} - F_{h_{i}}^{\star})\right\}\ind_{\mathcal{A}_{T,\mathrm{sp}}^{\star}}\mid E_{\bb{h}}^{\star}\right] \leq Ce^{-cr}.
\end{equation}
If $A_{i,d} = \varnothing$, the supremum in \eqref{eq:sup.localization.uniform.annulus} is taken to be $ - \infty$, so the left-hand side is zero. The common event satisfies, as $T \to \infty$,
\begin{equation}\label{eq:sup.localization.conditional.event}
\PP_{\star}((\mathcal{A}_{T,\mathrm{sp}}^{\star})^{c}\mid E_{\bb{h}}^{\star}) = \oo(1).
\end{equation}
Both estimates are uniform in $s \in [0,1]$, the representative positions $h_{1},\ldots,h_{k}$ in distinct blocks $B_{j_{1}}^{\e,s},\ldots,B_{j_{k}}^{\e,s}$, the allowed bin endpoints $x_{1},\ldots,x_{k}$, $i \in \{1,\ldots,k\}$, and the dyadic radii $d \leq 1$. The constants $a,c > 0$ depend only on $\theta$, and the constant $C$ may also depend on $k,\e,J,b$.
\end{lemma}

\begin{proof}[Proof of Lemma~\ref{lem:annular.supremum}]
The proof is divided into three steps. Step~1 bounds the probability of $(\mathcal{A}_{T,\mathrm{sp}}^{\mathrm{A}})^{c}$ and establishes the conditional Laplace estimate with expectations restricted to $\mathcal{A}_{T,\mathrm{sp}}^{\star}$. Step~2 partitions $A_{i,d}$ into intervals of length at most $\ell_{T}$, uses the derivative to control the field between points of each interval, and combines Lemma~\ref{lem:annular.covariance} with H\"older's inequality to bound the exponential moment of the supremum over the entire annulus. Step~3 bounds the conditional probability of the exceptional event $(\mathcal{A}_{T,\mathrm{sp}}^{\star})^{c}$ given the root-height event $E_{\bb{h}}^{\star}$.

\medskip
{\bf Step 1. The common event $\smash{\mathcal{A}_{T,\mathrm{sp}}^{\star}}$ and conditional Laplace estimate.} By \eqref{eq:scale.asymptotics} and \eqref{eq:sup.localization.radii}, there are $\OO(t)$ low-frequency sums $\smash{F_{d_{n}}^{\mathrm{A},\mathrm{lo}}}$ in \eqref{eq:sup.localization.common.event}. The coefficient and derivative assumptions of \citet[Lemma~4.5]{OuimetExpanding2026} are satisfied because, for $q \in \{0,1,2\}$,
\begin{equation}\label{eq:sup.localization.coefficients}
\left|\ell_{T}^{q}\frac{\partial^{q}}{\partial h^{q}}\left(\frac{p^{-\ii jh}}{jp^{j\sig}}\right)\right|
= \frac{(\ell_{T}j\log p)^{q}}{jp^{j\sig}}
= \frac{2^{q}}{j}(\eta_{\sig}j\log p)^{q}e^{-\eta_{\sig}j\log p}p^{-j/2}
\leq C_{q}p^{-j/2}.
\end{equation}
The cases $q = 0,1$ verify these assumptions for the Euler field $F^{\mathrm{A}}$, while $q = 1,2$ verify them for its scaled derivative $\ell_{T}(F^{\mathrm{A}})'$. The low-frequency sum $\smash{F_{d_{n}}^{\mathrm{A},\mathrm{lo}}}(h)$ contains only the terms with $j\log p \leq 1/d_{n}$. This restriction replaces the excluded coefficients by zero and therefore preserves the bound \eqref{eq:sup.localization.coefficients}. Applying \citet[Lemma~4.5, Eq.~(4.36)]{OuimetExpanding2026} simultaneously to the field $F^{\mathrm{A}}$, its scaled derivative $\ell_{T}(F^{\mathrm{A}})'$, and all the low-frequency sums $\smash{F_{d_{n}}^{\mathrm{A},\mathrm{lo}}}$ indexed by \eqref{eq:sup.localization.radii}, which together define the event $\mathcal{A}_{T,\mathrm{sp}}^{\mathrm{A}}$ in \eqref{eq:sup.localization.common.event}, gives
\begin{equation}\label{eq:sup.localization.exceptional}
\PP_{\mathrm{A}}((\mathcal{A}_{T,\mathrm{sp}}^{\mathrm{A}})^{c}) \leq \exp\{ - c_{1}t^{1 + 2\e_{0}} + C_{1}t\}.
\end{equation}
Since the event $\mathcal{A}_{T,\mathrm{sp}}^{\mathrm{A}}$ is common to all annuli, the probability of its complement needs to be accounted for only once, after summing the annular bounds.

Recall that, for $h \in A_{i,d}$ and $z \in I_{\theta}$,
\[
X^{\star} = F_{h_{i}}^{\star} - F_{d}^{\star,\mathrm{lo}}(h_{i}), \qquad
Y^{\star}(h) = F_{h}^{\star} - F_{d}^{\star,\mathrm{lo}}(h_{i}), \qquad
D^{\star}(z) = \ell_{T}(F^{\star})'(z).
\]
With $a$ chosen as in \eqref{eq:sup.localization.exponents}, let $Z^{\star}$ be one of $ - aX^{\star}$, $ - aX^{\star} + 2Y^{\star}(h)$, or $ - aX^{\star} + 2Y^{\star}(h) \pm 2D^{\star}(z)$, and define $\bb{c}_{Z} \leqdef \Cov_{\mathrm{H}}(\bb{F}^{\mathrm{H}},Z^{\mathrm{H}})$. These linear combinations inherit the coefficient bounds \eqref{eq:sup.localization.coefficients}, and their Euler versions are bounded by $Ct^{1 + \e_{0}}$ on the common event \eqref{eq:sup.localization.common.event}. Moreover, the eigenvalues of $\Sigma_{T}$ are comparable to $t$ by Lemma~\ref{lem:annular.covariance}, and the conditioning bins have fixed width $b > 0$. Applying the conditional Laplace estimate in \citet[Lemma~4.6, Eq.~(4.51)]{OuimetExpanding2026} therefore gives
\begin{equation}\label{eq:sup.localization.conditional.laplace}
\EE_{\star}[e^{Z^{\star}}\ind_{\mathcal{A}_{T,\mathrm{sp}}^{\star}}\mid E_{\bb{h}}^{\star}] \leq C\exp\left\{\bb{c}_{Z}^{\mathsf{T}}\bb{\lambda}_{T} + \frac{1}{2}\Var_{\mathrm{H}}(Z^{\mathrm{H}})\right\}.
\end{equation}
Here, the nonpositive term $ - \frac{1}{2}\bb{c}_{Z}^{\mathsf{T}}\Sigma_{T}^{-1}\bb{c}_{Z}$ in the exponent of that estimate has been omitted.

\medskip
{\bf Step 2. An exponential bound for a continuous annular supremum.} Set
\begin{equation}\label{eq:sup.localization.weights}
\Gamma_{a}(r,m_{-}) \leqdef \exp\left\{ - am_{-}r + \frac{a^{2}r}{4}\right\}, \qquad \Xi_{T}^{\mathrm{ann}}(d) \leqdef \begin{cases}
1, & d \leq t^{-3}, \\
\exp(C_{k,\e,J}\log t), & d > t^{-3}.
\end{cases}
\end{equation}
For $Z^{\star} = - aX^{\star}$, the variance estimate \eqref{eq:sup.localization.local.covariances} and the mean bounds \eqref{eq:sup.localization.small.radius.means} and \eqref{eq:sup.localization.large.radius.means}, all from Lemma~\ref{lem:annular.covariance}, bound the exponent in \eqref{eq:sup.localization.conditional.laplace} by
\[
 - a\mu_{X} + \frac{a^{2}}{2}\Var_{\mathrm{H}}(X^{\mathrm{H}}) \leq - am_{-}r + \frac{a^{2}r}{4} + \log \Xi_{T}^{\mathrm{ann}}(d) + C.
\]
For $Z^{\star} = - aX^{\star} + 2Y^{\star}(h) \pm 2D^{\star}(z)$, all cross-covariances in its variance are bounded by \eqref{eq:sup.localization.local.covariances} and \eqref{eq:sup.localization.derivative.covariances} of Lemma~\ref{lem:annular.covariance}. These estimates, the mean bounds \eqref{eq:sup.localization.small.radius.means} and \eqref{eq:sup.localization.large.radius.means}, and \eqref{eq:sup.localization.derivative.mean}, all from Lemma~\ref{lem:annular.covariance}, give
\[
\begin{aligned}
\frac{1}{2}\Var_{\mathrm{H}}(Z^{\mathrm{H}})
&= \frac{a^{2}}{2}\Var_{\mathrm{H}}(X^{\mathrm{H}}) + 2\Var_{\mathrm{H}}(Y^{\mathrm{H}}(h)) + \OO(1)
= \frac{a^{2}r}{4} + r + \OO(1), \\
\Cov_{\mathrm{H}}(\bb{F}^{\mathrm{H}},Z^{\mathrm{H}})^{\mathsf{T}}\bb{\lambda}_{T}
&= - a\mu_{X} + 2\mu_{Y}(h) \pm 2\Cov_{\mathrm{H}}(\bb{F}^{\mathrm{H}},D^{\mathrm{H}}(z))^{\mathsf{T}}\bb{\lambda}_{T} \\
&\leq - am_{-}r + \log \Xi_{T}^{\mathrm{ann}}(d) + C,
\end{aligned}
\]
where the constant defining $\Xi_{T}^{\mathrm{ann}}(d)$ in \eqref{eq:sup.localization.weights} is enlarged if necessary. Applying \eqref{eq:sup.localization.conditional.laplace} proves
\begin{equation}\label{eq:sup.localization.weighted.transforms}
\begin{aligned}
\EE_{\star}[e^{-aX^{\star}}\ind_{\mathcal{A}_{T,\mathrm{sp}}^{\star}}\mid E_{\bb{h}}^{\star}] &\leq C\Xi_{T}^{\mathrm{ann}}(d)\Gamma_{a}(r,m_{-}), \\
\EE_{\star}[e^{-aX^{\star} + 2Y^{\star}(h) \pm 2D^{\star}(z)}\ind_{\mathcal{A}_{T,\mathrm{sp}}^{\star}}\mid E_{\bb{h}}^{\star}] &\leq C\Xi_{T}^{\mathrm{ann}}(d)\Gamma_{a}(r,m_{-})e^{r}.
\end{aligned}
\end{equation}

Partition each interval component of $A_{i,d}$ into cells $B$ of length $\ell_{T}$, except for at most one shorter cell of positive length, and treat each singleton component as a separate cell. Since there are at most two components and their total length is at most $2d$, the number of cells is at most
\begin{equation}\label{eq:sup.localization.cell.count}
2 + \frac{2d}{\ell_{T}} \leq Ce^{r}.
\end{equation}
For a cell $B$ with length $L_{B} > 0$, choose $h_{B} \in B$. The fundamental theorem of calculus gives
\[
\sup_{h\in B}Y^{\star}(h) \leq Y^{\star}(h_{B}) + \ell_{T}^{-1}\int_{B}|D^{\star}(z)|\rd{z}.
\]
Because $L_{B}/\ell_{T} \leq 1$, Jensen's inequality then yields
\begin{equation}\label{eq:sup.localization.cell.jensen}
\begin{aligned}
e^{2\sup_{h\in B}Y^{\star}(h)}
&\leq e^{2Y^{\star}(h_{B})}\exp\left\{2\ell_{T}^{-1}\int_{B}|D^{\star}(z)|\rd{z}\right\} \\
&\leq e^{2Y^{\star}(h_{B})}\exp\left\{\frac{2}{L_{B}}\int_{B}|D^{\star}(z)|\rd{z}\right\}
\leq \frac{1}{L_{B}}\int_{B}e^{2Y^{\star}(h_{B}) + 2|D^{\star}(z)|}\rd{z}.
\end{aligned}
\end{equation}
The finite prime sums have continuous paths, so the same bound includes the cell endpoints. For a singleton cell $B = \{h_{B}\}$, its supremum is simply $Y^{\star}(h_{B})$.

Let $W^{\star} \leqdef e^{-aX^{\star}}\ind_{\mathcal{A}_{T,\mathrm{sp}}^{\star}}$. Multiplying \eqref{eq:sup.localization.cell.jensen} by $W^{\star}$, using $e^{2|z|} \leq e^{2z} + e^{-2z}$, and applying the second bound in \eqref{eq:sup.localization.weighted.transforms} gives
\begin{equation}\label{eq:sup.localization.cell.transform}
\begin{aligned}
\EE_{\star}[W^{\star}e^{2\sup_{h\in B}Y^{\star}(h)}\mid E_{\bb{h}}^{\star}]
&\leq \frac{1}{L_{B}}\int_{B}\sum_{\delta\in\{ - 1,1\}}\EE_{\star}[W^{\star}e^{2Y^{\star}(h_{B}) + 2\delta D^{\star}(z)}\mid E_{\bb{h}}^{\star}]\rd{z} \\
&\leq C\Xi_{T}^{\mathrm{ann}}(d)\Gamma_{a}(r,m_{-})e^{r}.
\end{aligned}
\end{equation}
For a singleton cell, the same expectation bound follows from \eqref{eq:sup.localization.conditional.laplace} with $Z^{\star} = - aX^{\star} + 2Y^{\star}(h_{B})$. Define $S^{\star} \leqdef \sum_{B}e^{2\sup_{h\in B}Y^{\star}(h)}$. Summing \eqref{eq:sup.localization.cell.transform} over the at most $Ce^{r}$ cells counted in \eqref{eq:sup.localization.cell.count} introduces an additional factor $e^{r}$ and gives the first inequality below:
\begin{equation}\label{eq:sup.localization.weighted.cell.sum}
\EE_{\star}[W^{\star}S^{\star}\mid E_{\bb{h}}^{\star}] \leq C\Xi_{T}^{\mathrm{ann}}(d)\Gamma_{a}(r,m_{-})e^{2r}, \qquad \EE_{\star}[W^{\star}\mid E_{\bb{h}}^{\star}] \leq C\Xi_{T}^{\mathrm{ann}}(d)\Gamma_{a}(r,m_{-}).
\end{equation}
The second inequality restates the first bound in \eqref{eq:sup.localization.weighted.transforms}.

The exact increment decomposition \eqref{eq:sup.localization.increment} and the definition of $S^{\star}$ imply
\begin{equation}\label{eq:sup.localization.increment.power}
\exp\left\{a\sup_{h\in A_{i,d}}(F_{h}^{\star} - F_{h_{i}}^{\star})\right\}\ind_{\mathcal{A}_{T,\mathrm{sp}}^{\star}} \leq W^{\star}(S^{\star})^{a/2} = (W^{\star}S^{\star})^{a/2}(W^{\star})^{1 - a/2}.
\end{equation}
Since $0 < a/2 < 1$ by \eqref{eq:sup.localization.exponents}, applying H\"older's inequality with exponents $2/a$ and $1/(1 - a/2)$ to \eqref{eq:sup.localization.increment.power}, and using \eqref{eq:sup.localization.weighted.cell.sum}, gives
\begin{equation}\label{eq:sup.localization.annular.exponential}
\begin{aligned}
&\EE_{\star}\left[\exp\left\{a\sup_{h\in A_{i,d}}(F_{h}^{\star} - F_{h_{i}}^{\star})\right\}\ind_{\mathcal{A}_{T,\mathrm{sp}}^{\star}}\mid E_{\bb{h}}^{\star}\right] \leq \left(\EE_{\star}[W^{\star}S^{\star}\mid E_{\bb{h}}^{\star}]\right)^{a/2}\left(\EE_{\star}[W^{\star}\mid E_{\bb{h}}^{\star}]\right)^{1 - a/2} \\
&\leq C\Xi_{T}^{\mathrm{ann}}(d)\Gamma_{a}(r,m_{-})e^{ar} = C\Xi_{T}^{\mathrm{ann}}(d)\exp\left\{ - \left(a(m_{-} - 1) - \frac{a^{2}}{4}\right)r\right\} = C\Xi_{T}^{\mathrm{ann}}(d)e^{-c_{0}r}.
\end{aligned}
\end{equation}
After H\"older's inequality, the spatial cell count contributes $e^{ar/2}$ and the exponential-moment bound for each cell contributes another $e^{ar/2}$, giving the combined factor $e^{ar}$.

If $d > t^{-3}$, then $r = s_{T} + \log d > s_{T} - 3\log t$ by \eqref{eq:scales}, and $s_{T}/t \to 1$ by \eqref{eq:scale.asymptotics}. Thus, \eqref{eq:sup.localization.weights} gives, for all sufficiently large $T$,
\begin{equation}\label{eq:sup.localization.large.radius.absorption}
\log \Xi_{T}^{\mathrm{ann}}(d) \leq C_{k,\e,J}\log t \leq \frac{c_{0}r}{2}, \qquad \Xi_{T}^{\mathrm{ann}}(d)e^{-c_{0}r} \leq e^{-c_{0}r/2}.
\end{equation}
For $d \leq t^{-3}$, $\Xi_{T}^{\mathrm{ann}}(d) = 1$ by \eqref{eq:sup.localization.weights}. With $c = c_{0}/2$, \eqref{eq:sup.localization.annular.exponential} and \eqref{eq:sup.localization.large.radius.absorption} therefore prove \eqref{eq:sup.localization.uniform.annulus} in both cases.

\medskip
{\bf Step 3. Bound on the conditional exceptional probability.} For the Euler field, the lower probability bound for conditioning bins of fixed width $b > 0$ in \citet[proof of Lemma~4.6, Eqs.~(4.73) and~(4.78)]{OuimetExpanding2026} gives
\begin{equation}\label{eq:sup.localization.bin.lower.bound}
\PP_{\mathrm{A}}(E_{\bb{h}}^{\mathrm{A}}) \geq e^{-C_{k,\e,J,b}t}.
\end{equation}
The required height and covariance conditions hold because $u_{T} + x_{j} = \OO(t)$ uniformly in $j \in \{1,\ldots,k\}$ by \eqref{eq:height.asymptotics} and \eqref{eq:scale.asymptotics}, and the eigenvalues of $\Sigma_{T}$ are comparable to $t$ by Lemma~\ref{lem:annular.covariance}. Combining \eqref{eq:sup.localization.bin.lower.bound} with \eqref{eq:sup.localization.exceptional} gives
\begin{equation}\label{eq:sup.localization.conditional.exceptional}
\PP_{\mathrm{A}}((\mathcal{A}_{T,\mathrm{sp}}^{\mathrm{A}})^{c}\mid E_{\bb{h}}^{\mathrm{A}}) \leq \frac{\PP_{\mathrm{A}}((\mathcal{A}_{T,\mathrm{sp}}^{\mathrm{A}})^{c})}{\PP_{\mathrm{A}}(E_{\bb{h}}^{\mathrm{A}})} \leq \exp\{ - c_{1}t^{1 + 2\e_{0}} + C_{k,\e,J,b}t\} = \oo(1).
\end{equation}
The bound in \eqref{eq:sup.localization.conditional.exceptional} is uniform over the representative positions $h_{1},\ldots,h_{k}$ and the allowed bin endpoints $x_{1},\ldots,x_{k}$, and is independent of the annular radius $d$. It proves \eqref{eq:sup.localization.conditional.event} for the Euler field. For the Haar field, $(\mathcal{A}_{T,\mathrm{sp}}^{\mathrm{H}})^{c} = \varnothing$, so its conditional probability is zero.
\end{proof}

\begin{lemma}[Conditional supremum localization]\label{lem:conditional.supremum.localization}
Under the preceding assumptions, let $a,c > 0$ be the constants from Lemma~\ref{lem:annular.supremum}. For $\star \in \{\mathrm{A}, \mathrm{H}\}$, $i \in \{1,\ldots,k\}$, $R \geq 1$, and fixed $x \in \R$, we have, as $T \to \infty$,
\begin{equation}\label{eq:conditional.supremum.localization}
\PP_{\star}\left( \sup_{\substack{h\in B_{j_{i}}^{\e,s} \\
|h - h_{i}|\geq R\ell_{T}}} (F_{h}^{\star} - u_{T}) \geq x \, \middle|\, E_{\bb{h}}^{\star}\right) \leq C_{k, \e,J,b}e^{-ax}R^{-c} + \oo(1).
\end{equation}
The bound is uniform in $s \in [0,1]$, $R \geq 1$, the representative positions $h_{1},\ldots,h_{k}$ in distinct blocks $B_{j_{1}}^{\e,s},\ldots,B_{j_{k}}^{\e,s}$, and the allowed bin endpoints $x_{1},\ldots,x_{k}$.
\end{lemma}

\begin{proof}[Proof of Lemma~\ref{lem:conditional.supremum.localization}]
On the event $E_{\bb{h}}^{\star}$, the conditioned root satisfies $F_{h_{i}}^{\star} - u_{T} \leq \sup J$. Thus
\[
\left\{\sup_{h\in A_{i,d}}(F_{h}^{\star} - u_{T}) \geq x\right\}\cap E_{\bb{h}}^{\star} \subseteq \left\{\sup_{h\in A_{i,d}}(F_{h}^{\star} - F_{h_{i}}^{\star}) \geq x - \sup J\right\}\cap E_{\bb{h}}^{\star}.
\]
Markov's inequality and the annular exponential bound \eqref{eq:sup.localization.uniform.annulus} from Lemma~\ref{lem:annular.supremum} therefore yield
\begin{equation}\label{eq:sup.localization.annular.probability}
\PP_{\star}\left(\left\{\sup_{h\in A_{i,d}}(F_{h}^{\star} - u_{T}) \geq x\right\} \cap \mathcal{A}_{T,\mathrm{sp}}^{\star}\,\middle|\,E_{\bb{h}}^{\star}\right) \leq Ce^{a(\sup J - x)}e^{-cr}.
\end{equation}

The block $\smash{B_{j_{i}}^{\e,s}} = [s + j_{i} + \e,s + j_{i} + 1 - \e] \cap I_{\theta}$ has diameter at most $1 - 2\e < 1$. Since $d_{n+1} = 2d_{n}$, the dyadic intervals $[d_{n},2d_{n}]$ cover distances without gaps. Starting with the interval containing $R\ell_{T}$ therefore covers all distances in $[R\ell_{T},1)$, so the points under consideration are covered by the annuli $A_{i,d_{n}}$ with
\[
n \geq n_{R}, \qquad n_{R} \leqdef \max\{0,\lfloor\log R/\log 2\rfloor\}, \qquad d_{n} \leq 1.
\]
Since $r = n\log 2$ by \eqref{eq:sup.localization.radii}, summing \eqref{eq:sup.localization.annular.probability} and using the geometric-series identity give
\begin{equation}\label{eq:sup.localization.annular.sum}
\begin{aligned}
&\PP_{\star}\left(\left\{\sup_{\substack{h\in B_{j_{i}}^{\e,s} \\ |h - h_{i}|\geq R\ell_{T}}}(F_{h}^{\star} - u_{T}) \geq x\right\} \cap \mathcal{A}_{T,\mathrm{sp}}^{\star}\,\middle|\,E_{\bb{h}}^{\star}\right) \\
&\quad\leq Ce^{a(\sup J - x)}\sum_{n\geq n_{R}}2^{-cn}
= Ce^{a(\sup J - x)}\frac{2^{-cn_{R}}}{1 - 2^{-c}}
\leq C_{k,\e,J,b}e^{-ax}R^{-c}.
\end{aligned}
\end{equation}
If no such annulus intersects the block $B_{j_{i}}^{\e,s}$, the event $\smash{\{\sup_{h\in B_{j_{i}}^{\e,s},\,|h - h_{i}|\geq R\ell_{T}}(F_{h}^{\star} - u_{T}) \geq x\}}$ in \eqref{eq:conditional.supremum.localization} is empty. Combining \eqref{eq:sup.localization.conditional.event} of Lemma~\ref{lem:annular.supremum} with \eqref{eq:sup.localization.annular.sum} proves \eqref{eq:conditional.supremum.localization}.
\end{proof}

\section{Random grids and the probability of no exceedances}\label{sec:void}

For $k \in \N$, let
\begin{equation}\label{max:eq:separated.domain}
\mathscr{D}_{k,T}(d) \leqdef \{\bb{h} \in I_{\theta}^{\times k}: |h_{i} - h_{j}| \geq d \text{ whenever }i \neq j\}.
\end{equation}
For $k = 1$, the separation condition is empty.

\begin{lemma}[Averaged separated kernels]\label{max:kernel.average}
For every fixed $k \in \N$ and $d > 0$, as $T \to \infty$,
\begin{equation}\label{max:eq:kernel.average}
\frac{1}{|I_{\theta}|^{k}}\int_{\mathscr{D}_{k,T}(d)} \mathcal{K}_{T}(\bb{h})\rd{\bb{h}} \to 1,
\end{equation}
where the kernel $\mathcal{K}_{T}(\bb{h})$ is defined in \eqref{eq:local.profile.kernel}. Furthermore, for some constant $C_{k,d}$,
\begin{equation}\label{max:eq:kernel.pointwise}
\sup_{\bb{h}\in\mathscr{D}_{k,T}(d)}\mathcal{K}_{T}(\bb{h}) \leq t^{C_{k,d}}
\end{equation}
for all sufficiently large $T$.
\end{lemma}

\begin{proof}[Proof of Lemma~\ref{max:kernel.average}]
The proof is divided into four steps. Step~1 decomposes $\mathcal{K}_{T}$ into low-prime and middle-prime factors with a uniformly negligible remainder on $\mathscr{D}_{k,T}(d)$, where $|h_{i} - h_{j}| \geq d$ for $i \neq j$. Step~2 shows that the normalized spatial average of the low-prime factor converges to one. Step~3 shows that the normalized spatial integral of the product of the low-prime and middle-prime factors over $\mathscr{D}_{k,T}(d)$ also converges to one. Step~4 establishes the pointwise bound \eqref{max:eq:kernel.pointwise} for the full kernel.

\medskip
{\bf Step 1. Separate the prime ranges.} Fix a sufficiently large $B > 0$, set $P = t^{B}$, and define the low-prime field and its normalized exponential by
\[
F_{h}^{\mathrm{H}, \leq P} \leqdef \Re\sum_{\substack{p\leq P \\
p^{j}\leq X}} \frac{U_{p}^{j}}{jp^{j(\sig + \ii h)}}, \qquad \Pi_{T,P, \lambda_{T}}(h) \leqdef \frac{e^{\lambda_{T} F_{h}^{\mathrm{H}, \leq P}}}{\EE_{\mathrm{H}}[e^{\lambda_{T} F_{0}^{\mathrm{H}, \leq P}}]}.
\]
Let
\[
\begin{aligned}
\mathcal{K}_{T}^{\mathrm{low}}(\bb{h}, \lambda_{T})
&\leqdef \EE_{\mathrm{H}}\left[\prod_{i=1}^{k}\Pi_{T,P, \lambda_{T}}(h_{i})\right], \qquad \Upsilon_{T}(z) \leqdef \sum_{P<p\leq\exp(t^{3})}p^{-1 - 2\eta_{\sig}}e^{\ii z\log p}, \\
\mathcal{K}_{T}^{\mathrm{mid}}(\bb{h}, \lambda_{T})
&\leqdef \exp\left\{\frac{\lambda_{T}^{2}}{2}\sum_{i<j}\Re\Upsilon_{T}(h_{i} - h_{j})\right\}.
\end{aligned}
\]
For each prime $p \leq X$, let
\[
F_{h,p}^{\mathrm{H}} \leqdef \Re\sum_{\substack{a\geq1 \\ p^{a}\leq X}}\frac{U_{p}^{a}}{ap^{a(\sig + \ii h)}}
\]
be its contribution to $F_{h}^{\mathrm{H}}$. Independence across primes and stationarity give
\[
\mathcal{K}_{T}(\bb{h}) = \prod_{p\leq X}\mathcal{K}_{T,p}(\bb{h}), \qquad
\mathcal{K}_{T,p}(\bb{h}) \leqdef \frac{\EE_{\mathrm{H}}\left[\exp\left\{\lambda_{T}\sum_{i=1}^{k}F_{h_{i},p}^{\mathrm{H}}\right\}\right]}{\prod_{i=1}^{k}\EE_{\mathrm{H}}[\exp\{\lambda_{T}F_{h_{i},p}^{\mathrm{H}}\}]}.
\]
Here, stationarity identifies each denominator factor with $\smash{\EE_{\mathrm{H}}[\exp\{\lambda_{T}F_{0,p}^{\mathrm{H}}\}]}$, as in \eqref{eq:local.profile.kernel}.

Each $F_{h,p}^{\mathrm{H}}$ is centered and satisfies $|F_{h,p}^{\mathrm{H}}| \leq Cp^{-1/2}$ uniformly in $h$. Since $\lambda_{T}$ is bounded, Taylor expansion of the logarithms of the numerator and denominator gives, for $p > P$,
\[
\begin{aligned}
\log\mathcal{K}_{T,p}(\bb{h})
&= \frac{\lambda_{T}^{2}}{2}\left\{\Var_{\mathrm{H}}\left(\sum_{i=1}^{k}F_{h_{i},p}^{\mathrm{H}}\right) - \sum_{i=1}^{k}\Var_{\mathrm{H}}(F_{h_{i},p}^{\mathrm{H}})\right\} + \OO_{k}(p^{-3/2}) \\
&= \lambda_{T}^{2}\sum_{1\leq i<j\leq k}\Cov_{\mathrm{H}}(F_{h_{i},p}^{\mathrm{H}},F_{h_{j},p}^{\mathrm{H}}) + \OO_{k}(p^{-3/2}).
\end{aligned}
\]
Thus the individual variance terms cancel, leaving the covariances between distinct positions. The remainder is $\OO_{k}(p^{-3/2})$ because the terms beyond quadratic order are bounded by a constant times $(p^{-1/2})^{3}$. Summing these errors over $p > P$ gives a total error bounded in absolute value by
\[
C_{k}\sum_{P<p\leq X}p^{-3/2}
\leq C_{k}\sum_{n>P}n^{-3/2}
\leq C_{k}P^{-1/2}.
\]
The high-prime estimate (4.116) and the logarithmic expansion preceding (4.121) in \citet[proof of Lemma~4.8]{OuimetExpanding2026} give
\begin{equation}\label{max:eq:kernel.factorization}
\begin{aligned}
\mathcal{K}_{T}(\bb{h}) &= \mathcal{K}_{T}^{\mathrm{low}}(\bb{h}, \lambda_{T})\mathcal{K}_{T}^{\mathrm{mid}}(\bb{h}, \lambda_{T})e^{e_{T}(\bb{h})}, \\
\sup_{\bb{h}\in\mathscr{D}_{k,T}(d)}|e_{T}(\bb{h})| &\leq C_{k,d}\{P^{-1/2} + t^{-3} + (1 + |I_{\theta}|)t^{C}\exp( - ct^{3/2})\} = \oo(1).
\end{aligned}
\end{equation}
Here $\exp(t^{3}) < X$ eventually by \eqref{eq:fixed.parameters}. The last error in \eqref{max:eq:kernel.factorization} tends to zero because $|I_{\theta}| = 2\exp(\theta t)$ and $t^{3/2}$ dominates $t$.

\medskip
{\bf Step 2. Average the low-prime factor.} Set
\[
\overline{\Pi}_{T,P, \lambda_{T}} \leqdef \frac{1}{|I_{\theta}|}\int_{I_{\theta}}\Pi_{T,P, \lambda_{T}}(h)\rd{h}.
\]
The low-prime averaging estimate in \citet[proof of Lemma~4.8, Eq.~(4.109)]{OuimetExpanding2026} states that, for every fixed $r \in \N$,
\begin{equation}\label{max:eq:kernel.low.moments}
\EE_{\mathrm{H}}\big[|\overline{\Pi}_{T,P, \lambda_{T}} - 1|^{r}\big] \to 0.
\end{equation}
Expanding the $k$th power and applying \eqref{max:eq:kernel.low.moments} at the orders $1, \ldots,k$ gives
\[
\left|\EE_{\mathrm{H}}[\overline{\Pi}_{T,P, \lambda_{T}}^{k}] - 1\right|
\leq \sum_{r=1}^{k}\binom{k}{r}\EE_{\mathrm{H}}\big[|\overline{\Pi}_{T,P, \lambda_{T}} - 1|^{r}\big] = \oo(1).
\]
Tonelli's theorem then identifies the required spatial integral:
\begin{equation}\label{max:eq:kernel.low.average}
\frac{1}{|I_{\theta}|^{k}}\int_{I_{\theta}^{\times k}}\mathcal{K}_{T}^{\mathrm{low}}(\bb{h}, \lambda_{T})\rd{\bb{h}} = \EE_{\mathrm{H}}\left[\prod_{i=1}^{k}\frac{1}{|I_{\theta}|}\int_{I_{\theta}}\Pi_{T,P, \lambda_{T}}(h_{i})\rd{h_{i}}\right] = \EE_{\mathrm{H}}[\overline{\Pi}_{T,P, \lambda_{T}}^{k}] = 1 + \oo(1).
\end{equation}
The primewise cumulant expansion \eqref{eq:local.profile.prime.expansion} gives
\begin{equation}\label{max:eq:kernel.low.pointwise}
\begin{aligned}
&\sup_{\bb{h}\in I_{\theta}^{\times k}}\mathcal{K}_{T}^{\mathrm{low}}(\bb{h}, \lambda_{T})
\leq C_{k}\exp\left(C_{k}\sum_{p\leq P}\frac{1}{p}\right)
\leq C_{k}(\log P)^{C_{k}} \\
&= C_{k}(B\log t)^{C_{k}}
\leq C_{k,B}(\log t)^{C_{k}}.
\end{aligned}
\end{equation}
All these estimates hold uniformly when $\lambda_{T}$ is replaced by any real $\lambda$ in a fixed bounded interval containing $\lambda_{T}$ for all sufficiently large $T$.

\medskip
{\bf Step 3. Insert the middle-prime factor and the separation constraint.} Combining the low-prime bound (4.109) with the middle-prime estimate (4.120) in \citet[proof of Lemma~4.8]{OuimetExpanding2026} gives constants $A,C_{k} > 0$, with $C_{k}$ independent of $B$, such that
\begin{equation}\label{max:eq:kernel.middle.error}
\frac{1}{|I_{\theta}|^{k}}\int_{\mathscr{D}_{k,T}(d)}\mathcal{K}_{T}^{\mathrm{low}}(\bb{h}, \lambda_{T})|\mathcal{K}_{T}^{\mathrm{mid}}(\bb{h}, \lambda_{T}) - 1|\rd{\bb{h}} \leq C(\log t)^{A}t^{C_{k}}\left\{t^{-B/2} + \left(\frac{\log t}{|I_{\theta}|}\right)^{1/2}\right\}.
\end{equation}
Choose $B > 2C_{k} + 2$. The two terms on the right of \eqref{max:eq:kernel.middle.error} satisfy
\[
(\log t)^{A}t^{C_{k} - B/2} = \oo(1), \qquad
(\log t)^{A + 1/2}t^{C_{k}}e^{-\theta t/2} = \oo(1).
\]
Consequently, \eqref{max:eq:kernel.low.average}, \eqref{max:eq:kernel.low.pointwise}, and \eqref{max:eq:kernel.middle.error} give the estimates needed below:
\begin{equation}\label{max:eq:kernel.inputs}
\begin{aligned}
&\frac{1}{|I_{\theta}|^{k}}\int_{I_{\theta}^{\times k}}\mathcal{K}_{T}^{\mathrm{low}}(\bb{h}, \lambda_{T})\rd{\bb{h}} = 1 + \oo(1), \qquad
\sup_{\bb{h}\in I_{\theta}^{\times k}}\mathcal{K}_{T}^{\mathrm{low}}(\bb{h}, \lambda_{T}) \leq C_{k,B}(\log t)^{C_{k}}, \\
&\frac{1}{|I_{\theta}|^{k}}\int_{\mathscr{D}_{k,T}(d)}\mathcal{K}_{T}^{\mathrm{low}}(\bb{h}, \lambda_{T})|\mathcal{K}_{T}^{\mathrm{mid}}(\bb{h}, \lambda_{T}) - 1|\rd{\bb{h}} = \oo(1).
\end{aligned}
\end{equation}
A union bound over the pairs of coordinates gives
\[
\Leb\big(I_{\theta}^{\times k}\setminus\mathscr{D}_{k,T}(d)\big)
\leq \sum_{i<j}\Leb\{\bb{h} \in I_{\theta}^{\times k}: |h_{i} - h_{j}| < d\}
\leq 2\binom{k}{2}d|I_{\theta}|^{k - 1}.
\]
The supremum bound in \eqref{max:eq:kernel.inputs} therefore implies
\begin{equation}\label{max:eq:kernel.removed.low}
\frac{1}{|I_{\theta}|^{k}}\int_{I_{\theta}^{\times k}\setminus\mathscr{D}_{k,T}(d)}\mathcal{K}_{T}^{\mathrm{low}}(\bb{h}, \lambda_{T})\rd{\bb{h}}
\leq \frac{C_{k,d,B}(\log t)^{C_{k}}}{|I_{\theta}|} = \oo(1).
\end{equation}
Subtracting the integral in \eqref{max:eq:kernel.removed.low} from the full low-prime integral in \eqref{max:eq:kernel.inputs}, and then using the last line of \eqref{max:eq:kernel.inputs}, yields
\begin{equation}\label{max:eq:kernel.separated.product}
\begin{aligned}
\frac{1}{|I_{\theta}|^{k}}\int_{\mathscr{D}_{k,T}(d)}\mathcal{K}_{T}^{\mathrm{low}}(\bb{h}, \lambda_{T})\rd{\bb{h}} &= 1 + \oo(1), \\
\frac{1}{|I_{\theta}|^{k}}\int_{\mathscr{D}_{k,T}(d)}\mathcal{K}_{T}^{\mathrm{low}}(\bb{h}, \lambda_{T})\mathcal{K}_{T}^{\mathrm{mid}}(\bb{h}, \lambda_{T})\rd{\bb{h}} &= 1 + \oo(1).
\end{aligned}
\end{equation}
Finally, the factorization and error bound in \eqref{max:eq:kernel.factorization}, together with \eqref{max:eq:kernel.separated.product}, give
\[
\begin{aligned}
&\frac{1}{|I_{\theta}|^{k}}\left|\int_{\mathscr{D}_{k,T}(d)}\mathcal{K}_{T}(\bb{h})\rd{\bb{h}} - \int_{\mathscr{D}_{k,T}(d)}\mathcal{K}_{T}^{\mathrm{low}}(\bb{h}, \lambda_{T})\mathcal{K}_{T}^{\mathrm{mid}}(\bb{h}, \lambda_{T})\rd{\bb{h}}\right| \\
&\qquad\leq \sup_{\bb{h}\in\mathscr{D}_{k,T}(d)}|e^{e_{T}(\bb{h})} - 1|\,
\frac{1}{|I_{\theta}|^{k}}\int_{\mathscr{D}_{k,T}(d)}\mathcal{K}_{T}^{\mathrm{low}}(\bb{h}, \lambda_{T})\mathcal{K}_{T}^{\mathrm{mid}}(\bb{h}, \lambda_{T})\rd{\bb{h}} = \oo(1).
\end{aligned}
\]
Together with the second line of \eqref{max:eq:kernel.separated.product}, this proves \eqref{max:eq:kernel.average}.

\medskip
{\bf Step 4. Bound the full kernel pointwise.} The primewise cumulant expansion from \citet[proof of Lemma~4.6, Eqs.~(4.54)--(4.55)]{OuimetExpanding2026} gives
\begin{equation}\label{max:eq:kernel.cumulant.expansion}
\log\mathcal{K}_{T}(\bb{h}) = \lambda_{T}^{2}\sum_{i<j}\Cov_{\mathrm{H}}(F_{h_{i}}^{\mathrm{H}},F_{h_{j}}^{\mathrm{H}}) + \OO_{k}(1).
\end{equation}
For $\bb{h} \in \mathscr{D}_{k,T}(d)$, applying the covariance bound of \citet[Lemma~4.3, Eqs.~(4.13)--(4.14)]{OuimetExpanding2026} and the boundedness of $\lambda_{T}$ in \eqref{eq:height.asymptotics} to \eqref{max:eq:kernel.cumulant.expansion} gives
\[
\log\mathcal{K}_{T}(\bb{h}) \leq C_{k,d}\log t, \qquad
\mathcal{K}_{T}(\bb{h}) \leq \exp(C_{k,d}\log t) = t^{C_{k,d}}.
\]
This proves \eqref{max:eq:kernel.pointwise}.
\end{proof}

Fix $0 < \e < 1/4$. Let $S$ be uniform on $[0,1]$ and independent of the Euler or Haar field $F^{\star}$. Define the blocks
\begin{equation}\label{max:eq:retained.blocks}
B_{j}(S) \leqdef I_{\theta}\cap[S + j + \e,S + j + 1 - \e], \qquad j \in \Z.
\end{equation}
Let $\mathcal{B}_{T}(S)$ be the finite collection of all nonempty blocks $B_{j}(S)$, including any partial blocks at the endpoints of $I_{\theta}$. Distinct blocks in this collection are separated by at least $2\e$. Write $\bigcup\mathcal{B}_{T}(S) \leqdef \bigcup_{B\in\mathcal{B}_{T}(S)}B$ for their union, and define
\begin{equation}\label{max:eq:distinct.block.domain}
\mathcal{R}_{k,T}(S) \leqdef \bigcup_{\substack{(B_{1}, \ldots,B_{k})\in\mathcal{B}_{T}(S)^{k} \\
B_{i}\text{ distinct}}} B_{1}\times\cdots\times B_{k}.
\end{equation}
The union is over ordered tuples of distinct blocks.

\begin{lemma}[Kernel integral bounds over distinct blocks]\label{max:retained.kernel}
For every fixed $k \in \N$,
\begin{equation}\label{max:eq:retained.kernel}
\begin{aligned}
(1 - 2k\e)_{+} &\leq \liminf_{T\to\infty} \frac{1}{|I_{\theta}|^{k}}\EE_{S}\left[\int_{\mathcal{R}_{k,T}(S)} \mathcal{K}_{T}(\bb{h})\rd{\bb{h}}\right] \\
&\leq \limsup_{T\to\infty} \frac{1}{|I_{\theta}|^{k}}\EE_{S}\left[\int_{\mathcal{R}_{k,T}(S)} \mathcal{K}_{T}(\bb{h})\rd{\bb{h}}\right] \leq 1.
\end{aligned}
\end{equation}
\end{lemma}

\begin{proof}[Proof of Lemma~\ref{max:retained.kernel}]
By \eqref{max:eq:retained.blocks}, distinct blocks in $\mathcal{B}_{T}(S)$ have distance at least $2\e$. Hence, by the definition \eqref{max:eq:separated.domain}, $\mathcal{R}_{k,T}(S) \subseteq \mathscr{D}_{k,T}(2\e)$ for every $S$, and positivity of $\mathcal{K}_{T}(\bb{h})$ gives
\begin{equation}\label{max:eq:block.kernel.upper}
\frac{1}{|I_{\theta}|^{k}}\EE_{S}\left[\int_{\mathcal{R}_{k,T}(S)}\mathcal{K}_{T}(\bb{h})\rd{\bb{h}}\right]
\leq \frac{1}{|I_{\theta}|^{k}}\int_{\mathscr{D}_{k,T}(2\e)}\mathcal{K}_{T}(\bb{h})\rd{\bb{h}} = 1 + \oo(1),
\end{equation}
where the last equality is \eqref{max:eq:kernel.average} of Lemma~\ref{max:kernel.average} with $d = 2\e$.

For the lower bound, fix $h \in I_{\theta}$. The fractional part of $h - S$ is uniform on $[0,1]$, and $h \in \bigcup\mathcal{B}_{T}(S)$ precisely when this fractional part lies in $[\e,1 - \e]$, apart from null endpoint events. Hence
\begin{equation}\label{max:eq:retention.probability}
\PP_{S}\{h \notin \bigcup\mathcal{B}_{T}(S)\} = 2\e.
\end{equation}
Every block in $\mathcal{B}_{T}(S)$ has length at most $1 - 2\e < 1$. Thus, when $\bb{h} \in \mathscr{D}_{k,T}(1)$ and $h_{i} \in \bigcup\mathcal{B}_{T}(S)$ for every $i \in \{1,\ldots,k\}$, the coordinates $h_{1},\ldots,h_{k}$ lie in distinct blocks. The union bound applied to the complementary events $\{h_{i} \notin \bigcup\mathcal{B}_{T}(S)\}$, together with \eqref{max:eq:retention.probability}, yields
\[
\PP_{S}\{\bb{h} \in \mathcal{R}_{k,T}(S)\} \geq 1 - \sum_{i=1}^{k}\PP_{S}\{h_{i} \notin \bigcup\mathcal{B}_{T}(S)\} = 1 - 2k\e.
\]
The left-hand side is also nonnegative, so its lower bound can be replaced by $(1 - 2k\e)_{+}$. Tonelli's theorem now gives
\begin{equation}\label{max:eq:block.kernel.lower}
\begin{aligned}
\frac{1}{|I_{\theta}|^{k}}\EE_{S}\left[\int_{\mathcal{R}_{k,T}(S)}\mathcal{K}_{T}(\bb{h})\rd{\bb{h}}\right]
&= \frac{1}{|I_{\theta}|^{k}}\int_{I_{\theta}^{\times k}}\mathcal{K}_{T}(\bb{h})\PP_{S}\{\bb{h} \in \mathcal{R}_{k,T}(S)\}\rd{\bb{h}} \\
&\geq \frac{(1 - 2k\e)_{+}}{|I_{\theta}|^{k}}\int_{\mathscr{D}_{k,T}(1)}\mathcal{K}_{T}(\bb{h})\rd{\bb{h}} = (1 - 2k\e)_{+}\{1 + \oo(1)\},
\end{aligned}
\end{equation}
where the last equality follows from \eqref{max:eq:kernel.average} of Lemma~\ref{max:kernel.average} with $d = 1$. Taking the upper limit in \eqref{max:eq:block.kernel.upper} and the lower limit in \eqref{max:eq:block.kernel.lower} as $T \to \infty$ proves \eqref{max:eq:retained.kernel}.
\end{proof}

Fix $q > 0$. Conditionally on $S$, attach independent variables $U_{B}$, uniform on $[0,q\ell_{T})$, to the blocks $B \in \mathcal{B}_{T}(S)$, independently of the Euler or Haar field $F^{\star}$. In block $B$, use the finite grid
\begin{equation}\label{max:eq:block.grid}
\mathscr{G}_{B} \leqdef B\cap(U_{B} + q\ell_{T}\Z).
\end{equation}
Set $\max\varnothing = - \infty$. For fixed real $x < H$, define
\begin{equation}\label{max:eq:band.count}
\mathfrak{W}_{T}(x,H;q, \e) \leqdef \sum_{B\in\mathcal{B}_{T}(S)} \ind_{\{u_{T} + x<\max_{h\in\mathscr{G}_{B}}F_{h}^{\star}\leq u_{T} + H\}}.
\end{equation}
For events and random variables involving these grids, $\PP_{\star}$ and $\EE_{\star}$ denote probability and expectation under the joint law of the field $F^{\star}$, the block shift $S$, and the grid offsets $U_{B}$.

Recall the limiting Gaussian profile $D_{\theta}$ from \eqref{eq:main.profile} and the constants $p(q)$ and $\kappa_{\theta}(q)$ from \eqref{eq:lattice.constant}. For $R > 0$, define
\[
p(q,R) \leqdef \PP\{D_{\theta}(nq) \leq 0 \text{ for every }n \in \Z \text{ with }0 < |nq| \leq R\}.
\]
By continuity from above,
\begin{equation}\label{max:eq:grid.constant}
p(q,R) \downarrow p(q) = q\kappa_{\theta}(q), \qquad R \to \infty.
\end{equation}

\begin{proposition}[Factorial moments of the band count]\label{max:factorial.moments}
For every fixed $q > 0$, $x < H$, $0 < \e < 1/4$, and $k \in \N$, let
\begin{equation}\label{max:eq:band.parameter}
\Lambda_{q,x,H} \leqdef \kappa_{\theta}(q)(e^{-\beta_{c}(\theta) x} - e^{-\beta_{c}(\theta) H}).
\end{equation}
Then
\begin{equation}\label{max:eq:factorial.bounds}
(1 - 2k\e)_{+}\Lambda_{q,x,H}^{k} \leq \liminf_{T\to\infty}\EE_{\star}[(\mathfrak{W}_{T})_{k}] \leq \limsup_{T\to\infty}\EE_{\star}[(\mathfrak{W}_{T})_{k}] \leq \Lambda_{q,x,H}^{k},
\end{equation}
where $(n)_{k} = n(n - 1)\cdots(n - k + 1)$.
\end{proposition}

\begin{proof}[Proof of Proposition~\ref{max:factorial.moments}]
The proof is divided into four steps. Step~1 expresses the $k$th factorial moment $\EE_{\star}[(\mathfrak{W}_{T})_{k}]$ of the band count $\mathfrak{W}_{T} = \mathfrak{W}_{T}(x,H;q,\e)$ as an integral over the possible positions $h_{1},\ldots,h_{k}$ of grid maximizers in distinct blocks. This is achieved by averaging the grid offsets $U_{B}$, which, conditionally on $S$, are independent and uniform on $[0,q\ell_{T})$. Step~2 shows that, for fixed $R > 0$, the contribution to this integral from tuples with at least one root $h_{i}$ within distance $R\ell_{T}$ of an endpoint of its block tends to zero as $T \to \infty$. Step~3 conditions on all root heights satisfying $u_{T} + x < F_{h_{i}}^{\star} \leq u_{T} + H$ and uses conditional localization to replace the requirement that each root maximize the field on the grid in its entire block by the requirement that it maximize the field over the grid points within distance $R\ell_{T}$. The resulting error vanishes as $T \to \infty$ and then $R \to \infty$. Step~4 evaluates the contribution of these local lattice maximizers to $\EE_{\star}[(\mathfrak{W}_{T})_{k}]$ using the lattice-bin asymptotic in Corollary~\ref{cor:local.profile.lattice}, then applies the kernel integral bounds in Lemma~\ref{max:retained.kernel} to obtain the stated factorial-moment bounds.

\medskip
{\bf Step 1. Express the factorial moment through block maximizers.} To rewrite the band count $\mathfrak{W}_{T}$ using the sums over maximizing grid positions in \eqref{max:eq:factorial.block.sum}, each block whose maximum lies in $(u_{T} + x,u_{T} + H]$ must contribute exactly one summand. This requires uniqueness of the grid maximizer: if two grid positions attained the same block maximum, summing over all maximizing positions would count that block twice.

Condition on the block shift $S$ and the grid offsets $(U_{B})_{B\in\mathcal{B}_{T}(S)}$, so that all grid positions are fixed. For two distinct positions with separation $z \neq 0$, their phases cannot agree at both primes $2$ and $3$. Indeed, this would give
\[
z\log2 \in 2\pi\Z, \qquad z\log3 \in 2\pi\Z,
\]
and hence make $\log2/\log3$ rational.

Under the Haar law, choose a prime at which the two phases differ and condition on all other prime phases. The difference of the two field values is then a nonconstant real-analytic function of the remaining prime's uniform angle. Its zero set has measure zero, so the two values are almost surely distinct.

Under the arithmetic law, the difference of the two field values is a real-analytic function of $\tau$. It is nonconstant because the differing phases give a nonzero coefficient at the frequency corresponding to that prime. Its zero set again has measure zero.

A finite union over pairs of grid positions, followed by averaging over $S$ and the grid offsets, proves that every nonempty block grid has an almost surely unique maximizer. This justifies the first identity in \eqref{max:eq:factorial.block.sum}, which counts each block through its single maximizing grid position.

For a block $B$ containing $h$, define the event that $h$ maximizes $F^{\star}$ on the grid $B\cap(h + q\ell_{T}\Z)$ by
\begin{equation}\label{max:eq:block.maximizer}
J_{T}(h,B) \leqdef \left\{F_{h}^{\star} = \max_{v\in B\cap(h + q\ell_{T}\Z)}F_{v}^{\star}\right\}.
\end{equation}
Define the event that the field values at all $k$ positions $h_{1},\ldots,h_{k}$ lie in $(u_{T} + x,u_{T} + H]$, and recall the mass of $(x,H]$ under $\Theta(\rd{y}) = \beta_{c}(\theta)e^{-\beta_{c}(\theta)y}\rd{y}$:
\[
E_{T}(\bb{h}) \leqdef \bigcap_{i=1}^{k}\{u_{T} + x < F_{h_{i}}^{\star} \leq u_{T} + H\}, \qquad
\Theta((x,H]) = e^{-\beta_{c}(\theta) x} - e^{-\beta_{c}(\theta) H}.
\]
For a nonnegative measurable $f$, averaging one grid offset gives the exact identity
\begin{equation}\label{max:eq:grid.offset.average}
\begin{aligned}
\EE_{U_{B}}\left[\sum_{h\in\mathscr{G}_{B}}f(h,\mathscr{G}_{B})\right]
&= \frac{1}{q\ell_{T}}\int_{0}^{q\ell_{T}}\sum_{n\in\Z}\ind_{B}(u + nq\ell_{T})f(u + nq\ell_{T},B\cap(u + q\ell_{T}\Z))\rd{u} \\
&= \frac{1}{q\ell_{T}}\int_{B}f(h,B\cap(h + q\ell_{T}\Z))\rd{h}.
\end{aligned}
\end{equation}
The substitution $h = u + nq\ell_{T}$ partitions $\R$ into the intervals $[nq\ell_{T},(n + 1)q\ell_{T})$ and leaves the shifted lattice unchanged.

Let $I_{B} \leqdef \ind_{\{u_{T} + x < \max_{h\in\mathscr{G}_{B}}F_{h}^{\star} \leq u_{T} + H\}}$. Uniqueness of the maximizer gives
\begin{equation}\label{max:eq:factorial.block.sum}
I_{B} = \sum_{h\in\mathscr{G}_{B}}\ind_{\{u_{T} + x < F_{h}^{\star} \leq u_{T} + H\}}\ind_{J_{T}(h,B)}, \qquad
(\mathfrak{W}_{T})_{k} = \sum_{\substack{(B_{1}, \ldots,B_{k})\in\mathcal{B}_{T}(S)^{k} \\
B_{i}\text{ distinct}}}\prod_{i=1}^{k}I_{B_{i}}.
\end{equation}

Now fix $S$ and an ordered tuple of distinct blocks $(B_{1},\ldots,B_{k}) \in \mathcal{B}_{T}(S)^{k}$. Condition on the field $F^{\star}$, leaving the offsets $U_{B_{1}},\ldots,U_{B_{k}}$ independent and uniform on $[0,q\ell_{T})$. Substituting the first identity in \eqref{max:eq:factorial.block.sum} into $\smash{\prod_{i=1}^{k}I_{B_{i}}}$ gives one sum over grid positions in each block. Applying \eqref{max:eq:grid.offset.average} successively to these $k$ sums replaces each sum and its offset average by an integral over $h_{i} \in B_{i}$, with prefactor $(q\ell_{T})^{-1}$. Consequently,
\[
\EE_{U_{B_{1}},\ldots,U_{B_{k}}}\left[\prod_{i=1}^{k}I_{B_{i}}\right]
= (q\ell_{T})^{-k}\int_{B_{1}\times\cdots\times B_{k}}
\ind_{E_{T}(\bb{h})}\prod_{i=1}^{k}\ind_{J_{T}(h_{i},B_{i})}\rd{\bb{h}}.
\]
Here the expectation averages only the indicated offsets, with $S$ and the field $(F_{h}^{\star})_{h\in I_{\theta}}$ fixed. After averaging the offset in block $B_{i}$, the grid containing the integration point $h_{i}$ is $B_{i}\cap(h_{i} + q\ell_{T}\Z)$, which explains the event $J_{T}(h_{i},B_{i})$ in the integrand.

Taking expectation with respect to the field turns the product of indicators into the probability of their intersection. Summing over the ordered tuples of distinct blocks in \eqref{max:eq:factorial.block.sum} combines the disjoint integration domains $B_{1}\times\cdots\times B_{k}$ into $\mathcal{R}_{k,T}(S)$, as defined in \eqref{max:eq:distinct.block.domain}. Finally, averaging over $S$ yields
\begin{equation}\label{max:eq:factorial.root.identity}
\EE_{\star}[(\mathfrak{W}_{T})_{k}] = (q\ell_{T})^{-k}\EE_{S}\left[\int_{\mathcal{R}_{k,T}(S)}\PP_{\star}\left(E_{T}(\bb{h})\cap\bigcap_{i=1}^{k}J_{T}(h_{i},B_{S}(h_{i}))\right)\rd{\bb{h}}\right],
\end{equation}
where $B_{S}(h)$ is the unique block in $\mathcal{B}_{T}(S)$ containing $h$. The probability inside the integral is taken over the field with $S$ and $\bb{h}$ fixed.

\medskip
{\bf Step 2. Remove roots near block endpoints.} Fix $R \geq \max\{1,q\}$. In each block $B \in \mathcal{B}_{T}(S)$, remove the points at distance at most $R\ell_{T}$ from either endpoint. If $V_{T,R}(S)$ denotes their union, then each of the $\OO(1 + |I_{\theta}|)$ blocks contributes at most $2R\ell_{T}$. Hence
\begin{equation}\label{max:eq:factorial.boundary.volume}
\Leb(V_{T,R}(S)) \leq 2R\ell_{T}|\mathcal{B}_{T}(S)| \leq C R\ell_{T}(1 + |I_{\theta}|).
\end{equation}
The tuples removed from $\mathcal{R}_{k,T}(S)$ in \eqref{max:eq:distinct.block.domain} are those with at least one coordinate in $V_{T,R}(S)$. Their $k$-dimensional volume satisfies
\begin{equation}\label{max:eq:factorial.boundary.tuple.volume}
\Leb\left(\left\{\bb{h} \in \mathcal{R}_{k,T}(S): h_{i} \in V_{T,R}(S)\text{ for some }i \in \{1,\ldots,k\}\right\}\right)
\leq k\Leb(V_{T,R}(S))|I_{\theta}|^{k - 1}.
\end{equation}

The separated-height estimate in \citet[proof of Lemma~4.8, Eq.~(4.133)]{OuimetExpanding2026} gives, for some fixed $\delta > 0$,
\begin{equation}\label{max:eq:separated.height.bound}
\PP_{\star}(E_{T}(\bb{h})) \leq C_{k,x,H}N_{T}^{-k}\mathcal{K}_{T}(\bb{h}) + \exp( - ct^{1 + \delta} + C_{k}t), \qquad \bb{h} \in \mathscr{D}_{k,T}(2\e).
\end{equation}
Restricting the integral in \eqref{max:eq:factorial.root.identity} to the discarded tuples and dropping the block-maximizer conditions bounds the integrand by $\PP_{\star}(E_{T}(\bb{h}))$. To estimate this contribution, use \eqref{max:eq:factorial.boundary.tuple.volume}, \eqref{max:eq:factorial.boundary.volume}, \eqref{max:eq:separated.height.bound}, \eqref{max:eq:kernel.pointwise} of Lemma~\ref{max:kernel.average}, and $\ell_{T}N_{T} = |I_{\theta}|$ from \eqref{eq:scales}. The contribution of the discarded tuples to \eqref{max:eq:factorial.root.identity} is
\begin{equation}\label{max:eq:factorial.boundary.error}
\begin{aligned}
&(q\ell_{T})^{-k}\EE_{S}\left[\int_{\mathcal{R}_{k,T}(S)} \ind_{\{\exists i:\,h_{i}\in V_{T,R}(S)\}} \PP_{\star}\left(E_{T}(\bb{h})\cap \bigcap_{i = 1}^{k}J_{T}(h_{i},B_{S}(h_{i}))\right)\rd{\bb{h}}\right] \\
&\quad\leq (q\ell_{T})^{-k}\EE_{S}\left[\int_{\mathcal{R}_{k,T}(S)} \ind_{\{\exists i:\,h_{i}\in V_{T,R}(S)\}} \PP_{\star}(E_{T}(\bb{h}))\rd{\bb{h}}\right] \\
&\quad\leq Ck(q\ell_{T})^{-k}N_{T}^{-k}t^{C_{k,\e}} R\ell_{T}(1 + |I_{\theta}|)|I_{\theta}|^{k - 1} + (q\ell_{T})^{-k}|I_{\theta}|^{k} \exp( - ct^{1 + \delta} + C_{k}t) \\
&\quad= Ckq^{-k}R\ell_{T}t^{C_{k,\e}} \frac{1 + |I_{\theta}|}{|I_{\theta}|} + q^{-k}N_{T}^{k}\exp( - ct^{1 + \delta} + C_{k}t) = \oo(1).
\end{aligned}
\end{equation}

\medskip
{\bf Step 3. Replace block maximizers by local maximizers.} Write
\begin{equation}\label{max:eq:interior.root.domain}
\mathcal{R}_{k,T,R}(S) \leqdef \mathcal{R}_{k,T}(S)\cap(I_{\theta}\setminus V_{T,R}(S))^{k}.
\end{equation}
By \eqref{max:eq:factorial.root.identity} and \eqref{max:eq:factorial.boundary.error}, for fixed $R$,
\begin{equation}\label{max:eq:factorial.interior.identity}
\EE_{\star}[(\mathfrak{W}_{T})_{k}] = (q\ell_{T})^{-k}\EE_{S}\left[\int_{\mathcal{R}_{k,T,R}(S)}\PP_{\star}\left(E_{T}(\bb{h})\cap\bigcap_{i = 1}^{k}J_{T}(h_{i},B_{S}(h_{i}))\right)\rd{\bb{h}}\right] + \oo(1).
\end{equation}
For a root $h_{i}$ of a tuple $\bb{h}$ in the domain \eqref{max:eq:interior.root.domain}, define
\[
J_{T,R}(h_{i}) \leqdef \{F_{h_{i} + nq\ell_{T}}^{\star} - F_{h_{i}}^{\star} \leq 0\text{ for every }n \in \Z\text{ with }0 < |nq| \leq R\}.
\]
Every tested point belongs to $B_{S}(h_{i})$, so $J_{T}(h_{i},B_{S}(h_{i}))$ from \eqref{max:eq:block.maximizer} implies $J_{T,R}(h_{i})$. To identify how the converse can fail, suppose that $J_{T,R}(h_{i})$ holds for every $i \in \{1,\ldots,k\}$ but $J_{T}(h_{j},B_{S}(h_{j}))$ fails for some $j$. Then a grid point $h_{j} + nq\ell_{T} \in B_{S}(h_{j})$ with $|nq| > R$ satisfies $\smash{F_{h_{j} + nq\ell_{T}}^{\star}} > \smash{F_{h_{j}}^{\star}}$. Define the event that some point $z$ in a root's block $B_{S}(h_{i})$, at distance greater than $R\ell_{T}$ from the root $h_{i}$, has a larger field value $F_{z}^{\star} > F_{h_{i}}^{\star}$:
\begin{equation}\label{max:eq:remote.block.exceedance}
A_{T,R}(\bb{h},S) \leqdef \bigcup_{i = 1}^{k}\left\{\sup_{\substack{z \in B_{S}(h_{i}) \\ |z - h_{i}| > R\ell_{T}}}(F_{z}^{\star} - F_{h_{i}}^{\star}) > 0\right\}.
\end{equation}
The preceding implication gives the event inclusion
\[
\left(\bigcap_{i = 1}^{k}J_{T,R}(h_{i})\right)\setminus\left(\bigcap_{i = 1}^{k}J_{T}(h_{i},B_{S}(h_{i}))\right) \subseteq A_{T,R}(\bb{h},S).
\]
On $E_{T}(\bb{h})$, every root $h_{i}$ satisfies $F_{h_{i}}^{\star} > u_{T} + x$. Hence the event $A_{T,R}(\bb{h},S)$ in \eqref{max:eq:remote.block.exceedance}, together with $E_{T}(\bb{h})$, implies
\[
\sup_{\substack{z \in B_{S}(h_{i}) \\ |z - h_{i}| \geq R\ell_{T}}}(F_{z}^{\star} - u_{T}) \geq x
\quad\text{for at least one }i \in \{1,\ldots,k\}.
\]
Apply Lemma~\ref{lem:conditional.supremum.localization} with $s = S$, $J = [x,H]$, $b = H - x$, and $x_{i} = x$ for every $i$. The strict lower inequalities in $E_{T}(\bb{h})$ do not affect the conditioning, since the root heights have no atoms by the real-analytic argument in Step~1 applied to $F_{h}^{\star} - y$ for fixed $y \in \R$.

The roots belong to distinct blocks in $\mathcal{B}_{T}(S)$, so all the hypotheses of Lemma~\ref{lem:conditional.supremum.localization} hold. A union bound over $i \in \{1,\ldots,k\}$ yields, uniformly for $\bb{h} \in \mathcal{R}_{k,T,R}(S)$ and $S \in [0,1]$,
\begin{equation}\label{max:eq:leader.conditional.localization}
\PP_{\star}(A_{T,R}(\bb{h},S)\mid E_{T}(\bb{h}))
\leq \sum_{i = 1}^{k}\PP_{\star}\left(\sup_{\substack{z \in B_{S}(h_{i}) \\ |z - h_{i}| \geq R\ell_{T}}}(F_{z}^{\star} - u_{T}) \geq x \,\middle|\, E_{T}(\bb{h})\right)
\leq C R^{-c} + \oo(1).
\end{equation}
Here $k$ and $x$ are fixed, so the factor $ke^{-ax}$ from \eqref{eq:conditional.supremum.localization} of Lemma~\ref{lem:conditional.supremum.localization} is included in $C$.

To integrate this conditional error, first bound the total contribution of the root-height event $E_{T}(\bb{h})$. The separated-height estimate \eqref{max:eq:separated.height.bound} and $\ell_{T}N_{T} = |I_{\theta}|$ from \eqref{eq:scales} give
\begin{equation}\label{max:eq:factorial.root.mass}
\begin{aligned}
&\ell_{T}^{-k}\int_{\mathscr{D}_{k,T}(2\e)}\PP_{\star}(E_{T}(\bb{h}))\rd{\bb{h}}
\leq \frac{C_{k,x,H}}{|I_{\theta}|^{k}}\int_{\mathscr{D}_{k,T}(2\e)}\mathcal{K}_{T}(\bb{h})\rd{\bb{h}} + N_{T}^{k}\exp( - ct^{1 + \delta} + C_{k}t) \\
&\leq C'_{k,x,H}.
\end{aligned}
\end{equation}
The first term is bounded by the kernel average \eqref{max:eq:kernel.average} of Lemma~\ref{max:kernel.average} with $d = 2\e$. The second tends to zero because $\log N_{T} = \OO(t)$ by \eqref{eq:scale.asymptotics}, whereas $t^{1 + \delta}/t \to \infty$.

For each tuple in $\mathcal{R}_{k,T,R}(S)$, the block-maximizer events $J_{T}(h_{i},B_{S}(h_{i}))$ in \eqref{max:eq:block.maximizer} imply the local events $J_{T,R}(h_{i})$. The difference of the corresponding joint probabilities is therefore nonnegative. The event inclusion above and \eqref{max:eq:leader.conditional.localization} bound it by
\[
\begin{aligned}
\PP_{\star}(E_{T}(\bb{h})\cap A_{T,R}(\bb{h},S))
&= \PP_{\star}(E_{T}(\bb{h}))\PP_{\star}(A_{T,R}(\bb{h},S)\mid E_{T}(\bb{h})) \\
&\leq \{CR^{-c} + \oo(1)\}\PP_{\star}(E_{T}(\bb{h})).
\end{aligned}
\]
Integrate this bound over $\mathcal{R}_{k,T,R}(S) \subseteq \mathscr{D}_{k,T}(2\e)$, average over $S$, and apply \eqref{max:eq:factorial.root.mass}. The inclusion holds because the coordinates of each tuple belong to distinct blocks in $\mathcal{B}_{T}(S)$, whose mutual distances are at least $2\e$ by \eqref{max:eq:retained.blocks}. With the prefactor $(q\ell_{T})^{-k}$ from \eqref{max:eq:factorial.interior.identity}, the result is
\begin{equation}\label{max:eq:leader.localization}
\begin{aligned}
&0
\leq (q\ell_{T})^{-k}\EE_{S}\left[\int_{\mathcal{R}_{k,T,R}(S)}\left\{\PP_{\star}\left(E_{T}(\bb{h})\cap\bigcap_{i = 1}^{k}J_{T,R}(h_{i})\right) \right. \right. \\
&\hspace{60mm} \left. \left. - \PP_{\star}\left(E_{T}(\bb{h})\cap\bigcap_{i = 1}^{k}J_{T}(h_{i},B_{S}(h_{i}))\right)\right\}\rd{\bb{h}}\right] \\
&\leq q^{-k}\{CR^{-c} + \oo(1)\}\ell_{T}^{-k}\int_{\mathscr{D}_{k,T}(2\e)}\PP_{\star}(E_{T}(\bb{h}))\rd{\bb{h}} \leq Cq^{-k}R^{-c} + \oo(1).
\end{aligned}
\end{equation}
For each fixed $R$, the $\oo(1)$ terms tend to zero as $T \to \infty$. The term $Cq^{-k}R^{-c}$ then tends to zero as $R \to \infty$, with $q$ fixed.

\medskip
{\bf Step 4. Evaluate the contribution of the local maximizers.} For every $S \in [0,1]$, each tuple $\bb{h} \in \mathcal{R}_{k,T,R}(S)$ from \eqref{max:eq:interior.root.domain} satisfies $|h_{i} - h_{j}| \geq 2\e$ for $i \neq j$ and $[h_{i} - R\ell_{T},h_{i} + R\ell_{T}] \subseteq I_{\theta}$ for every $i \in \{1,\ldots,k\}$. Indeed, its coordinates lie in distinct blocks, which are separated by at least $2\e$ by \eqref{max:eq:retained.blocks}, and each coordinate is more than $R\ell_{T}$ from the endpoints of its block. Thus each $R\ell_{T}$ neighborhood remains inside its block and hence inside $I_{\theta}$.

Let $r = \lfloor R/q\rfloor \geq 1$. The lattice sites tested in $J_{T,R}(h_{i})$ are exactly those indexed by $n \in \{ - r,\ldots, - 1,1,\ldots,r\}$, because
\[
\{n \in \Z:0 < |nq| \leq R\} = \{n \in \Z:1 \leq |n| \leq r\}.
\]
In particular, $p(q,qr) = p(q,R)$, where $p(q,R)$ is the finite-lattice probability used in \eqref{max:eq:grid.constant}. Corollary~\ref{cor:local.profile.lattice}, specifically \eqref{eq:local.profile.lattice.bin}, therefore applies with separation $d_{0} = 2\e$, these $2r$ offsets, and the height bin $[x,H]$ at every root, uniformly over $S \in [0,1]$ and $\bb{h} \in \mathcal{R}_{k,T,R}(S)$. The absence of atoms again permits replacing its closed lower endpoint by the open endpoint in $E_{T}(\bb{h})$.

With $R$ fixed, the uniformity in \eqref{eq:local.profile.lattice.bin} of Corollary~\ref{cor:local.profile.lattice} gives a deterministic nonnegative sequence $\e_{T}^{\mathrm{loc}} \to 0$, independent of $S$ and $\bb{h}$, such that, for every $S \in [0,1]$ and $\bb{h} \in \mathcal{R}_{k,T,R}(S)$,
\begin{equation}\label{max:eq:local.leader.approximation}
\left|\PP_{\star}\left(E_{T}(\bb{h})\cap\bigcap_{i = 1}^{k}J_{T,R}(h_{i})\right) - N_{T}^{-k}\mathcal{K}_{T}(\bb{h})\Theta((x,H])^{k}p(q,R)^{k}\right| \leq \e_{T}^{\mathrm{loc}}N_{T}^{-k}\mathcal{K}_{T}(\bb{h}).
\end{equation}
Each height bin contributes $\Theta((x,H]) = e^{-\beta_{c}(\theta) x} - e^{-\beta_{c}(\theta) H}$, and each limiting local profile contributes $p(q,R)$. Independence of the $k$ limiting profiles accounts for the product $p(q,R)^{k}$.

The local-maximizer probability in \eqref{max:eq:local.leader.approximation} is the integrand obtained from \eqref{max:eq:factorial.interior.identity} through the replacement justified by \eqref{max:eq:leader.localization}. To use this approximation in the factorial-moment integral, its error must remain negligible after integration and averaging over $S$. Integrating \eqref{max:eq:local.leader.approximation} over $\mathcal{R}_{k,T,R}(S)$ with respect to $\bb{h}$, averaging over $S$, multiplying by $(q\ell_{T})^{-k}$, and using $(q\ell_{T})^{-k}\smash{N_{T}^{-k}} = q^{-k}|I_{\theta}|^{-k}$ from \eqref{eq:scales} gives
\begin{equation}\label{max:eq:marked.local.input}
\begin{aligned}
&(q\ell_{T})^{-k}\EE_{S}\bigg[\int_{\mathcal{R}_{k,T,R}(S)}
\bigg|\PP_{\star}\bigg(E_{T}(\bb{h})\cap\bigcap_{i = 1}^{k}J_{T,R}(h_{i})\bigg) \\
&\qquad\hspace{30mm} - N_{T}^{-k}\mathcal{K}_{T}(\bb{h})\Theta((x,H])^{k}p(q,R)^{k}\bigg|\rd{\bb{h}}\bigg] \\
&\qquad\leq \frac{q^{-k}\e_{T}^{\mathrm{loc}}}{|I_{\theta}|^{k}} \EE_{S}\left[\int_{\mathcal{R}_{k,T,R}(S)}\mathcal{K}_{T}(\bb{h})\rd{\bb{h}}\right]
\leq \frac{q^{-k}\e_{T}^{\mathrm{loc}}}{|I_{\theta}|^{k}} \int_{\mathscr{D}_{k,T}(2\e)}\mathcal{K}_{T}(\bb{h})\rd{\bb{h}} = \oo(1).
\end{aligned}
\end{equation}
The second inequality uses positivity of $\mathcal{K}_{T}$ and the inclusion $\mathcal{R}_{k,T,R}(S) \subseteq \mathscr{D}_{k,T}(2\e)$ for every $S$. The last equality follows from the convergence of the normalized kernel integral to $1$ in \eqref{max:eq:kernel.average} of Lemma~\ref{max:kernel.average}, applied with $d = 2\e$, and $\e_{T}^{\mathrm{loc}} \to 0$, with $q$ fixed. Consequently, the triangle inequality and \eqref{max:eq:marked.local.input} give
\begin{equation}\label{max:eq:interior.kernel.approximation}
\begin{aligned}
&(q\ell_{T})^{-k}\EE_{S}\left[\int_{\mathcal{R}_{k,T,R}(S)}
\PP_{\star}\left(E_{T}(\bb{h})\cap\bigcap_{i = 1}^{k}J_{T,R}(h_{i})\right)\rd{\bb{h}}\right] \\
&\qquad= \left\{\frac{p(q,R)}{q}\Theta((x,H])\right\}^{k}
\frac{1}{|I_{\theta}|^{k}}
\EE_{S}\left[\int_{\mathcal{R}_{k,T,R}(S)}\mathcal{K}_{T}(\bb{h})\rd{\bb{h}}\right] + \oo(1).
\end{aligned}
\end{equation}

It remains to replace the interior domain $\mathcal{R}_{k,T,R}(S)$ by $\mathcal{R}_{k,T}(S)$ in the normalized kernel integral on the right-hand side of \eqref{max:eq:interior.kernel.approximation}. The difference $\mathcal{R}_{k,T}(S)\setminus\mathcal{R}_{k,T,R}(S)$ consists of tuples with some $h_{i} \in V_{T,R}(S)$. The volume bounds \eqref{max:eq:factorial.boundary.tuple.volume} and \eqref{max:eq:factorial.boundary.volume}, together with the pointwise kernel bound \eqref{max:eq:kernel.pointwise} of Lemma~\ref{max:kernel.average}, yield
\begin{equation}\label{max:eq:factorial.boundary.kernel}
\begin{aligned}
\frac{1}{|I_{\theta}|^{k}}\EE_{S}\left[\int_{\mathcal{R}_{k,T}(S)\setminus\mathcal{R}_{k,T,R}(S)}\mathcal{K}_{T}(\bb{h})\rd{\bb{h}}\right]
&\leq \frac{k t^{C_{k,\e}}}{|I_{\theta}|}\EE_{S}[\Leb(V_{T,R}(S))] \\
&\leq CkR\ell_{T}t^{C_{k,\e}}\frac{1 + |I_{\theta}|}{|I_{\theta}|} = \oo(1).
\end{aligned}
\end{equation}
Here $R$ is fixed, $|I_{\theta}| = 2e^{\theta t} \to \infty$, and $\ell_{T} \leq Ct^{\Phi}e^{-t}$ by \eqref{eq:fixed.parameters} and \eqref{eq:scales}. Thus the exponential decay in $\ell_{T}$ dominates the fixed power of $t$ in \eqref{max:eq:factorial.boundary.kernel}.

Starting from \eqref{max:eq:factorial.interior.identity}, replace the block-maximizer probabilities by the local-maximizer probabilities using \eqref{max:eq:leader.localization}. Apply \eqref{max:eq:interior.kernel.approximation} to express the resulting integral in terms of the kernel integral over $\mathcal{R}_{k,T,R}(S)$, and then replace this domain by $\mathcal{R}_{k,T}(S)$ using \eqref{max:eq:factorial.boundary.kernel}. For fixed $q$, $x$, $H$, $\e$, $k$, and $R$, these operations give
\[
\bigg|\EE_{\star}[(\mathfrak{W}_{T})_{k}] - \left\{\frac{p(q,R)}{q}\Theta((x,H])\right\}^{k}\frac{1}{|I_{\theta}|^{k}}\EE_{S}\left[\int_{\mathcal{R}_{k,T}(S)}\mathcal{K}_{T}(\bb{h})\rd{\bb{h}}\right]\bigg|
\leq Cq^{-k}R^{-c} + \oo(1).
\]
Taking the upper limit as $T \to \infty$ and then letting $R \to \infty$ proves
\begin{equation}\label{max:eq:factorial.kernel.comparison}
\begin{aligned}
&\lim_{R\to\infty}\limsup_{T\to\infty}\bigg|\EE_{\star}[(\mathfrak{W}_{T})_{k}] \\
&\qquad - \left\{\frac{p(q,R)}{q}(e^{-\beta_{c}(\theta) x} - e^{-\beta_{c}(\theta) H})\right\}^{k} \frac{1}{|I_{\theta}|^{k}}\EE_{S}\left[\int_{\mathcal{R}_{k,T}(S)}\mathcal{K}_{T}(\bb{h})\rd{\bb{h}}\right]\bigg| = 0.
\end{aligned}
\end{equation}

By \eqref{max:eq:grid.constant}, $p(q,R)/q \to \kappa_{\theta}(q)$. Hence the scalar coefficient in \eqref{max:eq:factorial.kernel.comparison} tends to $\Lambda_{q,x,H}^{k}$ from \eqref{max:eq:band.parameter}. The normalized kernel integral is nonnegative and has upper limit at most one by Lemma~\ref{max:retained.kernel}, so replacing the coefficient by $\Lambda_{q,x,H}^{k}$ introduces an error with upper limit at most
\[
\left|\left\{\frac{p(q,R)}{q}\Theta((x,H])\right\}^{k} - \Lambda_{q,x,H}^{k}\right| \to 0, \qquad R \to \infty.
\]
Consequently, \eqref{max:eq:factorial.kernel.comparison} implies
\[
\limsup_{T\to\infty}\left|\EE_{\star}[(\mathfrak{W}_{T})_{k}] - \frac{\Lambda_{q,x,H}^{k}}{|I_{\theta}|^{k}}\EE_{S}\left[\int_{\mathcal{R}_{k,T}(S)}\mathcal{K}_{T}(\bb{h})\rd{\bb{h}}\right]\right| = 0.
\]
Applying the lower and upper bounds in Lemma~\ref{max:retained.kernel} to this last comparison gives
\[
(1 - 2k\e)_{+}\Lambda_{q,x,H}^{k} \leq \liminf_{T\to\infty}\EE_{\star}[(\mathfrak{W}_{T})_{k}] \leq \limsup_{T\to\infty}\EE_{\star}[(\mathfrak{W}_{T})_{k}] \leq \Lambda_{q,x,H}^{k},
\]
which is \eqref{max:eq:factorial.bounds}.
\end{proof}

The next proposition converts the factorial-moment bounds in Proposition~\ref{max:factorial.moments} into an asymptotic formula for the probability that no random-grid value exceeds $u_{T} + x$. This prepares the passage to the continuum maximum in Proposition~\ref{max:mesh.to.continuum}. Recall that $\mathcal{B}_{T}(S)$ is the collection of nonempty blocks defined in \eqref{max:eq:retained.blocks}, and $\mathscr{G}_{B} = B\cap(U_{B} + q\ell_{T}\Z)$ is the randomly shifted grid in block $B$ defined in \eqref{max:eq:block.grid}.

\begin{proposition}[Probability of no exceedances on random grids]\label{max:grid.void}
Let
\begin{equation}\label{max:eq:grid.maximum}
M_{T,q,\e}^{\star} \leqdef \max_{B\in\mathcal{B}_{T}(S)}\max_{h\in\mathscr{G}_{B}}F_{h}^{\star}.
\end{equation}
For every fixed $q > 0$ and $x \in \R$,
\begin{equation}\label{max:eq:grid.void}
\lim_{\e\downarrow0}\limsup_{T\to\infty} \left|\PP_{\star}\{M_{T,q,\e}^{\star} \leq u_{T} + x\} - \exp\{ - \kappa_{\theta}(q)e^{-\beta_{c}(\theta) x}\}\right| = 0.
\end{equation}
\end{proposition}

\begin{proof}[Proof of Proposition~\ref{max:grid.void}]
Fix $H > \max\{x,0\}$ and let $\mathfrak{W}_{T} \leqdef \mathfrak{W}_{T}(x,H;q,\e)$ be the count in \eqref{max:eq:band.count}. For a nonnegative integer $n$,
\[
\sum_{k=0}^{2m + 1}( - 1)^{k}\binom{n}{k}
\leq \ind_{\{n = 0\}}
\leq \sum_{k=0}^{2m}( - 1)^{k}\binom{n}{k}, \qquad m \in \N_{0},
\]
where $\binom{n}{k} = 0$ for $k > n$ and $\binom{0}{0} = 1$.

Substitute $n = \mathfrak{W}_{T}$, use $\smash{\binom{\mathfrak{W}_{T}}{k}} = (\mathfrak{W}_{T})_{k}/k!$, where $(\mathfrak{W}_{T})_{k}$ denotes the $k$th falling factorial of $\mathfrak{W}_{T}$, and take expectations to obtain the Bonferroni bounds
\begin{equation}\label{max:eq:bonferroni}
\sum_{k=0}^{2m + 1}\frac{( - 1)^{k}}{k!}\EE_{\star}[(\mathfrak{W}_{T})_{k}]
\leq \PP_{\star}(\mathfrak{W}_{T} = 0)
\leq \sum_{k=0}^{2m}\frac{( - 1)^{k}}{k!}\EE_{\star}[(\mathfrak{W}_{T})_{k}],
\end{equation}
with $(\mathfrak{W}_{T})_{0} = 1$.

Applying \eqref{max:eq:factorial.bounds} of Proposition~\ref{max:factorial.moments} to the finitely many terms in \eqref{max:eq:bonferroni}, taking the appropriate upper and lower limits as $T \to \infty$, then letting $\e \downarrow 0$ and finally $m \to \infty$ gives
\begin{equation}\label{max:eq:band.void}
\lim_{\e\downarrow0}\limsup_{T\to\infty}|\PP_{\star}(\mathfrak{W}_{T} = 0) - e^{-\Lambda_{q,x,H}}| = 0.
\end{equation}

The events $\{\mathfrak{W}_{T} = 0\}$ and $\{M_{T,q,\e}^{\star} \leq u_{T} + x\}$ agree unless $M_{T,q,\e}^{\star} > u_{T} + H$. Since $\smash{M_{T,q,\e}^{\star} \leq M_{T}^{\star}}$, the uniform maximum upper-tail estimate $\limsup_{T\to\infty}\PP_{\star}(M_{T}^{\star} > u_{T} + H) \leq Ce^{-\beta_{c}(\theta)H}$ for $H \geq 0$, established in \citet[Lemma~4.7, Eq.~(4.82)]{OuimetExpanding2026}, gives
\begin{equation}\label{max:eq:remove.upper.band}
\limsup_{T\to\infty}\left|\PP_{\star}(\mathfrak{W}_{T} = 0) - \PP_{\star}(M_{T,q,\e}^{\star} \leq u_{T} + x)\right| \leq Ce^{-\beta_{c}(\theta) H},
\end{equation}
uniformly in $q$ and $\e$. For fixed $q > 0$, \eqref{max:eq:band.parameter} of Proposition~\ref{max:factorial.moments} gives $\Lambda_{q,x,H} \to \kappa_{\theta}(q)e^{-\beta_{c}(\theta)x}$ as $H \to \infty$. Combining this convergence with \eqref{max:eq:band.void} and \eqref{max:eq:remove.upper.band}, and taking limits in the order $T \to \infty$, $\e \downarrow 0$, and $H \to \infty$, proves \eqref{max:eq:grid.void}.
\end{proof}

\section{The Euler and Haar continuum maxima}\label{sec:continuum}

\subsection{Approximation by a spatial mesh}

Define the recentered continuum maximum by
\begin{equation}\label{max:recentered.continuum.maximum}
\widetilde{M}_{T}^{\star} = M_{T}^{\star} - u_{T} = \max_{h\in I_{\theta}}F_{h}^{\star} - u_{T}.
\end{equation}
The meshes used below may be deterministic or independent of the Euler or Haar field $F^{\star}$. For $0 < q \leq 1/2$, let $\mathcal{L}_{T,q} \subseteq I_{\theta}$ be a finite set. Suppose that $I_{\theta}$ has a finite cover by intervals $B \subseteq I_{\theta}$, each containing a chosen mesh point $g_{B} \in B \cap \mathcal{L}_{T,q}$. Assume that
\begin{equation}\label{max:mesh.cover}
|B| \leq 2q\ell_{T}, \qquad \#\{B\} \leq \frac{C N_{T}}{q},
\end{equation}
where the first bound holds for every interval $B$ in the cover and $\#\{B\}$ denotes the total number of intervals in that cover. Below, sums and unions indexed by $B$ range over these covering intervals. The constant $C$ is independent of $T$ and $q$. Set the recentered discrete maximum over the mesh $\mathcal{L}_{T,q}$ to be
\begin{equation}\label{max:recentered.mesh.maximum}
\widetilde{M}_{T,q}^{\star} = \max_{h\in\mathcal{L}_{T,q}}F_{h}^{\star} - u_{T}.
\end{equation}
Uniformly spaced meshes with spacing $q\ell_{T}$ satisfy \eqref{max:mesh.cover}. More generally, partition $I_{\theta}$ into intervals whose lengths are bounded below by a fixed positive constant, and place an independently shifted lattice of spacing $q\ell_{T}$ inside each interval, with the shifts also independent of the Euler or Haar field $F^{\star}$. The lower bound on the interval lengths ensures that there are $\OO(|I_{\theta}|)$ partition intervals and that each contains at least two lattice points for all sufficiently large $T$.

Assign to each lattice point its Voronoi cell $B$, clipped to the corresponding partition interval, and take that point as $g_{B}$. This cell consists of the positions at least as close to $g_{B}$ as to any other lattice point in the same partition interval; neighboring cells may share their midpoint boundary. Interior cells have length $q\ell_{T}$, and endpoint cells have length at most $3q\ell_{T}/2$. Moreover, each partition interval contributes at most its length divided by $q\ell_{T}$ plus one lattice point, so the total number of cells satisfies
\[
\#\{B\} \leq \frac{|I_{\theta}|}{q\ell_{T}} + C|I_{\theta}| \leq \frac{C'N_{T}}{q}.
\]
Thus \eqref{max:mesh.cover} holds. The shifts can be conditioned upon throughout the proof of the following lemma because they are independent of $F^{\star}$.

\begin{lemma}[Uniform approximation at a fixed height]\label{max:mesh.error}
For every fixed $x \in \R$ and $\e > 0$, there are $q_{0} > 0$, $c > 0$, and $C_{x} < \infty$ such that
\begin{equation}\label{max:mesh.error.bound}
\limsup_{T\to\infty} \PP_{\star}\big(\widetilde{M}_{T}^{\star} > x + \e, \widetilde{M}_{T,q}^{\star} \leq x\big) \leq C_{x} q^{-1}\exp( - c\e/q), \qquad 0 < q \leq q_{0}.
\end{equation}
The constants are uniform over the meshes satisfying \eqref{max:mesh.cover}, and the same conclusion holds after averaging over independent mesh shifts. The assertion also holds with the continuum maximum taken over a union of intervals in $I_{\theta}$ and the discrete maximum taken over the associated mesh points, provided the union and the collection $(B,g_{B})_{B}$ are deterministic or jointly independent of the Euler or Haar field $F^{\star}$. In this version, the cells need only cover that union and satisfy the two bounds in \eqref{max:mesh.cover}.
\end{lemma}

\begin{proof}[Proof of Lemma~\ref{max:mesh.error}]
The proof is divided into four steps. Step~1 bounds exponential moments of the scaled derivative $\ell_{T}(F^{\star})'$ after conditioning the root height $F_{g_{B}}^{\star}$ to lie in a unit interval. Step~2 uses these bounds to control the oscillation $\omega_{B} = \sup_{h\in B}|F_{h}^{\star} - F_{g_{B}}^{\star}|$ inside a cell $B$. Step~3 sums the resulting bounds over the cells and over height bins covering the nonnegative root heights $0 \leq F_{g_{B}}^{\star} \leq u_{T} + x$. Step~4 handles cells with negative root heights $F_{g_{B}}^{\star} < 0$: an excursion above $u_{T} + x + \e$ in such a cell requires $\omega_{B} > u_{T} + x + \e$, whose probability is controlled by applying Markov's inequality to a sufficiently high moment of $\omega_{B}$. Combining the contributions from nonnegative and negative root heights gives the claimed mesh approximation.

\medskip
{\bf Step 1. Conditional exponential moments of the scaled derivative.} Fix a mesh satisfying \eqref{max:mesh.cover}, and define
\[
Y_{T}^{\star}(z) \leqdef \ell_{T}(F^{\star})'(z), \qquad z \in I_{\theta}.
\]
The absolute-value estimates used to prove \eqref{eq:local.profile.derivative.bounds} in Lemma~\ref{lem:local.profile} are uniform over the spatial positions $h,z \in I_{\theta}$. They therefore give
\begin{equation}\label{max:derivative.variance}
\sup_{z\in I_{\theta}}\Var_{\mathrm{H}}(Y_{T}^{\mathrm{H}}(z)) \leq C, \qquad \sup_{h,z\in I_{\theta}}|\Cov_{\mathrm{H}}(F_{h}^{\mathrm{H}},Y_{T}^{\mathrm{H}}(z))| \leq C.
\end{equation}

Fix $0 < \e_{0} < 1/4$ and recall the field-and-derivative truncation event introduced in the proof of Lemma~\ref{lem:local.profile},
\[
\mathcal{A}_{T,\mathrm{der}}^{\mathrm{A}} =
\left\{\sup_{z\in I_{\theta}}\max\{|F_{z}^{\mathrm{A}}|,|Y_{T}^{\mathrm{A}}(z)|\} \leq t^{1 + \e_{0}}\right\},
\]
with $\mathcal{A}_{T,\mathrm{der}}^{\mathrm{H}}$ the whole Haar probability space. The coefficient bounds have already been verified in \eqref{eq:sup.localization.coefficients}, so the same exceptional-probability estimate used in that proof gives
\begin{equation}\label{max:derivative.exceptional}
\PP_{\mathrm{A}}((\mathcal{A}_{T,\mathrm{der}}^{\mathrm{A}})^{c}) \leq \exp\{ - ct^{1 + 2\e_{0}} + Ct\} = \oo(1).
\end{equation}

Let $v/s_{T}$ range over a fixed compact interval, and define the event that the root height $F_{g_{B}}^{\star}$ lies in $[v - 1,v]$ by
\begin{equation}\label{max:mesh.root.bin.event}
E_{B,v}^{\star} \leqdef \{F_{g_{B}}^{\star} \in v + [ - 1,0]\}.
\end{equation}
The bounds \eqref{max:derivative.variance} and $V_{T} \asymp s_{T}$ from \eqref{eq:field.variance} imply
\[
\frac{|v - 1|}{V_{T}}|\Cov_{\mathrm{H}}(F_{g_{B}}^{\mathrm{H}},Y_{T}^{\mathrm{H}}(z))| \leq C, \qquad 0 \leq \Var_{\mathrm{H}}(Y_{T}^{\mathrm{H}}(z)) - \frac{\Cov_{\mathrm{H}}(F_{g_{B}}^{\mathrm{H}},Y_{T}^{\mathrm{H}}(z))^{2}}{V_{T}} \leq C.
\]
These are the mean and variance expressions in the conditional transform estimate of \citet[Lemma~4.6, Eq.~(4.51)]{OuimetExpanding2026}. Applying that estimate to $\alpha Y_{T}^{\star}(z)$ and $ - \alpha Y_{T}^{\star}(z)$, for a fixed $\alpha > 0$, gives
\begin{equation}\label{max:derivative.conditional}
\EE_{\star}\left[(e^{\alpha Y_{T}^{\star}(z)} + e^{-\alpha Y_{T}^{\star}(z)})\ind_{\mathcal{A}_{T,\mathrm{der}}^{\star}} \,\middle|\, E_{B,v}^{\star}\right] \leq C_{\alpha}.
\end{equation}
The bound is uniform over $g_{B},z \in I_{\theta}$, all the indicated bins $[v - 1,v]$ with $v/s_{T}$ in the fixed compact interval, and the mesh parameter $q$.

\medskip
{\bf Step 2. Conditional oscillation inside a cell.} For a cell $B$ with root $g_{B}$, define
\[
\omega_{B} \leqdef \sup_{h\in B}|F_{h}^{\star} - F_{g_{B}}^{\star}|.
\]
Since the interval between $g_{B}$ and any $h \in B$ lies in $B$, the fundamental theorem of calculus gives
\begin{equation}\label{max:oscillation.integral}
|F_{h}^{\star} - F_{g_{B}}^{\star}| = \left|\ell_{T}^{-1}\int_{g_{B}}^{h}Y_{T}^{\star}(z)\rd{z}\right| \leq \ell_{T}^{-1}\int_{B}|Y_{T}^{\star}(z)|\rd{z}, \qquad h \in B.
\end{equation}
For $|B| > 0$, write $\varrho_{B} \leqdef |B|/(2q\ell_{T}) \in (0,1]$, using \eqref{max:mesh.cover}. Taking the supremum in \eqref{max:oscillation.integral} and applying Jensen's inequality to a convex combination with weights $\varrho_{B}$ and $1 - \varrho_{B}$ yields
\[
\begin{aligned}
\exp\left(\frac{\alpha\omega_{B}}{2q}\right)
&\leq \exp\left\{\varrho_{B}\frac{1}{|B|}\int_{B}\alpha|Y_{T}^{\star}(z)|\rd{z} + (1 - \varrho_{B})0\right\} \\
&\leq 1 - \varrho_{B} + \frac{\varrho_{B}}{|B|}\int_{B}e^{\alpha|Y_{T}^{\star}(z)|}\rd{z}
= 1 - \frac{|B|}{2q\ell_{T}} + \frac{1}{2q\ell_{T}}\int_{B}e^{\alpha|Y_{T}^{\star}(z)|}\rd{z}.
\end{aligned}
\]
Multiplying by $\ind_{\mathcal{A}_{T,\mathrm{der}}^{\star}}$, conditioning on $E_{B,v}^{\star}$ from \eqref{max:mesh.root.bin.event}, and using $e^{\alpha|y|} \leq e^{\alpha y} + e^{-\alpha y}$ and \eqref{max:derivative.conditional} gives
\begin{equation}\label{max:oscillation.exponential.moment}
\EE_{\star}\left[\exp\left(\frac{\alpha\omega_{B}}{2q}\right)\ind_{\mathcal{A}_{T,\mathrm{der}}^{\star}} \,\middle|\, E_{B,v}^{\star}\right] \leq 1 + \frac{C_{\alpha}|B|}{2q\ell_{T}} \leq 1 + C_{\alpha}.
\end{equation}
For a singleton cell, $\omega_{B} = 0$, so \eqref{max:oscillation.exponential.moment} holds directly. Markov's inequality applied to \eqref{max:oscillation.exponential.moment}, followed by an enlargement of $C_{\alpha}$, now yields
\begin{equation}\label{max:oscillation.conditional}
\PP_{\star}\left(\{\omega_{B} > s\}\cap\mathcal{A}_{T,\mathrm{der}}^{\star}\,\middle|\,E_{B,v}^{\star}\right) \leq C_{\alpha}e^{-\alpha s/(2q)}, \qquad s \geq 0.
\end{equation}

\medskip
{\bf Step 3. Sum over nonnegative root heights.} Let $w \leqdef u_{T} + x$. By \eqref{eq:scale.asymptotics} and \eqref{eq:height.asymptotics}, $w \geq ct$ for all sufficiently large $T$. The one-point bin and Haar tail estimates in \citet[proof of Lemma~4.7, Eqs.~(4.86) and~(4.88)]{OuimetExpanding2026}, together with the quantile identity \eqref{eq:height.quantile}, give, for $j \in \{0,\ldots,\lfloor w\rfloor\}$,
\begin{equation}\label{max:root.height.bin}
N_{T}\PP_{\star}(E_{B,w - j}^{\star})
= \frac{\PP_{\star}(E_{B,w - j}^{\star})}{\PP_{\mathrm{H}}(F_{0}^{\mathrm{H}} > u_{T})}
\leq C\frac{V_{T}^{-1/2}\exp\{ - (u_{T} + x - j)^{2}/(2V_{T})\}}{V_{T}^{-1/2}\exp\{ - u_{T}^{2}/(2V_{T})\}}
\leq C_{x}e^{Cj}.
\end{equation}
The bounds are uniform in $j$ because $0 \leq (w - j)/s_{T} \leq w/s_{T} = \OO(1)$. The last inequality uses $u_{T}^{2} - (u_{T} + x - j)^{2} = 2u_{T}(j - x) - (x - j)^{2} \leq 2u_{T}(j + |x|)$ and $u_{T}/V_{T} = \OO(1)$ from \eqref{eq:height.asymptotics}, with the dependence on the fixed $x$ absorbed into $C_{x}$.

The event in \eqref{max:mesh.error.bound} requires the field to be at most $w$ at every mesh point and to exceed $w + \e$ somewhere in at least one cell. If, for such a cell $B$, the value $F_{g_{B}}^{\star}$ at its mesh point lies in $[w - j - 1,w - j]$, so that the event $E_{B,w - j}^{\star}$ in \eqref{max:mesh.root.bin.event} occurs, then
\[
\omega_{B} \geq \sup_{h\in B}F_{h}^{\star} - F_{g_{B}}^{\star} > w + \e - (w - j) = \e + j.
\]
Multiplying the conditional bound \eqref{max:oscillation.conditional}, with $s = \e + j$, by the bin probability \eqref{max:root.height.bin} therefore gives
\[
\begin{aligned}
\PP_{\star}\left(E_{B,w - j}^{\star}\cap\{\omega_{B} > \e + j\}\cap\mathcal{A}_{T,\mathrm{der}}^{\star}\right)
&= \PP_{\star}(E_{B,w - j}^{\star})\PP_{\star}\left(\{\omega_{B} > \e + j\}\cap\mathcal{A}_{T,\mathrm{der}}^{\star}\,\middle|\,E_{B,w - j}^{\star}\right) \\
&\leq \frac{C_{x}}{N_{T}}e^{Cj}e^{-\alpha(\e + j)/(2q)}.
\end{aligned}
\]
The bins $[w - j - 1,w - j]$, for $0 \leq j \leq \lfloor w\rfloor$, cover $[0,w]$. A union bound over these bins and the at most $CN_{T}/q$ cells in \eqref{max:mesh.cover} yields
\begin{equation}\label{max:nonnegative.roots}
\begin{aligned}
&\PP_{\star}\left(\mathcal{A}_{T,\mathrm{der}}^{\star}\cap\bigcup_{B}\left\{0 \leq F_{g_{B}}^{\star} \leq w,\ \sup_{h\in B}F_{h}^{\star} > w + \e\right\}\right) \\
&\qquad\leq \frac{CN_{T}}{q}\frac{C_{x}}{N_{T}}\sum_{j=0}^{\lfloor w\rfloor}e^{Cj}e^{-\alpha(\e + j)/(2q)}
\leq \frac{C_{x}'}{q}e^{-\alpha\e/(2q)}\sum_{j=0}^{\infty}e^{(C - \alpha/(2q))j}.
\end{aligned}
\end{equation}
Take $\alpha = 1$ and choose $q_{0} \leq 1/2$ so that $e^{C - 1/(2q_{0})} \leq 1/2$. For every $0 < q \leq q_{0}$,
\begin{equation}\label{max:mesh.bin.sum}
\sum_{j=0}^{\infty}e^{(C - 1/(2q))j} = \frac{1}{1 - e^{C - 1/(2q)}} \leq 2,
\end{equation}
so that \eqref{max:nonnegative.roots} and \eqref{max:mesh.bin.sum} yield, for every $0 < q \leq q_{0}$,
\begin{equation}\label{max:nonnegative.roots.bound}
\PP_{\star}\left(\mathcal{A}_{T,\mathrm{der}}^{\star}\cap\bigcup_{B}\left\{0 \leq F_{g_{B}}^{\star} \leq w,\ \sup_{h\in B}F_{h}^{\star} > w + \e\right\}\right)
\leq \frac{2C_{x}'}{q}e^{-\e/(2q)}.
\end{equation}

\medskip
{\bf Step 4. Treat negative values at mesh points and combine the bounds.} The moment estimate of \citet[Lemma~4.2, Eq.~(4.9)]{OuimetExpanding2026}, applied to the scaled derivative with the coefficient bound in \citet[proof of Lemma~4.7, Eq.~(4.94)]{OuimetExpanding2026}, gives
\begin{equation}\label{max:derivative.high.moment}
\sup_{z\in I_{\theta}}\EE_{\star}[|Y_{T}^{\star}(z)|^{2D}] \leq (CD)^{D}, \qquad D \in \N, \quad D = \oo(t^{\phi}).
\end{equation}
Taking the supremum in \eqref{max:oscillation.integral} and applying Jensen's inequality to the probability measure $|B|^{-1}\rd{z}$ on a cell with $|B| > 0$ gives
\begin{equation}\label{max:oscillation.high.moment}
\begin{aligned}
\EE_{\star}[\omega_{B}^{2D}]
&\leq \left(\frac{|B|}{\ell_{T}}\right)^{2D}\EE_{\star}\left[\left\{\frac{1}{|B|}\int_{B}|Y_{T}^{\star}(z)|\rd{z}\right\}^{2D}\right] \\
&\leq \left(\frac{|B|}{\ell_{T}}\right)^{2D}\frac{1}{|B|}\int_{B}\EE_{\star}[|Y_{T}^{\star}(z)|^{2D}]\rd{z}
\leq (2q)^{2D}(CD)^{D},
\end{aligned}
\end{equation}
where the last inequality uses \eqref{max:mesh.cover} and \eqref{max:derivative.high.moment}. Singleton cells again satisfy \eqref{max:oscillation.high.moment} because $\omega_{B} = 0$.

If $\smash{F_{g_{B}}^{\star}} < 0$ and $\sup_{h\in B}F_{h}^{\star} > w + \e$, then $\omega_{B} > w + \e$. Choose $D = \lfloor a_{0}t^{2}\rfloor$ with $a_{0} > 0$ sufficiently small. This choice satisfies $D = \oo(t^{\phi})$ because $\phi > 7/2$ in \eqref{eq:fixed.parameters}. Since $w + \e \geq c_{1}t$ eventually and $q \leq 1/2$, choose $a_{0}$ so that
\[
\frac{4q^{2}CD}{(w + \e)^{2}} \leq \frac{Ca_{0}}{c_{1}^{2}} \leq \frac{1}{2}.
\]
Markov's inequality applied to \eqref{max:oscillation.high.moment}, followed by the cell count in \eqref{max:mesh.cover}, gives
\begin{equation}\label{max:negative.roots}
\sum_{B}\PP_{\star}(\omega_{B} > w + \e)
\leq \sum_{B}\frac{\EE_{\star}[\omega_{B}^{2D}]}{(w + \e)^{2D}}
\leq \frac{CN_{T}}{q}\left(\frac{4q^{2}CD}{(w + \e)^{2}}\right)^{D}
\leq \frac{CN_{T}}{q}e^{-ct^{2}} = \oo(1).
\end{equation}
The last equality holds for each fixed $q > 0$ because $\log N_{T} = \OO(t)$ by \eqref{eq:scale.asymptotics}.

If the field is at most $w$ at every mesh point but exceeds $w + \e$ somewhere in a cell $B$, then either $0 \leq F_{g_{B}}^{\star} \leq w$, as treated in Step~3, or $F_{g_{B}}^{\star} < 0$, in which case $\omega_{B} > w + \e$. Consequently, for every $0 < q \leq q_{0}$,
\[
\begin{aligned}
\PP_{\star}\left(\{\widetilde{M}_{T}^{\star} > x + \e\}\cap\{\widetilde{M}_{T,q}^{\star} \leq x\}\right)
&\leq \PP_{\star}((\mathcal{A}_{T,\mathrm{der}}^{\star})^{c}) + \frac{2C_{x}'}{q}e^{-\e/(2q)} + \sum_{B}\PP_{\star}(\omega_{B} > w + \e) \\
&= \frac{2C_{x}'}{q}e^{-\e/(2q)} + \oo(1),
\end{aligned}
\]
using \eqref{max:derivative.exceptional}, \eqref{max:nonnegative.roots.bound}, and \eqref{max:negative.roots}. The exceptional event is common to every cell and is therefore counted only once. Taking the limit superior as $T \to \infty$ proves \eqref{max:mesh.error.bound}.

Let $J_{T} \subseteq I_{\theta}$ be a union of intervals covered by cells $B$ satisfying \eqref{max:mesh.cover}. Since every point of $J_{T}$ belongs to at least one cell,
\[
\left\{\sup_{h\in J_{T}}F_{h}^{\star} > w + \e\right\}\cap\bigcap_{B}\{F_{g_{B}}^{\star} \leq w\}
\subseteq \bigcup_{B}\left\{F_{g_{B}}^{\star} \leq w, \sup_{h\in B}F_{h}^{\star} > w + \e\right\}.
\]
The same bounds therefore prove \eqref{max:mesh.error.bound} with $\widetilde{M}_{T}^{\star}$ replaced by $\sup_{h\in J_{T}}F_{h}^{\star} - u_{T}$ and $\smash{\widetilde{M}_{T,q}^{\star}}$ by $\max_{B}F_{g_{B}}^{\star} - u_{T}$. If $J_{T}$ and the collection $(B,g_{B})_{B}$ are random and jointly independent of $F^{\star}$, condition on their realization. This fixes the intervals and mesh points while leaving the Euler or Haar law unchanged. The probability estimates above are uniform over all such realizations satisfying \eqref{max:mesh.cover}, so averaging over these choices gives the same bound \eqref{max:mesh.error.bound} under the joint law of the field, the intervals, and the mesh.
\end{proof}

\begin{lemma}[Small spatial sets]\label{max:small.spatial.sets}
Let $J_{T} \subseteq I_{\theta}$ be a deterministic union of $k_{T}$ intervals. For every fixed $x \in \R$,
\begin{equation}\label{max:small.spatial.bound}
\limsup_{T\to\infty} \PP_{\star}\left(\sup_{h\in J_{T}}F_{h}^{\star} > u_{T} + x\right) \leq C_{x}\limsup_{T\to\infty} \left(\frac{|J_{T}|}{|I_{\theta}|} + \frac{k_{T}}{N_{T}}\right),
\end{equation}
provided $k_{T} = \OO(N_{T})$ as $T \to \infty$. The same bound holds for sets chosen independently of the Euler or Haar field $F^{\star}$, with the right-hand side replaced by
\[
C_{x}\limsup_{T\to\infty}\EE\left[\frac{|J_{T}|}{|I_{\theta}|} + \frac{k_{T}}{N_{T}}\right].
\]
\end{lemma}

\begin{proof}[Proof of Lemma~\ref{max:small.spatial.sets}]
The argument bounds exceedances in each short interval $B$ by separating the field values $F_{g_{B}}^{\star}$ at the mesh points $g_{B}$ into three ranges: above $u_{T} + x$, between $0$ and $u_{T} + x$, and below $0$. The conditional oscillation bound \eqref{max:oscillation.conditional} controls the second case, while the high-moment bound \eqref{max:oscillation.high.moment} controls the third.

Cover each component of $J_{T}$ by intervals $B$ of length at most $\ell_{T}$, and choose one endpoint $g_{B}$ of each interval as its mesh point. A component of length $a > 0$ requires at most $\lceil a/\ell_{T}\rceil \leq a/\ell_{T} + 1$ cells; a singleton component requires one. Therefore the total number $n_{T}$ of cells satisfies
\begin{equation}\label{max:small.spatial.cells}
n_{T} \leq \frac{|J_{T}|}{\ell_{T}} + k_{T} \leq N_{T} + k_{T} = \OO(N_{T}),
\end{equation}
where $|J_{T}| \leq |I_{\theta}|$ and $N_{T}\ell_{T} = |I_{\theta}|$ by \eqref{eq:scales}. Write
\[
w \leqdef u_{T} + x, \qquad \omega_{B} \leqdef \sup_{h\in B}|F_{h}^{\star} - F_{g_{B}}^{\star}|.
\]
Again, for all sufficiently large $T$, $w \geq ct > 0$ by \eqref{eq:height.asymptotics} and \eqref{eq:scale.asymptotics}.

First, the one-point tail estimate in \citet[proof of Lemma~4.7, Eq.~(4.88)]{OuimetExpanding2026}, normalized by \eqref{eq:height.quantile}, gives uniformly in the chosen endpoint $g_{B}$
\begin{equation}\label{max:small.spatial.high.roots}
\PP_{\star}(F_{g_{B}}^{\star} > w)
= \frac{1}{N_{T}}\frac{\PP_{\star}(F_{g_{B}}^{\star} > u_{T} + x)}{\PP_{\mathrm{H}}(F_{0}^{\mathrm{H}} > u_{T})}
\leq \frac{C}{N_{T}}\exp\left( - \frac{u_{T}x}{V_{T}} - \frac{x^{2}}{2V_{T}}\right)
\leq \frac{C_{x}}{N_{T}}.
\end{equation}
This bounds the contribution of cells $B$ for which the field value $F_{g_{B}}^{\star}$ already exceeds $w$.

Next, suppose $0 \leq F_{g_{B}}^{\star} \leq w$. The intervals $[w - j - 1,w - j]$, for $j \in \{0,\ldots,\lfloor w\rfloor\}$, cover $[0,w]$. On the event $E_{B,w - j}^{\star} = \{F_{g_{B}}^{\star} \in [w - j - 1,w - j]\}$ from \eqref{max:mesh.root.bin.event}, an exceedance inside $B$ implies
\[
\sup_{h\in B}F_{h}^{\star} > w \quad\Longrightarrow\quad \omega_{B} \geq \sup_{h\in B}F_{h}^{\star} - F_{g_{B}}^{\star} > w - (w - j) = j.
\]
The length bound $|B| \leq \ell_{T}$ permits use of \eqref{max:oscillation.conditional} with $q = 1/2$. Combining it with the bin estimate \eqref{max:root.height.bin}, for a fixed $\alpha > 0$, gives
\[
\begin{aligned}
\PP_{\star}\left(E_{B,w - j}^{\star}\cap\{\omega_{B} > j\}\cap\mathcal{A}_{T,\mathrm{der}}^{\star}\right)
&= \PP_{\star}(E_{B,w - j}^{\star})\PP_{\star}\left(\{\omega_{B} > j\}\cap\mathcal{A}_{T,\mathrm{der}}^{\star}\,\middle|\,E_{B,w - j}^{\star}\right) \\
&\leq \frac{C_{x,\alpha}}{N_{T}}e^{Cj}e^{-\alpha j}.
\end{aligned}
\]
Choose $\alpha > C + 1$, with $C$ as in \eqref{max:root.height.bin}. Summing the preceding bound over the bins yields
\begin{equation}\label{max:small.spatial.nonnegative.roots}
\PP_{\star}\left(\left\{\sup_{h\in B}F_{h}^{\star} > w,0 \leq F_{g_{B}}^{\star} \leq w\right\} \cap \mathcal{A}_{T,\mathrm{der}}^{\star}\right)
\leq \frac{C_{x}}{N_{T}}\sum_{j=0}^{\lfloor w\rfloor}e^{(C - \alpha)j}
\leq \frac{C_{x}}{N_{T}}\frac{1}{1 - e^{C - \alpha}}
\leq \frac{C_{x}'}{N_{T}}.
\end{equation}

Finally, if $\smash{F_{g_{B}}^{\star}} < 0$ and $\sup_{h\in B}F_{h}^{\star} > w$, then $\omega_{B} > w$. With $q = 1/2$, the moment bound \eqref{max:oscillation.high.moment} reads $\EE_{\star}[\omega_{B}^{2D}] \leq (CD)^{D}$. Choose $D = \lfloor a_{0}t^{2}\rfloor$ with $a_{0} > 0$ sufficiently small that $CD/w^{2} \leq 1/2$ for all large $T$. Markov's inequality and \eqref{max:small.spatial.cells} then give
\begin{equation}\label{max:small.spatial.negative.roots}
\sum_{B}\PP_{\star}\left(F_{g_{B}}^{\star} < 0,\sup_{h\in B}F_{h}^{\star} > w\right)
\leq \sum_{B}\PP_{\star}(\omega_{B} > w)
\leq n_{T}\left(\frac{CD}{w^{2}}\right)^{D}
\leq Cn_{T}e^{-ct^{2}} = \oo(1).
\end{equation}
Here $D = \oo(t^{\phi})$ by \eqref{eq:fixed.parameters}, as required for \eqref{max:oscillation.high.moment}, and $\log N_{T} = \OO(t)$ by \eqref{eq:scale.asymptotics}.

An exceedance anywhere in $J_{T}$ belongs to at least one covering cell. Outside the single exceptional event $(\mathcal{A}_{T,\mathrm{der}}^{\star})^{c}$, that cell is covered by one of the three cases above. A union bound using \eqref{max:small.spatial.high.roots}, \eqref{max:small.spatial.nonnegative.roots}, and \eqref{max:small.spatial.negative.roots} therefore gives
\[
\begin{aligned}
\PP_{\star}\left(\sup_{h\in J_{T}}F_{h}^{\star} > u_{T} + x\right)
&\leq \PP_{\star}((\mathcal{A}_{T,\mathrm{der}}^{\star})^{c}) + \frac{C_{x}n_{T}}{N_{T}} + Cn_{T}e^{-ct^{2}} \\
&\leq C_{x}\left(\frac{|J_{T}|}{N_{T}\ell_{T}} + \frac{k_{T}}{N_{T}}\right) + \oo(1),
\end{aligned}
\]
where the last inequality uses \eqref{max:derivative.exceptional} and \eqref{max:small.spatial.cells}. Substituting $N_{T}\ell_{T} = |I_{\theta}|$ from \eqref{eq:scales} and taking the limit superior as $T\to \infty$ proves \eqref{max:small.spatial.bound}.

For an independently chosen $J_{T}$, condition on $J_{T}$, apply the same uniform bounds to its covering cells, and average the resulting inequality.
\end{proof}

\subsection{Existence and positivity of the continuum constant}

\begin{lemma}[Tightness at the quantile level]\label{max:quantile.tightness}
For both $\star \in \{\mathrm{A}, \mathrm{H}\}$, the family $(\widetilde{M}_{T}^{\star})_{T}$ defined in \eqref{max:recentered.continuum.maximum} is tight as $T \to \infty$. Moreover, there is $C < \infty$ such that
\begin{equation}\label{max:upper.tail.input}
\limsup_{T\to\infty}\PP_{\star}(\widetilde{M}_{T}^{\star} > x) \leq Ce^{-\beta_{c}(\theta) x}, \qquad x \geq 0.
\end{equation}
\end{lemma}

\begin{proof}[Proof of Lemma~\ref{max:quantile.tightness}]
The upper-tail estimate \eqref{max:upper.tail.input} for $\widetilde{M}_{T}^{\star} = \sup_{h\in I_{\theta}}F_{h}^{\star} - u_{T}$ is the second bound in \citet[Lemma~4.7, Eq.~(4.82)]{OuimetExpanding2026}. Since $\beta_{c}(\theta) > 0$, it gives
\[
0 \leq \lim_{L\to\infty}\limsup_{T\to\infty}\PP_{\star}(\widetilde{M}_{T}^{\star} > L)
\leq \lim_{L\to\infty}Ce^{-\beta_{c}(\theta)L} = 0.
\]
It remains to prove the corresponding statement for $\PP_{\star}(\widetilde{M}_{T}^{\star} < - L)$. Fix $\beta > \beta_{c}(\theta)$ and define
\begin{equation}\label{max:exponential.mass.definition}
Z_{\beta,T}^{\star} \leqdef \ell_{T}^{-1}\int_{I_{\theta}}e^{\beta(F_{h}^{\star} - u_{T})}\rd{h}.
\end{equation}
The lower bound for this exponential mass in \citet[Corollary~4.15]{OuimetExpanding2026} states that
\begin{equation}\label{max:exponential.mass.lower}
\lim_{\delta\downarrow0}\limsup_{T\to\infty}\PP_{\star}(Z_{\beta,T}^{\star} < \delta) = 0.
\end{equation}
Thus $Z_{\beta,T}^{\star}$ is unlikely to be arbitrarily small. The remaining argument shows that the part of this integral contributed by heights below $u_{T} - L$ becomes small as $L \to \infty$.

For $L > 0$, define that contribution by
\begin{equation}\label{max:truncated.mass.definition}
Z_{\beta,T}^{\star,< - L} \leqdef \ell_{T}^{-1}\int_{I_{\theta}}e^{\beta(F_{h}^{\star} - u_{T})}\ind_{\{F_{h}^{\star} - u_{T}<- L\}}\rd{h}.
\end{equation}
Choose $c_{\beta}$ with $0 < c_{\beta} < \beta - \beta_{c}(\theta)$. The first bound in \citet[Lemma~4.7, Eq.~(4.82)]{OuimetExpanding2026} gives
\begin{equation}\label{max:truncated.mass.bound}
\limsup_{T\to\infty}\EE_{\star}[Z_{\beta,T}^{\star,< - L}] \leq Ce^{-c_{\beta}L}.
\end{equation}
Since $Z_{\beta,T}^{\star} = Z_{\beta,T}^{\star,< - L}$ on $\{\widetilde{M}_{T}^{\star} < - L\}$, Markov's inequality and \eqref{max:truncated.mass.bound} give
\begin{equation}\label{max:maximum.lower.tail.comparison}
\limsup_{T\to\infty}\PP_{\star}(\widetilde{M}_{T}^{\star} < - L) \leq \limsup_{T\to\infty}\PP_{\star}(Z_{\beta,T}^{\star} < \delta) + C\delta^{-1}e^{-c_{\beta}L}, \qquad \delta > 0.
\end{equation}
Letting $L \to \infty$ and then $\delta \downarrow 0$, using \eqref{max:exponential.mass.lower}, proves
\begin{equation}\label{max:maximum.lower.tightness}
\lim_{L\to\infty}\limsup_{T\to\infty}\PP_{\star}(\widetilde{M}_{T}^{\star} < - L) = 0.
\end{equation}
Together with \eqref{max:upper.tail.input}, this proves tightness.
\end{proof}

\begin{proposition}[The continuum maximum]\label{max:mesh.to.continuum}
The constants $\kappa_{\theta}(q)$ in \eqref{eq:lattice.constant} satisfy
\begin{equation}\label{max:continuum.constant}
\kappa_{\theta} \leqdef \lim_{q\downarrow0}\kappa_{\theta}(q) = \lim_{q\downarrow0}\frac{1}{q} \PP\{D_{\theta}(qj) \leq 0\text{ for all }j \in \Z\setminus\{0\}\} \in (0, \infty).
\end{equation}
For each $\star \in \{\mathrm{A}, \mathrm{H}\}$ and $x \in \R$, as $T \to \infty$,
\begin{equation}\label{max:continuum.limit}
\PP_{\star}(\widetilde{M}_{T}^{\star} \leq x) \to \exp\{ - \kappa_{\theta} e^{-\beta_{c}(\theta) x}\}.
\end{equation}
\end{proposition}

\begin{proof}[Proof of Proposition~\ref{max:mesh.to.continuum}]
The proof is divided into four steps. Step~1 bounds the grid constants $\kappa_{\theta}(q)$ uniformly in $q$ by comparing the grid maximum with the full continuum maximum. Step~2 covers each block $B \in \mathcal{B}_{T}(S)$ of length at least $q\ell_{T}$ by short intervals $B'$, each containing an associated grid point $g_{B'} \in \mathscr{G}_{B}$, applies the mesh approximation from Lemma~\ref{max:mesh.error}, and bounds the contribution of the shorter endpoint blocks separately. Step~3 combines this approximation with the probability bound for the omitted gaps to bound the distribution of $\smash{\widetilde{M}_{T}^{\star}}$ between two limiting grid distributions. Step~4 uses those bounds to prove that $\kappa_{\theta}(q)$ has a positive limit as $q \downarrow 0$, and then obtains \eqref{max:continuum.limit}.

Use the blocks in $\mathcal{B}_{T}(S)$ defined by \eqref{max:eq:retained.blocks} and the grids $\mathscr{G}_{B}$ from \eqref{max:eq:block.grid}, with independent offsets $U_{B}$. Denote the maximum over all block grids and the continuum maximum on $\bigcup\mathcal{B}_{T}(S)$, both recentered by $u_{T}$, by
\begin{equation}\label{max:retained.maxima.definition}
\widetilde{M}_{T,q,\e}^{\star} \leqdef M_{T,q,\e}^{\star} - u_{T}, \qquad
\widetilde{M}_{T,\e}^{\star} \leqdef \max_{h\in\bigcup\mathcal{B}_{T}(S)}F_{h}^{\star} - u_{T},
\end{equation}
where $\smash{M_{T,q,\e}^{\star}} = \max_{B\in\mathcal{B}_{T}(S)}\max_{h\in\mathscr{G}_{B}}F_{h}^{\star}$, as defined in \eqref{max:eq:grid.maximum} of Proposition~\ref{max:grid.void}.

Since every grid point belongs to a block in $\mathcal{B}_{T}(S)$ and every such block is contained in $I_{\theta}$, these maxima satisfy
\[
\widetilde{M}_{T,q,\e}^{\star} \leq \widetilde{M}_{T,\e}^{\star} \leq \widetilde{M}_{T}^{\star}.
\]
By \eqref{max:retained.maxima.definition}, the grid probability limit \eqref{max:eq:grid.void} in Proposition~\ref{max:grid.void} can equivalently be written as
\begin{equation}\label{max:retained.grid.limit}
\lim_{\e\downarrow0}\limsup_{T\to\infty} \left|\PP_{\star}(\widetilde{M}_{T,q,\e}^{\star} \leq x) - \exp\{-\kappa_{\theta}(q)e^{-\beta_{c}(\theta) x}\}\right| = 0.
\end{equation}

\medskip
{\bf Step 1. Bound $\kappa_{\theta}(q)$ uniformly in $q$.} The inequality $\widetilde{M}_{T,q,\e}^{\star} \leq \widetilde{M}_{T}^{\star}$ gives
\[
1 - \PP_{\star}(\widetilde{M}_{T,q,\e}^{\star} \leq x)
= \PP_{\star}(\widetilde{M}_{T,q,\e}^{\star} > x)
\leq \PP_{\star}(\widetilde{M}_{T}^{\star} > x).
\]
For fixed $x \geq 0$, apply the upper-tail bound \eqref{max:upper.tail.input} of Lemma~\ref{max:quantile.tightness} on the right and the grid limit \eqref{max:retained.grid.limit} on the left. Taking the limit superior as $T \to \infty$ and then letting $\e \downarrow 0$ yields
\begin{equation}\label{max:kappa.tail.bound}
1 - \exp\{ - \kappa_{\theta}(q)e^{-\beta_{c}(\theta)x}\} \leq Ce^{-\beta_{c}(\theta)x}, \qquad x \geq 0.
\end{equation}
For each fixed $q > 0$, the number $\kappa_{\theta}(q) = p(q)/q$ in \eqref{max:eq:grid.constant} is finite and nonnegative. Divide \eqref{max:kappa.tail.bound} by $y = e^{-\beta_{c}(\theta)x}$ and let $x \to \infty$, so that $y \downarrow 0$. Since
\[
\lim_{y\downarrow0}\frac{1 - e^{-\kappa_{\theta}(q)y}}{y} = \kappa_{\theta}(q),
\]
the constant $C$ in \eqref{max:upper.tail.input} of Lemma~\ref{max:quantile.tightness}, which also appears in \eqref{max:kappa.tail.bound} and does not depend on $q$, satisfies
\begin{equation}\label{max:kappa.upper}
0 \leq \kappa_{\theta}(q) \leq C.
\end{equation}

\medskip
{\bf Step 2. Compare the grid and continuum maxima on the union of the blocks.} Fix $x \in \R$ and $\delta > 0$. Apply \eqref{max:mesh.error.bound} of Lemma~\ref{max:mesh.error} with its parameters $x$ and $\e$ replaced by $x - \delta$ and $\delta$, respectively. Since $(x - \delta) + \delta = x$, this gives, for every mesh satisfying \eqref{max:mesh.cover} and all sufficiently small $q > 0$,
\[
\limsup_{T\to\infty}\PP_{\star}\left(\{\widetilde{M}_{T}^{\star} > x\}\cap\{\widetilde{M}_{T,q}^{\star} \leq x - \delta\}\right)
\leq C_{x - \delta}q^{-1}\exp( - c\delta/q).
\]
Thus the estimate controls the probability that the continuum maximum exceeds $u_{T} + x$ while the mesh maximum is at most $u_{T} + x - \delta$. We will use the version of this estimate for a union of covered intervals, as permitted by Lemma~\ref{max:mesh.error}, after verifying \eqref{max:mesh.cover} below.

Condition on $S$ and the offsets $(U_{B})_{B\in\mathcal{B}_{T}(S)}$, so that the blocks and their grids are fixed. Since these random choices are independent of $F^{\star}$, conditioning on them leaves the Euler or Haar law of $F^{\star}$ unchanged.

In a block $B = [a,b]$ of length at least $q\ell_{T}$, list the grid points as $g_{1} < \cdots < g_{m}$, where $m \geq 1$. Split $B$ at the midpoints between consecutive grid points by setting
\[
c_{0} = a, \qquad c_{m} = b, \qquad c_{j} = \frac{g_{j} + g_{j+1}}{2}, \qquad 1 \leq j < m.
\]
The resulting intervals $[c_{j-1},c_{j}]$, $1 \leq j \leq m$, cover $B$ and overlap only at their common endpoints. Each contains exactly one grid point, namely $g_{j}$. Denote any such interval by $B'$ and its grid point by $g_{B'}$. If $m = 1$, the construction gives $B' = B$ and $g_{B'} = g_{1}$.

If $B$ contains at least two grid points, an interior interval has length $q\ell_{T}$. Each endpoint interval consists of a segment of length at most $q\ell_{T}$ between the block endpoint and the nearest grid point, together with half the distance to the next grid point. If $B$ contains a single grid point, both endpoints of $B$ are within distance $q\ell_{T}$ of it. Thus in every case,
\[
|B'| \leq 2q\ell_{T}, \qquad g_{B'} \in B' \cap \mathscr{G}_{B}.
\]
A block contains at most $|B|/(q\ell_{T}) + 1$ grid points. Summing this bound over the blocks and using $|\mathcal{B}_{T}(S)| \leq C(1 + |I_{\theta}|)$ gives
\[
\#\{B'\}
\leq \frac{|I_{\theta}|}{q\ell_{T}} + C(1 + |I_{\theta}|)
\leq \frac{CN_{T}}{q}.
\]
The last inequality follows from $N_{T}\ell_{T} = |I_{\theta}|$ in \eqref{eq:scales}, $\ell_{T} \to 0$, and $|I_{\theta}| \to \infty$. Hence the intervals $B'$ and their associated grid points $g_{B'}$ satisfy \eqref{max:mesh.cover}, uniformly in $S$, the offsets $(U_{B})_{B\in\mathcal{B}_{T}(S)}$, and $0 < \e < 1/4$.

Every full block $B \in \mathcal{B}_{T}(S)$ has length $1 - 2\e > 1/2$. Therefore, for all sufficiently large $T$, only the at most two partial blocks at the endpoints of $I_{\theta}$ can have length below $q\ell_{T}$. Let $J_{T}$ denote the union of those shorter blocks, and let $k_{T}$ be their number. Then
\[
|J_{T}| \leq 2q\ell_{T}, \qquad k_{T} \leq 2, \qquad
\frac{|J_{T}|}{|I_{\theta}|} + \frac{k_{T}}{N_{T}} \leq \frac{2q + 2}{N_{T}} \to 0.
\]
The random set $J_{T}$ depends only on $S$ and is independent of the field. Applying Lemma~\ref{max:small.spatial.sets} and averaging over $S$ therefore gives
\begin{equation}\label{max:short.block.error}
\PP_{\star}\left(\sup_{h\in J_{T}}F_{h}^{\star} > u_{T} + x\right) = \oo(1).
\end{equation}

Consider the event $\{\smash{\widetilde{M}_{T,\e}^{\star}} > x,\smash{\widetilde{M}_{T,q,\e}^{\star}} \leq x - \delta\}$ for the maxima in \eqref{max:retained.maxima.definition}. Either the continuum exceedance occurs in $J_{T}$, or it occurs in one of the blocks covered by the intervals $B'$. In the second case, $F_{g_{B'}}^{\star} \leq u_{T} + x - \delta$ at every grid point $g_{B'}$, while $\sup_{h\in B'}F_{h}^{\star} > u_{T} + x$ for at least one interval $B'$. The union-of-intervals assertion in Lemma~\ref{max:mesh.error} applies to precisely this event. Combining that bound with \eqref{max:short.block.error} and averaging over $S$ and the offsets $(U_{B})_{B\in\mathcal{B}_{T}(S)}$ yields
\begin{equation}\label{max:retained.discretization}
\limsup_{T\to\infty} \PP_{\star}\big(\widetilde{M}_{T,\e}^{\star} > x, \widetilde{M}_{T,q,\e}^{\star} \leq x - \delta\big) \leq r_{q}(x, \delta), \qquad r_{q}(x, \delta) \leqdef C_{x - \delta}q^{-1}e^{-c\delta/q}.
\end{equation}
All bounds used in this step are uniform in $\e$. For fixed $x$ and $\delta$, the exponential factor in \eqref{max:retained.discretization} also gives $r_{q}(x,\delta) \to 0$ as $q \downarrow 0$.

\medskip
{\bf Step 3. Bound the distribution of the full continuum maximum.} Separating according to whether $\smash{\widetilde{M}_{T,\e}^{\star} \leq x}$ gives the event inclusion
\begin{equation}\label{max:retained.event.inclusion}
\{\widetilde{M}_{T,q,\e}^{\star} \leq x - \delta\}
\subseteq \{\widetilde{M}_{T,\e}^{\star} \leq x\}
\cup \{\widetilde{M}_{T,\e}^{\star} > x,\widetilde{M}_{T,q,\e}^{\star} \leq x - \delta\}.
\end{equation}
Taking probabilities in \eqref{max:retained.event.inclusion} and rearranging gives
\begin{equation}\label{max:retained.probability.bound}
\PP_{\star}(\widetilde{M}_{T,\e}^{\star} \leq x)
\geq \PP_{\star}(\widetilde{M}_{T,q,\e}^{\star} \leq x - \delta) - \PP_{\star}\left(\{\widetilde{M}_{T,\e}^{\star} > x\}\cap\{\widetilde{M}_{T,q,\e}^{\star} \leq x - \delta\}\right).
\end{equation}

To compare $\widetilde{M}_{T,\e}^{\star}$ with $\widetilde{M}_{T}^{\star}$, choose a measurable maximizer $h_{T}^{*}$ of the Euler or Haar field $F^{\star}$ on $I_{\theta}$, for example the leftmost one. This point depends only on the field and is independent of $S$, so \eqref{max:eq:retention.probability} gives
\[
\PP_{S}\left\{h_{T}^{*} \notin \bigcup\mathcal{B}_{T}(S)\right\} = 2\e
\]
with the field fixed. On the complementary event, $\smash{\widetilde{M}_{T,\e}^{\star} = \widetilde{M}_{T}^{\star}}$. Since $\smash{\widetilde{M}_{T,\e}^{\star}} \leq \smash{\widetilde{M}_{T}^{\star}}$, with equality whenever $h_{T}^{*} \in \bigcup\mathcal{B}_{T}(S)$, we obtain
\begin{equation}\label{max:eq:gap.loss}
0 \leq \PP_{\star}(\widetilde{M}_{T,\e}^{\star} \leq x) - \PP_{\star}(\widetilde{M}_{T}^{\star} \leq x)
\leq \EE_{\star}\left[\PP_{S}\left\{h_{T}^{*} \notin \bigcup\mathcal{B}_{T}(S)\right\}\right]
= 2\e.
\end{equation}
Combining \eqref{max:eq:gap.loss} with $\widetilde{M}_{T,q,\e}^{\star} \leq \widetilde{M}_{T}^{\star}$ gives
\begin{equation}\label{max:continuum.grid.comparison}
\PP_{\star}(\widetilde{M}_{T}^{\star} \leq x) \geq \PP_{\star}(\widetilde{M}_{T,\e}^{\star} \leq x) - 2\e, \qquad \PP_{\star}(\widetilde{M}_{T}^{\star} \leq x) \leq \PP_{\star}(\widetilde{M}_{T,q,\e}^{\star} \leq x).
\end{equation}
Inserting \eqref{max:retained.probability.bound} into the first inequality of \eqref{max:continuum.grid.comparison} yields, for each $T$,
\[
\PP_{\star}(\widetilde{M}_{T}^{\star} \leq x)
\geq \PP_{\star}(\widetilde{M}_{T,q,\e}^{\star} \leq x - \delta) - \PP_{\star}(\widetilde{M}_{T,\e}^{\star} > x,\widetilde{M}_{T,q,\e}^{\star} \leq x - \delta) - 2\e.
\]
Take the limit inferior in this inequality and the limit superior in the upper bound of \eqref{max:continuum.grid.comparison}. Use \eqref{max:retained.discretization} to control the error probability, and then let $\e \downarrow 0$ using \eqref{max:retained.grid.limit} at the heights $x - \delta$ and $x$. The resulting bounds are
\begin{equation}\label{max:distribution.sandwich}
\begin{aligned}
\exp\{ - \kappa_{\theta}(q)e^{-\beta_{c}(\theta)(x - \delta)}\} - r_{q}(x, \delta)
&\leq \liminf_{T\to\infty}\PP_{\star}(\widetilde{M}_{T}^{\star} \leq x) \\
&\leq \limsup_{T\to\infty}\PP_{\star}(\widetilde{M}_{T}^{\star} \leq x)
\leq \exp\{ - \kappa_{\theta}(q)e^{-\beta_{c}(\theta)x}\}.
\end{aligned}
\end{equation}
The two middle quantities concern only the full continuum maximum and therefore do not depend on $q$ or $\e$.

\medskip
{\bf Step 4. Prove existence and positivity of $\kappa_{\theta}$.} The uniform bounds \eqref{max:kappa.upper} allow the definitions
\[
K_{-} \leqdef \liminf_{q\downarrow0}\kappa_{\theta}(q), \qquad K_{+} \leqdef \limsup_{q\downarrow0}\kappa_{\theta}(q), \qquad 0 \leq K_{-} \leq K_{+} \leq C.
\]
Choose sequences $q_{n}^{-} \downarrow 0$ and $q_{n}^{+} \downarrow 0$ such that $\kappa_{\theta}(q_{n}^{-}) \to K_{-}$ and $\kappa_{\theta}(q_{n}^{+}) \to K_{+}$. In the lower bound of \eqref{max:distribution.sandwich} at $x = 0$, let $q = q_{n}^{-}$ and then $n \to \infty$. In its upper bound, let $q = q_{n}^{+}$ and then $n \to \infty$. The middle quantities are the same for the two sequences, and $r_{q_{n}^{-}}(0,\delta) \to 0$ by \eqref{max:retained.discretization}. Thus
\[
\exp\{ - K_{-}e^{\beta_{c}(\theta)\delta}\}
\leq \liminf_{T\to\infty}\PP_{\star}(\widetilde{M}_{T}^{\star} \leq 0)
\leq \limsup_{T\to\infty}\PP_{\star}(\widetilde{M}_{T}^{\star} \leq 0)
\leq e^{-K_{+}}.
\]
In particular, $\exp\{ - K_{-}e^{\beta_{c}(\theta)\delta}\} \leq e^{-K_{+}}$. Taking logarithms gives $K_{+} \leq e^{\beta_{c}(\theta)\delta}K_{-}$ for every $\delta > 0$. Letting $\delta \downarrow 0$ then gives $K_{+} \leq K_{-}$. Since the reverse inequality holds by definition,
\[
K_{-} = K_{+}, \qquad \kappa_{\theta} \leqdef \lim_{q\downarrow0}\kappa_{\theta}(q) \in [0,C].
\]

To show that this limit is positive, use the lower-tail tightness \eqref{max:maximum.lower.tightness} to choose $L > 1$ sufficiently large that
\[
\limsup_{T\to\infty}\PP_{\star}(\widetilde{M}_{T}^{\star} \leq - L) < \frac{1}{2}.
\]
Now apply the lower bound of \eqref{max:distribution.sandwich} with $x = - L$ and $\delta = 1$, and let $q \downarrow 0$. Since $\kappa_{\theta}(q) \to \kappa_{\theta}$ and $r_{q}(- L,1) \to 0$ by \eqref{max:retained.discretization}, it follows that
\[
\exp\{ - \kappa_{\theta}e^{\beta_{c}(\theta)(L + 1)}\}
\leq \liminf_{T\to\infty}\PP_{\star}(\widetilde{M}_{T}^{\star} \leq - L)
\leq \frac{1}{2}.
\]
Taking logarithms and rearranging gives
\[
\kappa_{\theta} \geq e^{-\beta_{c}(\theta)(L + 1)}\log2 > 0.
\]
This positivity, the upper bound \eqref{max:kappa.upper}, and the definition $\kappa_{\theta}(q) = p(q)/q$ in \eqref{eq:lattice.constant} prove \eqref{max:continuum.constant}. Finally, letting $q \downarrow 0$ and then $\delta \downarrow 0$ in \eqref{max:distribution.sandwich}, using \eqref{max:continuum.constant} and $r_{q}(x,\delta) \to 0$ from \eqref{max:retained.discretization}, proves \eqref{max:continuum.limit}.
\end{proof}

\section{An explicit recentering and the zeta maximum}\label{sec:recentering}

\subsection{Haar cumulants and explicit quantile recentering}

The following lemma identifies the contribution of the arithmetic factor $\mathfrak{a}$ defined in \eqref{eq:maximum.arithmetic.factor} of Theorem~\ref{thm:maximum} to the asymptotic expansion of the Haar cumulant generating function. The subsequent proposition combines this asymptotic with the one-point Haar tail estimate from \citet{OuimetExpanding2026} to evaluate the quantile $u_{T}$ explicitly.

\begin{lemma}[Haar cumulant asymptotics]\label{lem:haar.cumulants}
The arithmetic factor $\mathfrak{a}$ defined in \eqref{eq:maximum.arithmetic.factor} of Theorem~\ref{thm:maximum} is finite, positive, and continuous on $(0,\infty)$. As $T \to \infty$,
\begin{equation}\label{eq:maximum.mod.gaussian.s}
\log\EE_{\mathrm{H}}[e^{\lambda F_{0}^{\mathrm{H}}}] = \frac{\lambda^{2}s_{T}}{4} + \log\mathfrak{a}(\lambda/2) + \oo(1),
\end{equation}
locally uniformly for real $\lambda > 0$.
\end{lemma}

\begin{proof}[Proof of Lemma~\ref{lem:haar.cumulants}]
The proof is divided into three steps. Step~1 evaluates the moment $\EE_{\mathrm{H}}[|1 - U_{p}p^{-\sig}|^{-2z}]$ for each prime $p$ and shows that the product of these moments, each multiplied by $\smash{(1 - p^{-2\sig})^{z^{2}}}$, converges to the arithmetic factor $\mathfrak{a}(z)$ as $\sig \downarrow 1/2$. Step~2 proves that replacing the cutoff $p^{j} \leq X$ by the full sum over prime powers changes $\log\EE_{\mathrm{H}}[e^{\lambda F_{0}^{\mathrm{H}}}]$ by $\oo(1)$. Step~3 evaluates the full product using the pole of $\zeta$ at $1$.

Let $K$ be a fixed compact subinterval of $(0, \infty)$. All the estimates involving a real moment parameter $z$ below are uniform for $z \in K$.

\medskip
{\bf Step 1. The full prime factor.} Let $U$ be uniform on the complex unit circle. For $0 \leq r < 1$, the binomial series and orthogonality give
\[
(1 - rU)^{-z} = \sum_{j=0}^{\infty}\frac{z^{(j)}}{j!}r^{j}U^{j}, \qquad
\EE_{\mathrm{H}}[U^{j}\overline{U}^{\ell}] = \ind_{\{j=\ell\}}.
\]
Multiplying the series for $(1 - rU)^{-z}$ and $(1 - r\overline{U})^{-z}$ therefore yields
\begin{equation}\label{eq:maximum.prime.moment}
\EE_{\mathrm{H}}\left[|1 - rU|^{-2z}\right] = \sum_{j=0}^{\infty}\sum_{\ell=0}^{\infty}\frac{z^{(j)}z^{(\ell)}}{j!\ell!}r^{j + \ell}\EE_{\mathrm{H}}[U^{j}\overline{U}^{\ell}] = \sum_{j=0}^{\infty}\frac{(z^{(j)})^{2}}{(j!)^{2}}r^{2j}.
\end{equation}
The binomial series converge absolutely and uniformly in the phase of $U$ for each $r < 1$, so termwise multiplication and integration are justified. Taylor expansion of the Gauss hypergeometric function \citep[Section~15.2]{OlverEtAl2010} gives, uniformly for $z \in K$ and $0 \leq r \leq 2^{-1/2}$,
\[
\sum_{j=0}^{\infty}\frac{(z^{(j)})^{2}}{(j!)^{2}}r^{2j} = {}_{2}F_{1}(z,z;1;r^{2}) = 1 + z^{2}r^{2} + \OO_{K}(r^{4}).
\]
The logarithm of the factor $(1 - r^{2})^{z^{2}}$ cancels the term of order $r^{2}$ in the logarithm of the moment \eqref{eq:maximum.prime.moment}. More precisely, Taylor expansion at $r = 0$ gives
\[
z^{2}\log(1 - r^{2}) + \log\{1 + z^{2}r^{2} + \OO_{K}(r^{4})\}
= - z^{2}r^{2} + z^{2}r^{2} + \OO_{K}(r^{4}).
\]
Thus,
\begin{equation}\label{eq:maximum.prime.factor.bound}
\log\left\{(1 - r^{2})^{z^{2}} \sum_{j=0}^{\infty}\frac{(z^{(j)})^{2}}{(j!)^{2}}r^{2j}\right\} = \OO_{K}(r^{4}).
\end{equation}
For $r = p^{-\sig}$, the bound in \eqref{eq:maximum.prime.factor.bound} is at most $C_{K}p^{-2}$. Thus the product
\begin{equation}\label{eq:maximum.arithmetic.factor.sigma}
\mathfrak{a}_{\sig}(z) \leqdef \prod_{p}(1 - p^{-2\sig})^{z^{2}} \sum_{j=0}^{\infty}\frac{(z^{(j)})^{2}}{(j!)^{2}p^{2j\sig}}
\end{equation}
has an absolutely convergent logarithm, locally uniformly in $z$. In particular,
\[
|\log\mathfrak{a}_{\sig}(z)| \leq C_{K}\sum_{p}p^{-2} < \infty, \qquad z \in K,
\]
so $\mathfrak{a}_{\sig}(z)$ is positive and finite. For each fixed prime $p$, the logarithm of its factor in $\mathfrak{a}_{\sig}(z)$ converges to the corresponding logarithm in $\mathfrak{a}(z)$, uniformly for $z \in K$, as $\sig \downarrow 1/2$. The series over $j$ defining these factors converge uniformly for $z \in K$ and $\sig \geq 1/2$, by local uniform convergence of the Gauss hypergeometric series and $p^{-2\sig} \leq 1/2$.

For a fixed prime cutoff $P$, split the logarithmic products into $p \leq P$ and $p > P$. The difference over $p \leq P$ tends uniformly to zero, while \eqref{eq:maximum.prime.factor.bound} bounds the difference over $p > P$ by $2C_{K}\sum_{p>P}p^{-2}$. First letting $\sig \downarrow 1/2$ and then $P \to \infty$ proves
\begin{equation}\label{eq:maximum.arithmetic.factor.limit}
\sup_{z\in K}|\log\mathfrak{a}_{\sig}(z) - \log\mathfrak{a}(z)| \to 0.
\end{equation}
At $\sig = 1/2$, the same summable bound \eqref{eq:maximum.prime.factor.bound} shows that $z\mapsto \log\mathfrak{a}(z)$ is a sum of continuous real functions, converging uniformly on every compact subinterval of $(0,\infty)$. It is therefore finite and continuous on $(0,\infty)$, and exponentiation proves that $\mathfrak{a}(z)$ is finite, positive, and continuous there.

\medskip
{\bf Step 2. Removing the prime-power cutoff.} For each prime, define the centered real variables
\begin{equation}\label{eq:maximum.prime.variables}
B_{p,T} \leqdef \Re\sum_{j:p^{j}\leq X}\frac{U_{p}^{j}}{jp^{j\sig}}, \qquad B_{p, \infty} \leqdef \Re\sum_{j=1}^{\infty}\frac{U_{p}^{j}}{jp^{j\sig}} = - \log|1 - U_{p}p^{-\sig}|.
\end{equation}
An empty sum defining $B_{p,T}$ is zero. For either variable $B \in \{B_{p,T},B_{p,\infty}\}$ in \eqref{eq:maximum.prime.variables}, orthogonality gives $\EE_{\mathrm{H}}[B] = 0$, and $|B| \leq C\smash{p^{-1/2}}$. For $\lambda$ in a fixed compact real interval, the remainder in $e^{\lambda B} = 1 + \lambda B + \lambda^{2}B^{2}/2 + \OO(|B|^{3})$ is uniform. Taking expectations therefore gives
\[
\EE_{\mathrm{H}}[e^{\lambda B}]
= 1 + \frac{\lambda^{2}}{2}\EE_{\mathrm{H}}[B^{2}] + \OO(p^{-3/2}), \qquad
\EE_{\mathrm{H}}[B^{2}] = \OO(p^{-1}),
\]
uniformly for real $\lambda$ in that compact interval. For sufficiently large $p$, the quantity $\EE_{\mathrm{H}}[e^{\lambda B}] - 1$ is $\OO(p^{-1})$, so $\log(1 + y) = y + \OO(y^{2})$ adds an error $\OO(p^{-2})$. This error is absorbed by $\OO(p^{-3/2})$, giving
\begin{equation}\label{eq:maximum.prime.cumulant.remainder}
\log\EE_{\mathrm{H}} [e^{\lambda B}] = \frac{\lambda^{2}}{2}\EE_{\mathrm{H}} [B^{2}] + R_{p}(B, \lambda), \qquad |R_{p}(B, \lambda)| \leq Cp^{-3/2},
\end{equation}
for either variable $B \in \{B_{p,T},B_{p,\infty}\}$. For the finitely many small primes, the constant in the bound can be enlarged; their moment generating functions are positive and bounded away from zero on compact real parameter sets.

To compare the quadratic terms in \eqref{eq:maximum.prime.cumulant.remainder}, use
\[
\EE_{\mathrm{H}}[\Re(U_{p}^{j})\Re(U_{p}^{r})] = \frac{1}{2}\ind_{\{j=r\}}, \qquad j,r \in \N.
\]
For each prime $p$, the difference $B_{p,\infty} - B_{p,T}$ consists precisely of the terms $\Re(U_{p}^{j})/(jp^{j\sig})$ with $p^{j} > X$. By the orthogonality identity above, terms with distinct exponents are uncorrelated, and each omitted term has variance $1/(2j^{2}p^{2j\sig})$. The variance of $B_{p,\infty}$ therefore equals the variance of $B_{p,T}$ plus the sum of these omitted variances. Summing over primes gives
\begin{equation}\label{eq:maximum.omitted.variance}
\sum_{p}\left(\EE_{\mathrm{H}} [B_{p, \infty}^{2}] - \EE_{\mathrm{H}} [B_{p,T}^{2}]\right) = \frac{1}{2}\sum_{p,j:p^{j}>X}\frac{1}{j^{2}p^{2j\sig}} = \oo(1),
\end{equation}
where the final estimate is established below. For $j = 1$, the omitted sum is bounded by
\[
\sum_{p>X}p^{-1 - 2\eta_{\sig}} = \sum_{p>X}p^{-1 - \ell_{T}} \leq C\frac{X^{-\ell_{T}}}{\ell_{T}} = \exp\{ - \ell_{T}\log X + s_{T} + \OO(1)\} = \oo(1),
\]
because $\ell_{T}\log X/t \to \infty$ and $s_{T} \sim t$ by \eqref{eq:fixed.parameters}, \eqref{eq:scales}, and \eqref{eq:scale.asymptotics}. For $j \geq 2$, dominated convergence applies to the summable majorant $p^{-j}/j^{2}$:
\[
\sum_{p}\sum_{\substack{j\geq2\\p^{j}>X}}\frac{1}{j^{2}p^{2j\sig}}
\leq \sum_{p}\sum_{j\geq2}\frac{\ind_{\{p^{j}>X\}}}{j^{2}p^{j}}
\to 0.
\]
This proves \eqref{eq:maximum.omitted.variance}. For a fixed prime $p$, let $J_{p,T} \leqdef \lfloor\log X/\log p\rfloor$. The omitted prime powers also satisfy the uniform bound
\[
|B_{p,\infty} - B_{p,T}|
\leq \sum_{j>J_{p,T}}p^{-j/2}
= \frac{p^{-(J_{p,T} + 1)/2}}{1 - p^{-1/2}}
\to 0.
\]
For fixed $p$, both variables are uniformly bounded in $T$, so the differences between their second moments and between their moment generating functions tend to zero as $T \to \infty$, the latter uniformly for $\lambda$ in the fixed compact interval. Their moment generating functions are bounded below by $1$ by Jensen's inequality, since both variables are centered. Their logarithms therefore converge uniformly as well. By the definition of $R_{p}(B,\lambda)$ in \eqref{eq:maximum.prime.cumulant.remainder},
\[
\sup_{\lambda}|R_{p}(B_{p,T},\lambda) - R_{p}(B_{p,\infty},\lambda)| \to 0, \qquad
\sup_{\lambda}|R_{p}(B_{p,T},\lambda) - R_{p}(B_{p,\infty},\lambda)| \leq Cp^{-3/2},
\]
where both suprema are over that compact interval. Summability of $p^{-3/2}$ permits dominated convergence in the sum over primes. Together with the omitted-variance estimate \eqref{eq:maximum.omitted.variance}, this gives
\begin{equation}\label{eq:maximum.log.moment.comparison}
\begin{aligned}
&\sum_{p}\left|\log\EE_{\mathrm{H}}[e^{\lambda B_{p,T}}] - \log\EE_{\mathrm{H}}[e^{\lambda B_{p,\infty}}]\right| \\
&\qquad\leq \frac{\lambda^{2}}{2}\sum_{p}\left(\EE_{\mathrm{H}}[B_{p,\infty}^{2}] - \EE_{\mathrm{H}}[B_{p,T}^{2}]\right)
+ \sum_{p}|R_{p}(B_{p,T},\lambda) - R_{p}(B_{p,\infty},\lambda)|
= \oo(1).
\end{aligned}
\end{equation}
By the definition of the Haar field in \eqref{eq:Haar.field}, $F_{0}^{\mathrm{H}} = \sum_{p\leq X}B_{p,T}$. The variables $B_{p,T}$ are independent over $p$, so
\[
\log\EE_{\mathrm{H}}[e^{\lambda F_{0}^{\mathrm{H}}}]
= \sum_{p\leq X}\log\EE_{\mathrm{H}}[e^{\lambda B_{p,T}}]
= \sum_{p}\log\EE_{\mathrm{H}}[e^{\lambda B_{p,T}}],
\]
where the last equality uses $B_{p,T} = 0$ for $p > X$. The comparison of logarithmic moments in \eqref{eq:maximum.log.moment.comparison} therefore proves
\begin{equation}\label{eq:maximum.cutoff.moment}
\log\EE_{\mathrm{H}}[e^{\lambda F_{0}^{\mathrm{H}}}] = \sum_{p}\log\EE_{\mathrm{H}} [e^{\lambda B_{p, \infty}}] + \oo(1).
\end{equation}
The series on the right is convergent for each $T$, since $\sig > 1/2$ and its quadratic terms are summable.

\medskip
{\bf Step 3. Evaluating the full product.} For $\lambda = 2z$, the identity $B_{p,\infty} = - \log|1 - U_{p}p^{-\sig}|$ in \eqref{eq:maximum.prime.variables} and the prime moment \eqref{eq:maximum.prime.moment}, applied with $r = p^{-\sig}$, give
\[
\EE_{\mathrm{H}}[e^{2zB_{p,\infty}}]
= \EE_{\mathrm{H}}\left[|1 - U_{p}p^{-\sig}|^{-2z}\right]
= \sum_{j=0}^{\infty}\frac{(z^{(j)})^{2}}{(j!)^{2}p^{2j\sig}}.
\]
Taking logarithms in the product defining $\mathfrak{a}_{\sig}(z)$ in \eqref{eq:maximum.arithmetic.factor.sigma} and using the Euler product for $\zeta(2\sig)$, valid because $2\sig > 1$, yields
\[
\sum_{p}\log\EE_{\mathrm{H}}[e^{2zB_{p,\infty}}]
= - z^{2}\sum_{p}\log(1 - p^{-2\sig}) + \log\mathfrak{a}_{\sig}(z)
= z^{2}\log\zeta(2\sig) + \log\mathfrak{a}_{\sig}(z).
\]
Substituting this identity into \eqref{eq:maximum.cutoff.moment} proves
\begin{equation}\label{eq:maximum.moment.euler.product}
\log\EE_{\mathrm{H}}[e^{2zF_{0}^{\mathrm{H}}}] = z^{2}\log\zeta(2\sig) + \log\mathfrak{a}_{\sig}(z) + \oo(1).
\end{equation}
The Laurent expansion $\zeta(1 + \ell_{T}) = \ell_{T}^{-1} + \OO(1)$ \citep[Eq.~(25.2.4)]{OlverEtAl2010}, together with $\ell_{T} \to 0$ and $s_{T} = \log(1/\ell_{T})$ from \eqref{eq:scales} and \eqref{eq:scale.asymptotics}, gives
\begin{equation}\label{eq:maximum.log.zeta.expansion}
\log\zeta(1 + \ell_{T})
= \log(1/\ell_{T}) + \log\{1 + \OO(\ell_{T})\}
= s_{T} + \oo(1).
\end{equation}
Substitute this expansion and the convergence of $\log\mathfrak{a}_{\sig}(z)$ in \eqref{eq:maximum.arithmetic.factor.limit} into \eqref{eq:maximum.moment.euler.product}. Since $2\sig = 1 + \ell_{T}$ and $z = \lambda/2$, the result is
\[
\log\EE_{\mathrm{H}}[e^{\lambda F_{0}^{\mathrm{H}}}]
= \frac{\lambda^{2}s_{T}}{4} + \log\mathfrak{a}(\lambda/2) + \oo(1),
\]
which is \eqref{eq:maximum.mod.gaussian.s}. Every error is uniform for $z \in K$, and $K$ was an arbitrary compact subinterval of $(0,\infty)$, proving the stated local uniformity in $\lambda > 0$.
\end{proof}

\begin{proposition}[Explicit quantile recentering]\label{prop:maximum.recentering}
Let $\theta > 0$, and recall that $m(\theta) = \sqrt{1 + \theta}$ and $\beta_{c}(\theta) = 2m(\theta)$. The quantile $u_{T}$ defined in \eqref{eq:height.quantile} satisfies, as $T \to \infty$,
\begin{equation}\label{eq:maximum.explicit.recentering}
u_{T} = \sqrt{s_{T}\log N_{T}} - \frac{\log(4\pi\log N_{T}) - 2\log\mathfrak{a}(m(\theta))}{2\beta_{c}(\theta)} + \oo(1).
\end{equation}
Equivalently, since $\log N_{T} = s_{T} + \theta t + \log 2$, the quantile satisfies, as $T \to \infty$,
\begin{equation}\label{eq:maximum.explicit.recentering.second}
u_{T} = \sqrt{s_{T}(s_{T} + \theta t)} - \frac{\log s_{T}}{4m(\theta)} + \frac{1}{2m(\theta)}\log\left(\frac{\mathfrak{a}(m(\theta))}{m(\theta)\sqrt{\pi}}\right) + \oo(1).
\end{equation}
\end{proposition}

\begin{proof}[Proof of Proposition~\ref{prop:maximum.recentering}]
The proof is divided into two steps. Step~1 evaluates the constant difference between $V_{T}$ and $s_{T}/2$, then uses the cumulant asymptotic in Lemma~\ref{lem:haar.cumulants} to rewrite the one-point tail probability in \eqref{eq:maximum.imported.scalar.tail} in terms of $s_{T}$ and $\mathfrak{a}$. Step~2 substitutes $v = u_{T}$ into that tail formula and solves the defining equation $N_{T}\PP_{\mathrm{H}}(F_{0}^{\mathrm{H}} > u_{T}) = 1$ from \eqref{eq:height.quantile}, keeping all terms that contribute a nonvanishing correction to $u_{T}$.

\medskip
{\bf Step 1. The scalar tail in the $s_{T}$ normalization.} The scalar tail estimate established in \citet[proof of Lemma~4.7, Eqs.~(4.85) and~(4.88)]{OuimetExpanding2026} is
\begin{equation}\label{eq:maximum.imported.scalar.tail}
\PP_{\mathrm{H}}(F_{0}^{\mathrm{H}} > v) = \frac{\mathfrak{F}(v/V_{T})}{(v/V_{T})\sqrt{2\pi V_{T}}} \exp\left( - \frac{v^{2}}{2V_{T}}\right)\{1 + \oo(1)\},
\end{equation}
uniformly when $v/s_{T}$ ranges over a fixed compact subinterval of $(0, \infty)$. The factor in this formula is
\begin{equation}\label{eq:maximum.mod.factor.definition}
\mathfrak{F}(\lambda) \leqdef \lim_{T\to\infty}e^{-\lambda^{2}V_{T}/2}\EE_{\mathrm{H}}[e^{\lambda F_{0}^{\mathrm{H}}}], \qquad \lambda > 0.
\end{equation}
Recall from \eqref{eq:field.variance} that $V_{T} = \smash{\frac{1}{2}\sum_{p^{j}\leq X}j^{-2}p^{-2j\sig}}$. To rewrite \eqref{eq:maximum.imported.scalar.tail} in terms of $s_{T}$, we first determine the constant correction to the leading term $s_{T}/2$:
\begin{equation}\label{eq:maximum.variance.constant}
V_{T} = \frac{s_{T}}{2} + v_{0} + \oo(1), \qquad v_{0} \leqdef \frac{1}{2}\sum_{p}\sum_{j\geq2} \left(\frac{1}{j^{2}} - \frac{1}{j}\right)p^{-j}.
\end{equation}
To prove \eqref{eq:maximum.variance.constant}, the omitted-variance estimate \eqref{eq:maximum.omitted.variance} and $2\sig = 1 + \ell_{T}$ from \eqref{eq:scales} give
\[
2V_{T} = \sum_{p}\sum_{j=1}^{\infty}\frac{p^{-j(1 + \ell_{T})}}{j^{2}} + \oo(1).
\]
Compare this sum with the logarithm of the Euler product,
\[
\log\zeta(1 + \ell_{T}) = \sum_{p}\sum_{j=1}^{\infty}\frac{p^{-j(1 + \ell_{T})}}{j},
\]
which has coefficient $1/j$ in place of $1/j^{2}$. Subtracting cancels every term with $j = 1$, and dominated convergence with the summable majorant $p^{-j}$ for $j \geq 2$ gives
\[
2V_{T} - \log\zeta(1 + \ell_{T})
= \sum_{p}\sum_{j\geq2}\left(\frac{1}{j^{2}} - \frac{1}{j}\right)p^{-j(1 + \ell_{T})} + \oo(1)
\to \sum_{p}\sum_{j\geq2}\left(\frac{1}{j^{2}} - \frac{1}{j}\right)p^{-j} = 2v_{0}.
\]
Together with the expansion \eqref{eq:maximum.log.zeta.expansion} from the proof of Lemma~\ref{lem:haar.cumulants}, this proves \eqref{eq:maximum.variance.constant}.

Subtracting $\lambda^{2}V_{T}/2$ from the cumulant expansion \eqref{eq:maximum.mod.gaussian.s} in Lemma~\ref{lem:haar.cumulants} and using \eqref{eq:maximum.variance.constant} gives
\[
\log\EE_{\mathrm{H}}[e^{\lambda F_{0}^{\mathrm{H}}}] - \frac{\lambda^{2}V_{T}}{2}
= \log\mathfrak{a}(\lambda/2) - \frac{\lambda^{2}v_{0}}{2} + \oo(1).
\]
Thus the factor $\mathfrak{F}(\lambda)$ defined in \eqref{eq:maximum.mod.factor.definition} is
\begin{equation}\label{eq:maximum.mod.factor.identity}
\mathfrak{F}(\lambda) = \mathfrak{a}(\lambda/2)\exp( - \lambda^{2}v_{0}/2).
\end{equation}

Let $z = v/s_{T}$ range over a fixed compact subinterval of $(0,\infty)$. The variance expansion \eqref{eq:maximum.variance.constant} first gives
\[
\frac{1}{2V_{T}}
= \frac{1}{s_{T}}\left(1 + \frac{2v_{0} + \oo(1)}{s_{T}}\right)^{-1}
= \frac{1}{s_{T}} - \frac{2v_{0}}{s_{T}^{2}} + \oo(s_{T}^{-2}).
\]
Multiplying by $2v$ and $-v^{2}$, respectively, gives the first two expansions below. For the third, use $\mathfrak{F}(\lambda) = \mathfrak{a}(\lambda/2)\exp( - \lambda^{2}v_{0}/2)$ from \eqref{eq:maximum.mod.factor.identity} and continuity of $\log\mathfrak{a}$ on compact positive intervals. The resulting expansions are uniform in the specified range of $z$:
\begin{equation}\label{eq:maximum.tail.normalization}
\begin{aligned}
\frac{v}{V_{T}} &= 2z + \OO(s_{T}^{-1}), \qquad
- \frac{v^{2}}{2V_{T}} = - \frac{v^{2}}{s_{T}} + 2z^{2}v_{0} + \oo(1), \\
\log\mathfrak{F}(v/V_{T}) &= \log\mathfrak{a}(z) - 2z^{2}v_{0} + \oo(1), \qquad
\frac{v}{V_{T}}\sqrt{2\pi V_{T}} = 2z\sqrt{\pi s_{T}}\{1 + \oo(1)\}.
\end{aligned}
\end{equation}
In particular, the two terms involving $v_{0}$ cancel in the logarithm of the numerator of \eqref{eq:maximum.imported.scalar.tail}:
\[
\log\mathfrak{F}(v/V_{T}) - \frac{v^{2}}{2V_{T}}
= \log\mathfrak{a}(z) - 2z^{2}v_{0} - \frac{v^{2}}{s_{T}} + 2z^{2}v_{0} + \oo(1)
= \log\mathfrak{a}(z) - \frac{v^{2}}{s_{T}} + \oo(1).
\]
Exponentiating this identity and using the denominator expansion (the fourth expansion) in \eqref{eq:maximum.tail.normalization} rewrites the tail estimate \eqref{eq:maximum.imported.scalar.tail} as
\begin{equation}\label{eq:maximum.scalar.tail.s}
\PP_{\mathrm{H}}(F_{0}^{\mathrm{H}} > v) = \frac{\mathfrak{a}(v/s_{T})}{2(v/s_{T})\sqrt{\pi s_{T}}} \exp\left( - \frac{v^{2}}{s_{T}}\right)\{1 + \oo(1)\}.
\end{equation}

\medskip
{\bf Step 2. Inverting the quantile equation.} By \eqref{eq:height.asymptotics}, $u_{T}/s_{T} \to m(\theta) > 0$, so $u_{T}/s_{T} \in [m(\theta)/2,2m(\theta)]$ for all sufficiently large $T$. This is a fixed compact subinterval of $(0,\infty)$, so \eqref{eq:maximum.scalar.tail.s} applies with $v = u_{T}$. Multiplying that tail formula by $N_{T}$ and using the quantile normalization \eqref{eq:height.quantile} gives
\[
1 = \frac{N_{T}\mathfrak{a}(u_{T}/s_{T})}{2(u_{T}/s_{T})\sqrt{\pi s_{T}}}\exp\left( - \frac{u_{T}^{2}}{s_{T}}\right)\{1 + \oo(1)\}.
\]
All factors are positive for sufficiently large $T$. Taking logarithms, with $\log\{1 + \oo(1)\} = \oo(1)$, yields
\[
0 = \log N_{T} - \frac{u_{T}^{2}}{s_{T}} - \frac{1}{2}\log s_{T}
+ \log\mathfrak{a}(u_{T}/s_{T}) - \log\{2(u_{T}/s_{T})\sqrt{\pi}\} + \oo(1).
\]
The positivity and continuity of $\mathfrak{a}$ from Lemma~\ref{lem:haar.cumulants} therefore imply
\begin{equation}\label{eq:maximum.quantile.square}
\frac{u_{T}^{2}}{s_{T}} = \log N_{T} - \frac{1}{2}\log s_{T} + \log\mathfrak{a}(m(\theta)) - \log(2m(\theta)\sqrt{\pi}) + \oo(1).
\end{equation}
Let $w_{T} \leqdef \sqrt{s_{T}\log N_{T}}$, so that $w_{T}^{2}/s_{T} = \log N_{T}$. Subtracting this identity from \eqref{eq:maximum.quantile.square} gives
\[
\frac{(u_{T} - w_{T})(u_{T} + w_{T})}{s_{T}}
= - \frac{1}{2}\log s_{T} + \log\left(\frac{\mathfrak{a}(m(\theta))}{2m(\theta)\sqrt{\pi}}\right) + \oo(1).
\]
Solving for the difference $u_{T} - w_{T}$ yields
\begin{equation}\label{eq:maximum.quantile.difference}
u_{T} - w_{T}
= \frac{s_{T}}{u_{T} + w_{T}}\left\{ - \frac{1}{2}\log s_{T} + \log\left(\frac{\mathfrak{a}(m(\theta))}{2m(\theta)\sqrt{\pi}}\right) + \oo(1)\right\}
= \OO(\log s_{T}),
\end{equation}
because $u_{T} \asymp w_{T} \asymp s_{T}$ by \eqref{eq:height.asymptotics} and \eqref{eq:scale.asymptotics}. To replace the coefficient $s_{T}/(u_{T} + w_{T})$ in \eqref{eq:maximum.quantile.difference} by its limit, a rate is needed because this coefficient multiplies $\log s_{T}$. By $N_{T} = |I_{\theta}|/\ell_{T}$ and $s_{T} = \log(1/\ell_{T})$ from \eqref{eq:scales},
\[
\log N_{T} = s_{T} + \theta t + \log2, \qquad
\frac{\log N_{T}}{s_{T}} - m(\theta)^{2}
= \frac{\theta(t - s_{T}) + \log2}{s_{T}}
= \OO\left(\frac{\log t}{t}\right),
\]
where the last equality uses $s_{T} = t + \OO(\log t)$ from \eqref{eq:scale.asymptotics}. Taking square roots and using $u_{T} - w_{T} = \OO(\log s_{T})$ from \eqref{eq:maximum.quantile.difference} gives
\[
\frac{w_{T}}{s_{T}} = \sqrt{\frac{\log N_{T}}{s_{T}}}
= m(\theta) + \OO\left(\frac{\log t}{t}\right), \qquad
\frac{u_{T} - w_{T}}{s_{T}} = \OO\left(\frac{\log t}{t}\right).
\]
Thus $(u_{T} + w_{T})/s_{T} = 2m(\theta) + \OO(\log t/t)$. Since $m(\theta) > 0$, taking the reciprocal gives
\begin{equation}\label{eq:maximum.quantile.denominator}
\frac{s_{T}}{u_{T} + w_{T}} = \frac{1}{2m(\theta)} + \OO\left(\frac{\log t}{t}\right).
\end{equation}
Substituting \eqref{eq:maximum.quantile.denominator} into \eqref{eq:maximum.quantile.difference} produces an error $\OO((\log t)^{2}/t) = \oo(1)$ and proves
\begin{equation}\label{eq:maximum.quantile.expanded}
u_{T} = \sqrt{s_{T}\log N_{T}} - \frac{\log s_{T}}{4m(\theta)} + \frac{1}{2m(\theta)}\log\left(\frac{\mathfrak{a}(m(\theta))}{2m(\theta)\sqrt{\pi}}\right) + \oo(1).
\end{equation}
To obtain the first stated form \eqref{eq:maximum.explicit.recentering}, use $\log N_{T}/s_{T} \to m(\theta)^{2}$ from \eqref{eq:scale.asymptotics} to write
\[
\log(4\pi\log N_{T})
= \log s_{T} + \log\left(4\pi\frac{\log N_{T}}{s_{T}}\right)
= \log s_{T} + 2\log(2m(\theta)\sqrt{\pi}) + \oo(1).
\]
Since $\beta_{c}(\theta) = 2m(\theta)$ by \eqref{eq:height.asymptotics}, substituting this identity into \eqref{eq:maximum.quantile.expanded} proves \eqref{eq:maximum.explicit.recentering}.

For the second form \eqref{eq:maximum.explicit.recentering.second}, the identity $\log N_{T} = s_{T} + \theta t + \log2$ gives the exact numerator in the following difference of square roots:
\begin{equation}\label{eq:maximum.recentering.log.two}
\sqrt{s_{T}\log N_{T}} - \sqrt{s_{T}(s_{T} + \theta t)} = \frac{s_{T}\log2}{\sqrt{s_{T}\log N_{T}} + \sqrt{s_{T}(s_{T} + \theta t)}} = \frac{\log2}{2m(\theta)} + \oo(1).
\end{equation}
The last equality in \eqref{eq:maximum.recentering.log.two} follows by dividing the denominator by $s_{T}$: both square roots, divided by $s_{T}$, tend to $m(\theta)$. Inserting \eqref{eq:maximum.recentering.log.two} into \eqref{eq:maximum.quantile.expanded} combines the constant terms as
\[
\frac{\log2}{2m(\theta)} + \frac{1}{2m(\theta)}\log\left(\frac{\mathfrak{a}(m(\theta))}{2m(\theta)\sqrt{\pi}}\right)
= \frac{1}{2m(\theta)}\log\left(\frac{\mathfrak{a}(m(\theta))}{m(\theta)\sqrt{\pi}}\right).
\]
This is precisely the correction in \eqref{eq:maximum.explicit.recentering.second}.
\end{proof}

\subsection{Proof of Theorem~\ref{thm:maximum}}

Proposition~\ref{max:mesh.to.continuum} proves $\kappa_{\theta} \in (0,\infty)$ in \eqref{max:continuum.constant} and the convergence
\[
\PP_{\star}(M_{T}^{\star} - u_{T} \leq x) \to \exp\{ - \kappa_{\theta}e^{-\beta_{c}(\theta)x}\}, \qquad x \in \R,
\]
for the Euler and Haar fields, as stated in \eqref{max:continuum.limit}. To compare $M_{T}^{\mathrm{A}}$ with $M_{T}^{\zeta}$, define
\begin{equation}\label{eq:maximum.field.approximation}
\Delta_{T}^{\mathrm{fld}} \leqdef \sup_{h\in I_{\theta}} \left|\log|\zeta(\sig + \ii\tau + \ii h)| - F_{h}^{\mathrm{A}}\right|.
\end{equation}
By \citet[Proposition~3.3, Eq.~(3.8)]{OuimetExpanding2026}, applied under \eqref{eq:fixed.parameters}, $|M_{T}^{\zeta} - M_{T}^{\mathrm{A}}| \leq \Delta_{T}^{\mathrm{fld}} = \oo_{\PP}(1)$. Thus \eqref{max:continuum.limit} of Proposition~\ref{max:mesh.to.continuum} implies \eqref{eq:maximum.limit} by Slutsky's theorem.

For $\eta_{\sig} = c\,t^{A}/\log T$ satisfying \eqref{eq:fixed.parameters}, write $s_{T} = t - d_{T}$, where $d_{T} \leqdef A\log t + \log(2c) = \OO(\log t)$. Taylor expansion gives
\begin{equation}\label{eq:maximum.displacement.expansion}
\sqrt{s_{T}(s_{T} + \theta t)} = m(\theta)t - \frac{2 + \theta}{2m(\theta)}d_{T} + \OO(d_{T}^{2}/t).
\end{equation}
Since $d_{T}^{2}/t = \oo(1)$ and $\log s_{T} = \log t + \oo(1)$, substituting the square-root expansion into \eqref{eq:maximum.explicit.recentering.second} of Proposition~\ref{prop:maximum.recentering} proves \eqref{eq:maximum.displacement.recentering}.

Define the explicit recentering by
\[
a_{T} \leqdef \sqrt{s_{T}\log N_{T}} - \frac{\log(4\pi\log N_{T}) - 2\log\mathfrak{a}(m(\theta))}{2\beta_{c}(\theta)}.
\]
Since \eqref{eq:maximum.explicit.recentering} of Proposition~\ref{prop:maximum.recentering} gives $a_{T} - u_{T} \to 0$, Slutsky's theorem and \eqref{eq:maximum.limit} yield
\[
\PP\big(\beta_{c}(\theta)(M_{T}^{\zeta} - a_{T}) - \log\kappa_{\theta} \leq x\big) \to \exp( - e^{-x}), \qquad x \in \R.
\]
The same conclusion holds for the Euler and Haar maxima by \eqref{max:continuum.limit} of Proposition~\ref{max:mesh.to.continuum}. \qed

\section*{Statement of AI use}
\addcontentsline{toc}{section}{Statement of AI use}

The present paper builds on \citet{OuimetExpanding2026} and relies extensively on its results and methods. A rough draft of that work, containing all its main results, was written in 2019--2021 during my postdoctoral appointment at the California Institute of Technology. ChatGPT~6 Astra subsequently assisted in correcting and completing some of its mathematical arguments.

For the present paper, ChatGPT 6 Astra assisted in the discovery and writing of the mathematical arguments. I have completely reworked and independently validated all proofs by hand, and I take full responsibility for their accuracy and rigor.

\section*{Funding}
\addcontentsline{toc}{section}{Funding}

F.\ Ouimet is supported by the Fonds de recherche du Québec -- Nature et technologies (FRQNT) through the Relève professorale grant doi:10.69777/377977, and by the Natural Sciences and Engineering Research Council of Canada (NSERC) through Discovery Grant RGPIN-2026-04471 and Discovery Launch Supplement DGECR-2026-00449.

\renewcommand{\bibsection}{\section*{References}\addcontentsline{toc}{section}{References}}

\setlength{\bibsep}{0pt plus 0ex}

\bibliographystyle{authordate1-doi-mr}
\bibliography{bib}

@article{Aidekon2013,
    AUTHOR       = {A{\"i}d{\'e}kon, E.},
    TITLE        = {Convergence in law of the minimum of a branching random walk},
    JOURNAL      = {Ann. Probab.},
    FJOURNAL     = {The Annals of Probability},
    VOLUME       = {41},
    NUMBER       = {3A},
    YEAR         = {2013},
    PAGES        = {1362--1426},
    DOI          = {10.1214/12-AOP750},
    MRNUMBER     = {3098680},
}

@article{arXiv:2104.07403,
    AUTHOR       = {Amzallag, E. and Arguin, Louis-Pierre and Bailey, E. C. and Hui, K. and Rao, R.},
    TITLE        = {Evidence of random matrix corrections for the large deviations of {S}elberg's central limit theorem},
    JOURNAL      = {Exp. Math.},
    FJOURNAL     = {Experimental Mathematics},
    VOLUME       = {33},
    NUMBER       = {1},
    YEAR         = {2024},
    PAGES        = {123--135},
    DOI          = {10.1080/10586458.2021.2011806},
    MRNUMBER     = {4713867},
}

@incollection{arXiv:1601.00582,
    AUTHOR       = {Arguin, Louis-Pierre},
    EDITOR       = {Contucci, P. and Giardin{\`a}, C.},
    TITLE        = {Extrema of log-correlated random variables: principles and examples},
    BOOKTITLE    = {Advances in Disordered Systems, Random Processes and Some Applications},
    YEAR         = {2017},
    PAGES        = {166--204},
    PUBLISHER    = {Cambridge University Press},
    ADDRESS      = {Cambridge},
    ISBN         = {978-1-107-12410-3},
    DOI          = {10.1017/9781316403877.005},
    MRNUMBER     = {3644278},
}

@article{ArguinBailey2023,
    AUTHOR       = {Arguin, Louis-Pierre and Bailey, E. C.},
    TITLE        = {Large deviation estimates of {S}elberg's central limit theorem and applications},
    JOURNAL      = {Int. Math. Res. Not. IMRN},
    FJOURNAL     = {International Mathematics Research Notices. IMRN},
    VOLUME       = {2023},
    NUMBER       = {23},
    YEAR         = {2023},
    PAGES        = {20574--20612},
    DOI          = {10.1093/imrn/rnad176},
    MRNUMBER     = {4675078},
}

@article{ArguinBailey2025,
    AUTHOR       = {Arguin, Louis-Pierre and Bailey, E. C.},
    TITLE        = {Lower bounds for the large deviations of {S}elberg's central limit theorem},
    JOURNAL      = {Mathematika},
    FJOURNAL     = {Mathematika. A Journal of Pure and Applied Mathematics},
    VOLUME       = {71},
    NUMBER       = {1},
    YEAR         = {2025},
    PAGES        = {Paper No. e70002, 40 pp.},
    DOI          = {10.1112/mtk.70002},
    MRNUMBER     = {4838549},
}

@misc{ArguinBaileyRoberts2026,
    AUTHOR       = {Arguin, Louis-Pierre and Bailey, E. C. and Roberts, A.},
    TITLE        = {Conditional upper bounds for large deviations and moments of the {R}iemann zeta function},
    YEAR         = {2026},
    HOWPUBLISHED = {Preprint, arXiv:2604.25579},
    EPRINT       = {2604.25579},
    DOI          = {10.48550/arXiv.2604.25579},
}

@article{MR3594368,
    AUTHOR       = {Arguin, Louis-Pierre and Belius, D. and Bourgade, P.},
    TITLE        = {Maximum of the characteristic polynomial of random unitary matrices},
    JOURNAL      = {Comm. Math. Phys.},
    FJOURNAL     = {Communications in Mathematical Physics},
    VOLUME       = {349},
    NUMBER       = {2},
    YEAR         = {2017},
    PAGES        = {703--751},
    DOI          = {10.1007/s00220-016-2740-6},
    MRNUMBER     = {3594368},
}

@article{ABBRS_2019,
    AUTHOR       = {Arguin, Louis-Pierre and Belius, D. and Bourgade, P. and Radziwi{\l\l}, M. and Soundararajan, K.},
    TITLE        = {Maximum of the {R}iemann zeta function on a short interval of the critical line},
    JOURNAL      = {Comm. Pure Appl. Math.},
    FJOURNAL     = {Communications on Pure and Applied Mathematics},
    VOLUME       = {72},
    NUMBER       = {3},
    YEAR         = {2019},
    PAGES        = {500--535},
    DOI          = {10.1002/cpa.21791},
    MRNUMBER     = {3911893},
}

@article{MR3619786,
    AUTHOR       = {Arguin, Louis-Pierre and Belius, D. and Harper, A. J.},
    TITLE        = {Maxima of a randomized {R}iemann zeta function, and branching random walks},
    JOURNAL      = {Ann. Appl. Probab.},
    FJOURNAL     = {The Annals of Applied Probability},
    VOLUME       = {27},
    NUMBER       = {1},
    YEAR         = {2017},
    PAGES        = {178--215},
    DOI          = {10.1214/16-AAP1201},
    MRNUMBER     = {3619786},
}

@misc{ABR-2020,
    AUTHOR       = {Arguin, Louis-Pierre and Bourgade, P. and Radziwi{\l\l}, M.},
    TITLE        = {The {F}yodorov--{H}iary--{K}eating conjecture. {I}},
    YEAR         = {2020},
    HOWPUBLISHED = {Preprint, arXiv:2007.00988},
    EPRINT       = {2007.00988},
    DOI          = {10.48550/arXiv.2007.00988},
}

@misc{ArguinBourgadeRadziwill2023,
    AUTHOR       = {Arguin, Louis-Pierre and Bourgade, P. and Radziwi{\l\l}, M.},
    TITLE        = {The {F}yodorov--{H}iary--{K}eating conjecture. {II}},
    YEAR         = {2023},
    HOWPUBLISHED = {Preprint, arXiv:2307.00982},
    EPRINT       = {2307.00982},
    DOI          = {10.48550/arXiv.2307.00982},
}

@misc{ArguinCreighton2026,
    AUTHOR       = {Arguin, Louis-Pierre and Creighton, N.},
    TITLE        = {Lower bounds for the large deviations and moments of the {R}iemann zeta function on the critical line},
    YEAR         = {2026},
    HOWPUBLISHED = {Preprint, arXiv:2603.01711},
    EPRINT       = {2603.01711},
    DOI          = {10.48550/arXiv.2603.01711},
}

@article{ADH_2021,
    AUTHOR       = {Arguin, Louis-Pierre and Dubach, G. and Hartung, L.},
    TITLE        = {Maxima of a random model of the {R}iemann zeta function over intervals of varying length},
    JOURNAL      = {Ann. Inst. Henri Poincar{\'e} Probab. Stat.},
    FJOURNAL     = {Annales de l'Institut Henri Poincar{\'e} Probabilit{\'e}s et Statistiques},
    VOLUME       = {60},
    NUMBER       = {1},
    YEAR         = {2024},
    PAGES        = {588--611},
    DOI          = {10.1214/22-AIHP1323},
    MRNUMBER     = {4718391},
}

@article{ArguinHamdan2026,
    AUTHOR       = {Arguin, Louis-Pierre and Hamdan, J.},
    TITLE        = {The {F}yodorov--{H}iary--{K}eating conjecture on mesoscopic intervals},
    JOURNAL      = {Ann. Probab.},
    FJOURNAL     = {The Annals of Probability},
    VOLUME       = {54},
    NUMBER       = {5},
    YEAR         = {2026},
    PAGES        = {2173--2215},
    DOI          = {10.1214/25-AOP1783},
}

@article{AOR2019,
    AUTHOR       = {Arguin, Louis-Pierre and Ouimet, F. and Radziwi{\l\l}, M.},
    TITLE        = {Moments of the {R}iemann zeta function on short intervals of the critical line},
    JOURNAL      = {Ann. Probab.},
    FJOURNAL     = {The Annals of Probability},
    VOLUME       = {49},
    NUMBER       = {6},
    YEAR         = {2021},
    PAGES        = {3106--3141},
    DOI          = {10.1214/21-AOP1524},
    MRNUMBER     = {4348686},
}

@article{MR2678896,
    AUTHOR       = {Bourgade, P.},
    TITLE        = {Mesoscopic fluctuations of the zeta zeros},
    JOURNAL      = {Probab. Theory Related Fields},
    FJOURNAL     = {Probability Theory and Related Fields},
    VOLUME       = {148},
    NUMBER       = {3-4},
    YEAR         = {2010},
    PAGES        = {479--500},
    DOI          = {10.1007/s00440-009-0237-3},
    MRNUMBER     = {2678896},
}

@article{BramsonDingZeitouni2016,
    AUTHOR       = {Bramson, M. and Ding, J. and Zeitouni, O.},
    TITLE        = {Convergence in law of the maximum of the two-dimensional discrete {Gaussian} free field},
    JOURNAL      = {Comm. Pure Appl. Math.},
    FJOURNAL     = {Communications on Pure and Applied Mathematics},
    VOLUME       = {69},
    NUMBER       = {1},
    YEAR         = {2016},
    PAGES        = {62--123},
    DOI          = {10.1002/cpa.21621},
    MRNUMBER     = {3433630},
}

@misc{Chang2024,
    AUTHOR       = {Chang, C. K.},
    TITLE        = {Hybrid statistics of a random model of zeta over intervals of varying length},
    YEAR         = {2024},
    HOWPUBLISHED = {Preprint, arXiv:2404.08575},
    EPRINT       = {2404.08575},
    DOI          = {10.48550/arXiv.2404.08575},
}

@article{MR3848391,
    AUTHOR       = {Chhaibi, R. and Madaule, T. and Najnudel, J.},
    TITLE        = {On the maximum of the {${\rm C}\beta {\rm E}$} field},
    JOURNAL      = {Duke Math. J.},
    FJOURNAL     = {Duke Mathematical Journal},
    VOLUME       = {167},
    NUMBER       = {12},
    YEAR         = {2018},
    PAGES        = {2243--2345},
    DOI          = {10.1215/00127094-2018-0016},
    MRNUMBER     = {3848391},
}

@article{ChiariniCiprianiHazraInterfaces2016,
    AUTHOR       = {Chiarini, A. and Cipriani, A. and Hazra, R. S.},
    TITLE        = {Extremes of some {G}aussian random interfaces},
    JOURNAL      = {J. Stat. Phys.},
    FJOURNAL     = {Journal of Statistical Physics},
    VOLUME       = {165},
    NUMBER       = {3},
    YEAR         = {2016},
    PAGES        = {521--544},
    DOI          = {10.1007/s10955-016-1634-5},
    MRNUMBER     = {3562423},
}

@article{ChiariniCiprianiHazra2016,
    AUTHOR       = {Chiarini, A. and Cipriani, A. and Hazra, R. S.},
    TITLE        = {Extremes of the supercritical {G}aussian free field},
    JOURNAL      = {ALEA Lat. Am. J. Probab. Math. Stat.},
    FJOURNAL     = {ALEA. Latin American Journal of Probability and Mathematical Statistics},
    VOLUME       = {13},
    NUMBER       = {2},
    YEAR         = {2016},
    PAGES        = {711--724},
    DOI          = {10.30757/ALEA.v13-28},
}

@article{DingRoyZeitouni2017,
    AUTHOR       = {Ding, J. and Roy, R. and Zeitouni, O.},
    TITLE        = {Convergence of the centered maximum of log-correlated {Gaussian} fields},
    JOURNAL      = {Ann. Probab.},
    FJOURNAL     = {The Annals of Probability},
    VOLUME       = {45},
    NUMBER       = {6A},
    YEAR         = {2017},
    PAGES        = {3886--3928},
    DOI          = {10.1214/16-AOP1152},
    MRNUMBER     = {3729618},
}

@misc{FelsHartungLouidor2024,
    AUTHOR       = {Fels, M. and Hartung, L. and Louidor, O.},
    TITLE        = {{G}aussian free field on the tree subject to a hard wall {II}: Asymptotics},
    YEAR         = {2024},
    HOWPUBLISHED = {Preprint, arXiv:2409.00422},
    EPRINT       = {2409.00422},
    DOI          = {10.48550/arXiv.2409.00422},
}

@article{FyodorovHiaryKeating2012,
    AUTHOR       = {Fyodorov, Y. V. and Hiary, G. A. and Keating, J. P.},
    TITLE        = {Freezing transition, characteristic polynomials of random matrices, and the {R}iemann zeta function},
    JOURNAL      = {Phys. Rev. Lett.},
    FJOURNAL     = {Physical Review Letters},
    VOLUME       = {108},
    NUMBER       = {17},
    YEAR         = {2012},
    PAGES        = {170601, 4 pp.},
    DOI          = {10.1103/PhysRevLett.108.170601},
}

@article{MR3151088,
    AUTHOR       = {Fyodorov, Y. V. and Keating, J. P.},
    TITLE        = {Freezing transitions and extreme values: random matrix theory, $\zeta(\frac{1}{2} + i t)$ and disordered landscapes},
    JOURNAL      = {Philos. Trans. Roy. Soc. A},
    FJOURNAL     = {Philosophical Transactions of the Royal Society A. Mathematical, Physical and Engineering Sciences},
    VOLUME       = {372},
    NUMBER       = {2007},
    YEAR         = {2014},
    PAGES        = {20120503, 32 pp.},
    DOI          = {10.1098/rsta.2012.0503},
    MRNUMBER     = {3151088},
}

@misc{arXiv:1906.05783,
    AUTHOR       = {Harper, A. J.},
    TITLE        = {On the partition function of the {R}iemann zeta function, and the {F}yodorov--{H}iary--{K}eating conjecture},
    YEAR         = {2019},
    HOWPUBLISHED = {Preprint, arXiv:1906.05783},
    EPRINT       = {1906.05783},
    DOI          = {10.48550/arXiv.1906.05783},
}

@book{Kallenberg2021,
    AUTHOR       = {Kallenberg, O.},
    TITLE        = {Foundations of modern probability},
    EDITION      = {Third},
    SERIES       = {Probability Theory and Stochastic Modelling},
    VOLUME       = {99},
    YEAR         = {2021},
    PAGES        = {xii+946},
    PUBLISHER    = {Springer},
    ADDRESS      = {Cham},
    ISBN         = {978-3-030-61870-4; 978-3-030-61871-1},
    DOI          = {10.1007/978-3-030-61871-1},
}

@article{KeatingSnaith2000,
    AUTHOR       = {Keating, J. P. and Snaith, N. C.},
    TITLE        = {Random matrix theory and {$\zeta(1/2 + it)$}},
    JOURNAL      = {Comm. Math. Phys.},
    FJOURNAL     = {Communications in Mathematical Physics},
    VOLUME       = {214},
    NUMBER       = {1},
    YEAR         = {2000},
    PAGES        = {57--89},
    DOI          = {10.1007/s002200000261},
    MRNUMBER     = {1794265},
}

@article{KistlerSchmidt2015,
    AUTHOR       = {Kistler, N. and Schmidt, M. A.},
    TITLE        = {From {D}errida's random energy model to branching random walks: from 1 to 3},
    JOURNAL      = {Electron. Commun. Probab.},
    FJOURNAL     = {Electronic Communications in Probability},
    VOLUME       = {20},
    YEAR         = {2015},
    PAGES        = {no. 47, 12 pp.},
    DOI          = {10.1214/ECP.v20-4189},
    MRNUMBER     = {3358969},
}

@article{Madaule2015,
    AUTHOR       = {Madaule, T.},
    TITLE        = {Maximum of a log-correlated {Gaussian} field},
    JOURNAL      = {Ann. Inst. Henri Poincar{\'e} Probab. Stat.},
    FJOURNAL     = {Annales de l'Institut Henri Poincar{\'e} Probabilit{\'e}s et Statistiques},
    VOLUME       = {51},
    NUMBER       = {4},
    YEAR         = {2015},
    PAGES        = {1369--1431},
    DOI          = {10.1214/14-AIHP633},
    MRNUMBER     = {3414451},
}

@book{MR2378655,
    AUTHOR       = {Montgomery, H. L. and Vaughan, R. C.},
    TITLE        = {Multiplicative number theory. {I}. {Classical} theory},
    SERIES       = {Cambridge Studies in Advanced Mathematics},
    VOLUME       = {97},
    YEAR         = {2007},
    PAGES        = {xviii+552},
    PUBLISHER    = {Cambridge University Press},
    ADDRESS      = {Cambridge},
    ISBN         = {978-0-521-84903-6},
    DOI          = {10.1017/CBO9780511618314},
    MRNUMBER     = {2378655},
}

@article{MR3851835,
    AUTHOR       = {Najnudel, J.},
    TITLE        = {On the extreme values of the {R}iemann zeta function on random intervals of the critical line},
    JOURNAL      = {Probab. Theory Related Fields},
    FJOURNAL     = {Probability Theory and Related Fields},
    VOLUME       = {172},
    NUMBER       = {1-2},
    YEAR         = {2018},
    PAGES        = {387--452},
    DOI          = {10.1007/s00440-017-0812-y},
    MRNUMBER     = {3851835},
}

@book{OlverEtAl2010,
    EDITOR       = {Olver, F. W. J. and Lozier, D. W. and Boisvert, R. F. and Clark, C. W.},
    TITLE        = {{NIST} handbook of mathematical functions},
    YEAR         = {2010},
    PAGES        = {xv+951},
    PUBLISHER    = {U.S. Department of Commerce, National Institute of Standards and Technology; Cambridge University Press},
    ADDRESS      = {Washington, DC; Cambridge},
    ISBN         = {978-0-521-19225-5},
    MRNUMBER     = {2723248},
}

@misc{OuimetHaar2026,
    AUTHOR       = {Ouimet, F.},
    TITLE        = {Overlap distributions under the {G}ibbs measure of a {H}aar model of $\log|\zeta|$ on and near the critical line},
    YEAR         = {2026},
    HOWPUBLISHED = {Preprint, arXiv:2610.14901},
    EPRINT       = {2610.14901},
    DOI          = {10.48550/arXiv.2610.14901},
}

@misc{OuimetExpanding2026,
    AUTHOR       = {Ouimet, F.},
    TITLE        = {{Poisson--Dirichlet statistics for the log-modulus of the Riemann zeta function on expanding intervals near the critical line}},
    YEAR         = {2026},
    HOWPUBLISHED = {Preprint, arXiv:2610.22821},
    EPRINT       = {2610.22821},
    DOI          = {10.48550/arXiv.2610.22821},
}

@article{MR3848227,
    AUTHOR       = {Paquette, E. and Zeitouni, O.},
    TITLE        = {The maximum of the {CUE} field},
    JOURNAL      = {Int. Math. Res. Not. IMRN},
    FJOURNAL     = {International Mathematics Research Notices. IMRN},
    VOLUME       = {2018},
    NUMBER       = {16},
    YEAR         = {2018},
    PAGES        = {5028--5119},
    DOI          = {10.1093/imrn/rnx033},
    MRNUMBER     = {3848227},
}

@article{PaquetteZeitouni2025,
    AUTHOR       = {Paquette, E. and Zeitouni, O.},
    TITLE        = {The extremal landscape for the {C}{$\beta$}{E} ensemble},
    JOURNAL      = {Forum Math. Sigma},
    FJOURNAL     = {Forum of Mathematics. Sigma},
    VOLUME       = {13},
    YEAR         = {2025},
    PAGES        = {Paper No. e1, 126 pp.},
    DOI          = {10.1017/fms.2024.129},
    MRNUMBER     = {4849370},
}

@article{MR3832861,
    AUTHOR       = {Radziwi{\l\l}, M. and Soundararajan, K.},
    TITLE        = {{Selberg}'s central limit theorem for $\log|\zeta(1/2 + i t)|$},
    JOURNAL      = {Enseign. Math.},
    FJOURNAL     = {L'Enseignement Math\'{e}matique},
    VOLUME       = {63},
    NUMBER       = {1-2},
    YEAR         = {2017},
    PAGES        = {1--19},
    DOI          = {10.4171/LEM/63-1/2-1},
    MRNUMBER     = {3832861},
}

@article{MR0020594,
    AUTHOR       = {Selberg, A.},
    TITLE        = {Contributions to the theory of the {R}iemann zeta-function},
    JOURNAL      = {Arch. Math. Naturvid.},
    FJOURNAL     = {Archiv for Mathematik og Naturvidenskab},
    VOLUME       = {48},
    NUMBER       = {5},
    YEAR         = {1946},
    PAGES        = {89--155},
    MRNUMBER     = {0020594},
}

@inproceedings{MR1220477,
    AUTHOR       = {Selberg, A.},
    EDITOR       = {Bombieri, E. and Perelli, A. and Salerno, S. and Zannier, U.},
    TITLE        = {Old and new conjectures and results about a class of {D}irichlet series},
    BOOKTITLE    = {Proceedings of the {A}malfi {C}onference on {A}nalytic {N}umber {T}heory ({M}aiori, 1989)},
    YEAR         = {1992},
    PAGES        = {367--385},
    PUBLISHER    = {Universit\`a di Salerno},
    ADDRESS      = {Salerno},
    MRNUMBER     = {1220477},
}

\end{document}